\documentclass[bj]{imsart}

\RequirePackage{amsthm,amsmath,amsfonts,amssymb,bm}
\RequirePackage[numbers,square]{natbib}  % AMH: square specified to differentiate citation [1] from equation (1)
\RequirePackage[colorlinks,citecolor=blue,urlcolor=blue]{hyperref}
\RequirePackage{graphicx}
\RequirePackage{multirow}
\RequirePackage{multicol}
\RequirePackage{anyfontsize}
\RequirePackage{booktabs,color}
\RequirePackage{float}

\startlocaldefs
\newtheorem{thm}{Theorem}
\newtheorem{lem}{Lemma}
\newtheorem{coro}{Corollary}

\theoremstyle{definition}

\newtheorem{exa}{Example}

\endlocaldefs

\begin{document}

\begin{frontmatter}
\title{Inference of Random Effects for Linear Mixed-Effects Models with a Fixed Number of Clusters}
%\title{A sample article title with some additional note\thanksref{t1}}
\runtitle{Estimation of linear mixed-effects models}
%\thankstext{T1}{A sample additional note to the title.}

\begin{aug}
\author[A]{\fnms{Chih-Hao} \snm{Chang}\ead[label=e1,mark]{jhow@nuk.edu.tw}},
\author[B]{\fnms{Hsin-Cheng} \snm{Huang}\ead[label=e2]{hchuang@stat.sinica.edu.tw}}
\and
\author[C]{\fnms{Ching-Kang} \snm{Ing}\ead[label=e3]{cking@stat.nthu.edu.tw}}
%%%%%%%%%%%%%%%%%%%%%%%%%%%%%%%%%%%%%%%%%%%%%%
%% Addresses                                %%
%%%%%%%%%%%%%%%%%%%%%%%%%%%%%%%%%%%%%%%%%%%%%%
\address[A]{Institute of Statistics,
National University of Kaohsiung, Kaohsiung, Taiwan.
\printead{e1}}

\address[B]{Institute of Statistical Science,
Academia Sinica, Taipei, Taiwan.
\printead{e2}}
\address[C]{Institute of Statistics,
	National Tsing Hua University, HsinChu, Taiwan.
\printead{e3}}
\end{aug}

\begin{abstract}

We consider a linear mixed-effects model with a clustered structure, where the
parameters are estimated using maximum likelihood (ML) based on possibly
unbalanced data. Inference with this model is typically done based on 
asymptotic theory, assuming that the number of clusters tends to infinity with
the sample size. However, when the number of clusters is fixed, classical
asymptotic theory developed under a divergent number of
clusters is no longer valid and can lead to erroneous conclusions. In this
paper, we establish the asymptotic properties of the ML estimators of 
random-effects parameters under a general setting, which can be applied to
conduct valid statistical inference with fixed numbers of clusters. 
Our asymptotic theorems allow both fixed effects and random effects to
be misspecified, and the dimensions of both effects to go to infinity with the
sample size.
\end{abstract}


\end{frontmatter}

\section{Introduction}

Over the past several decades, linear mixed-effects models have been broadly
applied to clustered data~\cite{Longford1993}, longitudinal 
data~\cite{Laird1982,Verbeke2000}, spatial data~\cite{Mardia1984}, and data
in scientific fields~\cite{Jiang2007,Jiang2017}, particularly due
to their usefulness
in modeling data with clustered structures.
Model parameters are traditionally estimated, for example, via
minimum norm quadratic, maximum likelihood (ML), and 
restricted ML (REML) methods. ML and REML
estimators are compared in Gumedze and Dunne~\cite{Gumedze2011}.

Estimating random-effects variances in mixed-effects models is usually more
challenging than estimating fixed-effects parameters.
Although desired asymptotic properties have been developed for ML and
REML estimators of random-effects variances~\cite{Hartley1967,Harville1977,Miller1977}, these are mainly
obtained under the mathematical device of requiring the number of clusters
(denoted as $m$) to grow to infinity with the sample size (denoted as $N$) and
the numbers of fixed effects and random effects (denoted as $p$ and $q$) to be
fixed.
In fact, most asymptotic results for likelihood ratio tests and model selection
in linear mixed-effects models are established under a similar mathematical
device; see Self and Liang~\cite{Self1987}, Stram and Lee~\cite{Stram1994}, 
Crainiceanu and Ruppert~\cite{Crainiceanu2004}, Pu and Niu~\cite{Pu2006}, 
Fan and Li~\cite{Fan2012}, and Peng and Lu~\cite{Peng2012}.  However, in many practical
situations, we are faced with a small $m$, which does not grow to infinity with
$N$.
As pointed out by McNeish and Stapleton~\cite{McNeish2016a} and Huang~\cite{Huang2018},
data collected in the fields of education or developmental psychology typically
have a small number of clusters, corresponding, for example, to classrooms or
schools.
Unfortunately, to the best of our knowledge,
no theoretical justification has been provided for random-effects
estimators when $m$ is fixed.

As shown by Maas and Hox~\cite{Maas2004}, Bell et al.~\cite{Bell2014}, 
and McNeish and Stapleton~\cite{McNeish2016b},
for a linear mixed-effects model with few clusters, random-effects
variances are not well estimated by either ML or REML.
This is because when $m$ is fixed,
the Fisher information for random-effects variances fails to grow with $N$,
and hence the corresponding ML estimators do not achieve consistency.
A similar difficulty arises in a spatial-regression model of Chang et al.~\cite{Chang2017}
under the fixed domain asymptotics,
in which the spatial covariance parameters cannot be consistently estimated.
A direct impact of this difficulty is that the classical central limit theorem
established under $m\rightarrow\infty$ for the ML (or REML) estimators
\cite{Hartley1967,Harville1977,Miller1977} is no longer valid.
Consequently, statistical inference based on the asymptotic results for
$m\rightarrow\infty$ can be misleading.

In this article, we focus on the ML estimators in linear mixed-effects models
with possibly unbalanced data.
We first develop the asymptotic properties of the ML estimators, without
assuming that fixed- and random-effects models are correctly specified, $p$ and
$q$ are fixed, or $m\rightarrow\infty$. Based on the asymptotic properties of
the ML estimators, we provide, for the first time in the mixed-effects models
literature, the asymptotic valid confidence intervals for random-effects
variances when $m$ is fixed. In addition, we present an example illustrating
that empirical best linear unbiased predictors (BLUPs) of random
effects (which are the BLUPs with the unknown parameters replaced by their ML
estimators) compare favorably to least squares (LS) predictors even when
the ML estimators are not consistent; see Section~\ref{section:correct} for
details. Also note that our asymptotic theorems allow both fixed- and
random-effects models to be misspecified.
Consequently, our results are crucial to facilitate further studies on model selection
for linear mixed-effects models with fixed $m$,
in which investigating the impact of model misspecification is indispensable.

This article is organized as follows. Section~\ref{section:models and loss}
introduces the linear mixed-effects model and the regularity conditions. The
asymptotic results for the ML estimators are given in 
Section~\ref{section:MLE}. Section~\ref{section:simulation} describes simulation studies
that confirm our asymptotic theory,
including a comparison between the conventional confidence intervals and the
proposed ones for random-effects variances. A brief discussion is given in
Section~\ref{section:discussion}. The proofs of all the theoretical results are
deferred to the online supplementary material.

\section{Linear Mixed-Effects Models}
\label{section:models and loss}

Consider a set of observations with $m$ clusters, $\{(\bm{y}_i,\bm{X}_i,\bm{Z}_i)\}_{i=1}^m$,
where $\bm{y}_i=(y_{i,1},\dots,y_{i,n_i})'$ is the response vector,
$\bm{X}_{i}$ and $\bm{Z}_i$ are $n_i\times p$ and $n_i\times q$ design matrices
of $p$~and~$q$ covariates with the $(j,k)$-th entries $x_{i,j,k}$ and
$z_{i,j,k}$, respectively, and $n_i$ is the number of observations in cluster
$i$.
A general linear mixed-effects model can be written as
\begin{align}
\bm{y}_i = \bm{X}_i\bm\beta + \bm{Z}_i\bm{b}_i + \bm{\epsilon}_i;\quad i=1,\dots,m,
\label{data:general}
\end{align}
where $\bm\beta=(\beta_1,\dots,\beta_p)'$ is the $p$-vector of fixed effects,
$\bm{b}_i=(b_{i,1},\dots,b_{i,q})'\sim
N(\bm{0},\mathrm{diag}(\sigma^2_1,\dots,\sigma^2_q))$ is the $q$-vector of
random effects, $\bm{\epsilon}_i\sim N(\bm{0},v^2\bm{I}_{n_i})$, and
$\bm{I}_{n_i}$ is the $n_i$-dimensional identity matrix.
Here $\{\bm{b}_i\}$ and $\{\bm{\epsilon}_i\}$ are mutually independent.
Let $\bm{y}$, $\bm{X}$, $\bm{b}$, and $\bm{\epsilon}$ be obtained by stacking
$\{\bm{y}_i\}$, $\{\bm{X}_i\}$, $\{\bm{b}_i\}$, and $\{\bm{\epsilon}_i\}$.
Also let $\bm{Z}=\mathrm{diag}(\bm{Z}_1,\dots,\bm{Z}_m)$ be the block diagonal matrix
with diagonal blocks $\{\bm{Z}_i\}$ and dimension $N\times(mq)$,
where $N=n_1+\cdots+n_m$ is the total sample size.
Let $\theta_k=\sigma^2_k/v^2$; $k=1,\dots,q$ and $\bm{D}=\mathrm{diag}(\theta_1,\dots,\theta_q)$.
Then we can rewrite \eqref{data:general} as
\begin{align}
\bm{y}=\bm{X}\bm\beta+\bm{Z}\bm{b}+\bm{\epsilon}\sim N(\bm{X}\bm\beta,v^2\bm{H}),
\label{data}
\end{align}
where $\bm{H}=\bm{R}+\bm{I}_N$, $\bm{R}=\mathrm{diag}(\bm{R}_1,\dots,\bm{R}_m)$,
and $\bm{R}_i=\bm{Z}_i\bm{D}\bm{Z}'_i$; $i=1,\dots,m$.

Let $\mathcal{A}\times\mathcal{G}\subset 2^{\{1,\dots,p\}}\times
2^{\{1,\dots,q\}}$ be the set of candidate models
with $\alpha\in\mathcal{A}$ and $\gamma\in\mathcal{G}$ corresponding
to the fixed-effects and random-effects covariates indexed by $\alpha$ and
$\gamma$, respectively.
Then a linear mixed-effects model corresponding to
$(\alpha,\gamma)\in\mathcal{A}\times\mathcal{G}$ can be written as
\begin{equation}
\bm{y}=\bm{X}(\alpha)\bm{\beta}(\alpha)+\bm{Z}(\gamma)\bm{b}(\gamma)+\bm{\epsilon}.
\label{lmm}
\end{equation}
\noindent For $i=1,\dots,m$, let $\bm{Z}_i(\gamma)$ be the sub-matrix of $\bm{Z}_i$
and $\bm{b}_i(\gamma)$ be the sub-vector of $\bm{b}_i$
corresponding to $\gamma$.
Then for $\gamma\in\mathcal{G}$,
\begin{align}
\bm{R}_i(\gamma,\bm\theta(\gamma))
\equiv&~\mathrm{var}(\bm{Z}_i(\gamma)\bm{b}_i(\gamma))
= \sum_{k\in\gamma}\theta_k\bm{z}_{i,k}\bm{z}_{i,k}',\notag
\end{align}
\noindent where $\bm{z}_{i,k}$ is the $k$-th column of $\bm{Z}_i$ and $\bm\theta(\gamma)$
is the parameter vector of $\theta_k$; $k\in\gamma$.
In other words, under $(\alpha,\gamma)\in\mathcal{A}\times\mathcal{G}$,
\begin{align}
\bm{y}\sim N(\bm{X}(\alpha)\bm{\beta}(\alpha),v^2\bm{H}(\gamma,\bm\theta)),
\label{model:linear mixed effects}
\end{align}
\noindent where
\begin{align}
\bm{H}(\gamma,\bm\theta)
=&~ \bm{R}(\gamma,\bm\theta)+\bm{I}_N,
\label{matrix:H}\\
\bm{R}(\gamma,\bm\theta)
=&~ \mathrm{diag}(\bm{R}_1(\gamma,\bm\theta),\dots,\bm{R}_m(\gamma,\bm\theta))
=\sum_{i=1}^m\sum_{k\in\gamma}\theta_k\bm{h}_{i,k}\bm{h}_{i,k}',\notag
\end{align}
\noindent $\bm{h}_{i,k}=(\bm{0}_{n_1}',\dots,\bm{0}_{n_{k-1}}',\bm{z}_{i,k}',\bm{0}_{n_{k+1}}',\dots,\bm{0}_{n_m}')'$,
and $\bm{0}_{n_i}$ is the $n_i$-vector of zeros. Here, for notational
simplicity, we suppress the dependence of $\bm{\theta}$ on $\gamma$.

For $(\alpha,\gamma)\in\mathcal{A}\times\mathcal{G}$,
let $p(\alpha)$ be the dimension of $\alpha$ and let $q(\gamma)$ be the
dimension of $\gamma$.
Assume that the true model of $\bm{y}$ is
\begin{align}
\bm{y}\sim N(\bm\mu_{0},v_0^2\bm{H}_{0}),
\label{data:linear mixed-effects}
\end{align}
\noindent where $\bm\mu_0$ is the underlying mean trend, $v_0^2>0$ is the true
value of $v^2$,
$\bm{H}_0=\bm{R}_0+\bm{I}_N$, $\bm{R}_0=\mathrm{diag}(\bm{Z}_1\bm{D}_0\bm{Z}'_1,\dots,\bm{Z}_m\bm{D}_0\bm{Z}_m')$,
and $\bm{D}_{0}=\mathrm{diag}(\theta_{1,0},\dots,\theta_{q,0})$ for some $\theta_{k,0}
\geq 0$; $k=1,\dots,q$. Similarly, let $v_0^2\bm{D}_0 = \mathrm{diag}(\sigma_{1,0}^2,\dots,\sigma_{q,0}^2)$
with $\sigma_{k,0}^2\geq 0$ being the true values of $\sigma_k^2$, for $k=1,\dots,q$.
We say that a fixed-effects model $\alpha$ is correct if there exists
$\bm{\beta}(\alpha)\in\mathbb{R}^{p(\alpha)}$ such that
$\bm{\mu}_{0}=\bm{X}(\alpha)\bm{\beta}(\alpha)$.
Similarly, a random-effects model $\gamma$ is correct if
$\{k:\theta_{k,0}>0,\,k=1,\dots,q\}\subset\gamma$.
Let $\mathcal{A}_0$ and $\mathcal{G}_0$ denote the sets of all correct
fixed-effects and random-effects models, respectively.
A linear mixed-effects model $(\alpha,\gamma)$ is said to be correct if
$(\alpha,\gamma)\in\mathcal{A}_0\times\mathcal{G}_0$.
We denote the smallest correct model by $(\alpha_0,\gamma_0)$, which satisfies
\begin{align*}
p_0\equiv p(\alpha_0)
=&~ \inf_{\alpha\in\mathcal{A}_0} p(\alpha),\\
q_0\equiv q(\gamma_0)
=&~ \inf_{\gamma\in\mathcal{G}_0} q(\gamma),
\end{align*}
\noindent where $p_0>0$ and $q_0>0$ are assumed fixed.

Given a model $(\alpha,\gamma)\in\mathcal{A}\times\mathcal{G}$, the covariance
parameters consist of $\bm\theta$ and $v^2$.
We estimate these by ML.
We assume that $\bm{X}$ and $\bm{Z}$ are of full column rank.
The ML estimators $\hat{\bm{\theta}}(\alpha,\gamma)$ and
$\hat{v}^2(\alpha,\gamma)$ of $\bm{\theta}$ and $v^2$
based on model $(\alpha,\gamma)\in\mathcal{A}\times\mathcal{G}$
can be obtained by minimizing the negative twice profile log-likelihood
function:
\begin{align}
\begin{split}
-2\log L(\bm\theta,v^2;\alpha,\gamma)
=&~ N\log(2\pi) + N\log(v^2) + \log\det(\bm{H}(\gamma,\bm\theta))\\
&~  +\frac{\bm{y}'\bm{H}^{-1}(\gamma,\bm\theta)\bm{A}(\alpha,\gamma;\bm\theta)\bm{y}}{v^2},\\
\end{split}
\label{fn:likelihood}
\end{align}
\noindent where
\begin{align}
\bm{A}(\alpha,\gamma;\bm\theta)
\equiv &~ \bm{I}_{N}-\bm{M}(\alpha,\gamma;\bm\theta),
\label{fn:A}\\
\bm{M}(\alpha,\gamma;\bm\theta)
\equiv &~ \bm{X}(\alpha)(\bm{X}(\alpha)'\bm{H}^{-1}(\gamma,\bm\theta)
\bm{X}(\alpha))^{-1}\bm{X}(\alpha)'\bm{H}^{-1}(\gamma,\bm\theta).
\label{fn:M}
\end{align}
\noindent Note that $\bm{M}^2(\alpha,\gamma;\bm\theta)=\bm{M}(\alpha,\gamma;\bm\theta)$,
$\bm{M}(\alpha,\gamma;\bm\theta)\bm{X}(\alpha)=\bm{X}(\alpha)$ and
\begin{align*}
\bm{M}(\alpha,\gamma;\bm\theta)'\bm{H}^{-1}(\gamma,\bm\theta)\bm{M}(\alpha,\gamma;\bm\theta)
=&~\bm{H}^{-1}(\gamma,\bm\theta)\bm{M}(\alpha,\gamma;\bm\theta).
%  \bm{A}(\alpha,\gamma;\bm\theta)'\bm{H}^{-1}(\gamma,\bm\theta)\bm{A}(\alpha,\gamma;\bm\theta)
% =&~\bm{H}^{-1}(\gamma,\bm\theta)\bm{A}(\alpha,\gamma;\bm\theta).
\end{align*}
\noindent For model $(\alpha,\gamma)\in\mathcal{A}\times\mathcal{G}$,
the ML estimator of $\bm\beta(\alpha)$ is given by
\begin{align}
\hat{\bm{\beta}}(\alpha,\gamma;\hat{\bm{\theta}})
=(\bm{X}(\alpha)'\bm{H}^{-1}(\gamma,\hat{\bm\theta})\bm{X}(\alpha))^{-1}\bm{X}(\alpha)'\bm{H}^{-1}(\gamma,\hat{\bm\theta})\bm{y},
\label{eq:blue}
\end{align}
\noindent where $\hat{\bm{\theta}}=\hat{\bm{\theta}}(\alpha,\gamma)$ satisfies
\begin{align*}
(\hat{\bm\theta}(\alpha,\gamma),\hat{v}^2(\alpha,\gamma))
=&~ \operatorname*{arg\,min}_{\bm\theta\in[0,\infty)^{q(\gamma)},v^2\in(0,\infty)}\{-2\log L(\bm\theta,v^2;\alpha,\gamma)\}.
\end{align*}
\noindent Then the ML estimator of $\sigma_k^2$ is
\[
\hat{\sigma}_k^2(\alpha,\gamma)=\hat{\theta}_k(\alpha,\gamma)\hat{v}^2(\alpha,\gamma);\quad k\in\gamma,
\]
where $\hat{\theta}_k(\alpha,\gamma)$ is the ML estimator of $\theta_k$ based on model $(\alpha,\gamma)$.

To establish the asymptotic theory for the ML estimators of the parameters in
linear mixed-effects models, we impose regularity conditions on
covariates of fixed effects and random effects.
\begin{enumerate}
	\item[(A0)] Let $n_{\min} = \displaystyle\min_{i=1,\dots,m}n_i$.
	Assume that $p=c_p+o(n_{\min}^\tau)$ and $q=c_q+o(n_{\min}^\tau)$, for some
	constant $\tau\in[0,1/2)$,
	where $c_p>0$ and $c_q>0$.
	\item[(A1)] With $\tau$ given in (A0), there exist
	constants $\xi\in(2\tau,1]$ and $d_{i,j}>0$; $i=1,\dots,m$, $j=1,\dots,p$,
	with $0<\inf\{d_{i,j}\}\leq\sup\{d_{i,j}\}<\infty$
	such that for $i=1,\dots,m$ and $1\leq j,j^*\leq p$,
	\begin{align*}
	\bm{x}_{i,j}'\bm{x}_{i,j^*}
	=&~\left\{
	\begin{array}{ll}
	d_{i,j} n_i^{\xi}+o(n_i^\xi); & \mbox{if }j=j^*,\\
	o(n_i^{\xi-\tau}); & \mbox{if }j\neq j^*\:,
	\end{array}
	\right.
	\end{align*}
	\noindent where $\bm{x}_{i,j}$ is the $j$-th column of $\bm{X}_i$, for
	$i=1,\dots,m$ and $j=1,\dots,p$.
	\item[(A2)] With $\tau$ given in (A0), there exist constants
	$\ell\in(2\tau,1]$ and $c_{i,k}>0$; $i=1,\dots,m$, $k=1,\dots,q$,
	with $0<\inf\{c_{i,k}\}\leq\sup\{c_{i,k}\}<\infty$
	such that for $i=1,\dots,m$ and $1\leq k,k^*\leq q$,
	\begin{align*}
	\bm{z}_{i,k}'\bm{z}_{i,k^*}
	=&~\left\{
	\begin{array}{ll}
	c_{i,k} n_i^{\ell}+o(n_i^\ell); & \mbox{if }k=k^*,\\
	o(n_i^{\ell-\tau}); & \mbox{if }k\neq k^*\:.
	\end{array}
	\right.
	\end{align*}
	\item [(A3)] For $i=1,\dots,m$, $j=1,\dots,p$, and $k=1,\dots,q$,
	\begin{align*}
	\bm{x}_{i,j}'\bm{z}_{i,k}
	=&~ o(n_i^{(\xi+\ell)/2-\tau}),
	\end{align*}
	\noindent where $\tau$, $\xi$, and $\ell$ are given in (A0), (A1), and (A2),
	respectively.
\end{enumerate}

Condition~(A0) allows the numbers of fixed effects and random effects (i.e.,
$p$ and $q$) to go to infinity with $n_{\min}$ at a certain rate.
Conditions (A1)--(A3) impose correlation constraints on $\{\bm{x}_{i,j}\}$
and $\{\bm{z}_{i,k}\}$.
For example, Condition~(A2) implies that the maximum eigenvalue satisfies
$\lambda_{\max}(\bm{Z}_i\bm{D}\bm{Z}_i')=O(n_i^\ell)$, which is similar to an
assumption given in Condition~3 of Fan and Li~\cite{Fan2012}.

\setcounter{equation}{0}
\section{Asymptotic Properties}
\label{section:MLE}

In this section, we investigate the asymptotic properties of the ML estimators
of $v^2$ and $\{\sigma_k^2:k\in\gamma\}$
for any $(\alpha,\gamma)\in\mathcal{A}\times\mathcal{G}$.
We allow $p$~and~$q$ to go to infinity with the sample size~$N$.
In addition, as we allow $m$ to be fixed, we must account for the fact that
$\{\sigma_k^2:k\in\gamma\}$ may not be estimated consistently.

\subsection{Asymptotics under correct specification}
\label{section:correct}

In this subsection, we consider a correct (but possibly overfitted) model
$(\alpha,\gamma)\in\mathcal{A}_0\times\mathcal{G}_0$.
We derive not only the convergence rates for the ML estimators of $v^2$ and
$\{\sigma_k^2:k\in\gamma\}$, but also their asymptotic distributions.

\begin{thm}
	Consider the data generated from \eqref{data} with the true parameters given
	by \eqref{data:linear mixed-effects}. Let
	$(\alpha,\gamma)\in\mathcal{A}_0\times\mathcal{G}_0$ be a correct model
	defined in (\ref{model:linear mixed effects}). Denote
	$\hat{\sigma}_k^2(\alpha,\gamma)$ and $\hat{v}^2(\alpha,\gamma)$ to be the
	ML estimators of $\sigma_k^2$ and $v^2$, respectively.
	Assume that (A0)--(A3) hold. Then
	\begin{align}
	\hat{v}^2(\alpha,\gamma)
	=&~ v_0^2 +O_p\Big(\frac{p+mq}{N}\Big)+O_p(N^{-1/2}),
	\label{appendix:thm:mle:correct:eq1}\\
	\hat{\sigma}_k^2(\alpha,\gamma)
	=&~\left\{
	\begin{array}{ll}
	\displaystyle\frac{1}{m}\sum_{i=1}^m b_{i,k}^2 + O_p\bigg(\frac{1}{m}\sum_{i=1}^mn_i^{-\ell/2}\bigg);  &
	\mbox{if }k\in\gamma\cap\gamma_0, \\
	O_p\big(n_{\max}^{-\ell}\big); & \mbox{if }k\in\gamma\setminus\gamma_0,
	\end{array}
	\right.
	\label{appendix:thm:mle:correct:eq2}
	\end{align}
	\noindent where $\displaystyle n_{\max} = \max_{i=1,\ldots,m}n_i$.
	In addition, if $p+mq=o\big(N^{1/2}\big)$, then
	\begin{align*}
	N^{1/2}\big(\hat{v}^2(\alpha,\gamma) - v_0^2\big) \xrightarrow{d} N\big(0,2v_0^4\big),\quad\mbox{as }N\rightarrow\infty.
	\end{align*}
	\label{appendix:theorem:MLE}
\end{thm}

When $m$ is fixed and $(\alpha,\gamma)\in\mathcal{A}_0\times\mathcal{G}_0$, it
follows from \eqref{appendix:thm:mle:correct:eq2} that
$\hat{\sigma}_k^2(\alpha,\gamma)$ does not converge to $\sigma_{k,0}^2$, for
$k\in\gamma\cap\gamma_0$. This is because the data do not contain enough
information for $\{\sigma_k^2:k\in\gamma\cap\gamma_0\}$.
Nevertheless, $\hat{\sigma}^2_k(\alpha,\gamma)$ converges to
$\sigma_{k,0}^2=0$, for $k\in\gamma\setminus\gamma_0$, at a rate
$n_{\max}^{-\ell}$, which can be faster than $N^{-1/2}$.
On the other hand, when $m\rightarrow\infty$,
by applying the law of large numbers and the  central limit theorem to
$b_{i,k}$; $i=1,\dots,m,\,k\in\gamma_0$, we immediately have the following
corollary.

\begin{coro}
	Under the assumptions of Theorem~\ref{appendix:theorem:MLE},
	$\hat{\sigma}_k^2(\alpha,\gamma)\xrightarrow{p}\sigma_{k,0}^2$ as
	$m\rightarrow\infty$, for $k\in\gamma$.
	If, in addition, $m=o(n_{\min}^\ell)$, then
	\begin{align*}
	m^{1/2}(\hat{\sigma}_k^2(\alpha,\gamma) - \sigma_{k,0}^2)
	\xrightarrow{d}N(0,2\sigma_{k,0}^4);\quad k\in\gamma\cap\gamma_0,\quad\mbox{as }N\rightarrow\infty.
	\end{align*}
	\label{coro:normality sigma k}
\end{coro}

From Corollary~\ref{coro:normality sigma k},
for $k\in\gamma_0$, we obtain
a $100(1-\alpha)\%$ confidence interval of $\sigma^2_{k,0}$\,:
\begin{align}
\bigg(\hat{\sigma}_k^2(\alpha,\gamma)
-\bigg(\frac{2\hat{\sigma}^4_k(\alpha,\gamma)}{m}\bigg)^{1/2}
\zeta_{1-\alpha/2},\,
\hat{\sigma}_k^2(\alpha,\gamma)
-\bigg(\frac{2\hat{\sigma}^4_k(\alpha,\gamma)}{m}\bigg)^{1/2}\zeta_{\alpha/2}\bigg),
\label{CI:standard normal}
\end{align}
\noindent where $\zeta_{a}$ is the $(100a)$-th percentile of the standard
normal distribution.
Although this confidence interval is commonly applied in practice (e.g.,
Maas and Hox~\cite{Maas2004}; McNeish and Stapleton~\cite{McNeish2016b}), it is only valid when $m$ is large, as
detailed in a simulation experiment of Section~\ref{section:fixed m}.
Thanks to Theorem~\ref{appendix:theorem:MLE}, we can derive a $100(1-\alpha)\%$
confidence interval of $\sigma_{k,0}^2$ valid for a fixed $m$.

\begin{thm}
	Under the assumptions of Theorem~\ref{appendix:theorem:MLE}, suppose that
	$m$ is fixed.
	Then for $k\in\gamma\cap\gamma_0$,
	% $m\hat{\sigma}_k^2(\alpha,\gamma)/\sigma_{k,0}^2\xrightarrow{d} \chi^2_m$ as $N\rightarrow\infty$,
	% where $\chi^2_m$ denotes the chi-square distribution with degrees of freedom $m$.
	a $100(1-\alpha)\%$ confidence interval of $\sigma^2_k$ is
	\begin{align}
	\bigg(\frac{m\hat{\sigma}_k^2(\alpha,\gamma)}{\chi^2_{m,1-\alpha/2}},
	\frac{m\hat{\sigma}_k^2(\alpha,\gamma)}{\chi^2_{m,\alpha/2}}\bigg),
	\label{CI:chisquare}
	\end{align}
	\noindent where $\chi^2_{m,a}$ denotes the $(100a)$-th percentile of the
	chi-square distribution on $m$ degrees of freedom.
	\label{theorem:CI}
\end{thm}

Note that the length of the confidence interval of $\sigma_{k,0}^2$
in \eqref{CI:chisquare} does not shrink to $0$
as $N\rightarrow\infty$, which is not surprising due to the fact that
$\hat{\sigma}_k^2(\alpha,\gamma)$ is not a consistent estimator of $\sigma_k^2$
when $m$ is fixed, for $k\in\gamma\cap\gamma_0$.

We close this section by mentioning that although a fixed $m$ hinders us from
consistently estimating $\sigma_k^2$, the empirical BLUPs of random
effects, based on the ML estimator of $\sigma_k^2$, are still asymptotically
more efficient than the LS predictors, as illustrated in the following example.

\begin{exa}
Consider model~\eqref{data} with $p=0$, $q=1$, $n_1=\cdots=n_m=n$ and $m>1$
fixed. Assume that (A2) holds with $c_{1,1}=\cdots = c_{m,1}=1$ and $\ell=1$.
Let $\tilde{\bm{b}}_{i}$ be the LS predictor of $\bm{b}_i$ and
$\hat{\bm{b}}_{i}(\sigma_1^2,v^2)$ be the BLUP of $\bm{b}_i$ given
$(\sigma_1^2,v^2)$. 
Define
\begin{align*}
D(\sigma_1^2,v^2)
\equiv &~ 
\sum_{i=1}^m\big\|\bm{Z}_{i}\big(\tilde{\bm{b}}_{i}-\bm{b}_{i})\big\|^2-
\sum_{i=1}^m
\big\|\bm{Z}_{i}\big(\hat{\bm{b}}_{i}(\sigma_1^2,v^2)-\bm{b}_{i}\big)\big\|^2.
\end{align*}
\noindent Then, we show in Appendix~\ref{appendix:proofs} of the supplementary
material that 
\begin{align}
nD(\hat{\sigma}_1^2,\hat{v}^2)
=&~ G_{n,m}+o_p(1),
\notag
\end{align}
\noindent where $\hat\sigma_1^2$ and $\hat{v}^2$ 
are   % AMH: check
the ML estimators of
$\sigma_1^2$ and $v^2$, and $G_{n,m}$ is some random variable depending on
$n,m$. Moreover, it is shown in the same appendix that the moments of $G_{n,m}$
do not exist for $m\leq 4$ and
\begin{align}
\mathrm{E}(G_{n,m})
=&~\frac{m(m-4)v_0^4}{(m-2)\sigma_{1,0}^2}
\label{eq:limiting prob}
\end{align} 
\noindent for $m>4$. Equation~\eqref{eq:limiting prob} reveals that for any
fixed $m>4$, the empirical BLUP,
$\bm{Z}_i\hat{\bm{b}}_i(\hat\sigma_1^2,\hat{v}^2)$ of $\bm{Z}_i\bm{b}_i$, is
asymptotically more efficient than its LS counterpart,
$\bm{Z}_i\tilde{\bm{b}}_i$, even when $\hat\sigma_1^2$ is not a consistent
estimator of $\sigma_1^2$. In addition, the advantage of the former over the
latter rapidly increases with $m$. 
\label{prop:prediction}
\end{exa}

\subsection{Asymptotics under misspecification}
\label{section:incorrect}

In this subsection, we consider a misspecified model
$(\alpha,\gamma)\in(\mathcal{A}\times\mathcal{G})
\setminus(\mathcal{A}_0\times\mathcal{G}_0)$. We derive not only the
convergence rates for $\hat{v}^2(\alpha,\gamma)$ and
$\{\hat{\sigma}_k^2(\alpha,\gamma):k\in\gamma\}$, but also their asymptotic
distributions.
These results are crucial for developing model selection consistency and
efficiency in linear mixed-effects models under fixed $m$; see
Chang et al.~\cite{Chang2020}.

We start by investigating the asymptotic properties for the ML estimators of
$v^2$ and $\{\sigma_k^2:k\in\gamma\}$ for
$(\alpha,\gamma)\in\mathcal{A}_0\times(\mathcal{G}\setminus\mathcal{G}_0)$
under a misspecified random-effects model.

\begin{thm}
	Under the assumptions of Theorem~\ref{appendix:theorem:MLE}, except
	that $(\alpha,\gamma)\in\mathcal{A}_0\times(\mathcal{G}\setminus\mathcal{G}_0)$,
	\begin{align}
	\begin{split}
	\hat{v}^2(\alpha,\gamma)
	=&~ v_0^2 + \frac{1}{N}\sum_{i=1}^m \bigg(n_i^{\ell} \sum_{k\in\gamma_0\setminus\gamma}c_{i,k}b_{i,k}^2\bigg)
	+ o_p\bigg(\frac{1}{N}\sum_{i=1}^mn_i^\ell\bigg)\\
	&~  +O_p\Big(\frac{p+mq}{N}\Big) +O_p(N^{-1/2})
	\end{split}
	\label{appendix:thm:mle2:correct:eq1}
	\end{align}
	\noindent and
	\begin{align}
	\hat{\sigma}_k^2(\alpha,\gamma)
	=&~ \left\{
	\begin{array}{ll}
	\displaystyle\frac{1}{m}\sum_{i=1}^m b_{i,k}^2+\displaystyle o_p(a_N(\xi,\ell))+o_p(1);  &
	\mbox{if }k\in\gamma\cap\gamma_0, \\
	\displaystyle o_p(a_N(\xi,\ell))+o_p(1); & \mbox{if }k\in\gamma\setminus\gamma_0,
	\label{appendix:thm:mle2:correct:eq2}
	\end{array}
	\right.
	\end{align}
	\noindent where  $a_N(\xi,\ell)=\displaystyle \bigg(\frac{\sum_{i=1}^mn_i^{\ell}}{\sum_{i=1}^mn_i^{\xi}}\bigg)\bigg(\frac{\sum_{i=1}^mn_i^{\xi-\ell}}{m}\bigg)$. 
	In addition, if $\ell<1$, then
	\[
	\hat{v}^2(\alpha,\gamma)\xrightarrow{p}v_0^2,\quad\mbox{as }N\rightarrow\infty.
	\]
	Furthermore, if $\ell\in(0,1/2)$ and $p+mq=o\big(N^{1/2}\big)$, then
	\begin{align*}
	N^{1/2}(\hat{v}^2(\alpha,\gamma) - v_0^2) \xrightarrow{d} N(0,2v_0^4),\quad\mbox{as }N\rightarrow\infty.
	\end{align*}
	\label{appendix:theorem:MLE 2}
\end{thm}

Note that $\displaystyle\frac{1}{N}\sum_{i=1}^m \bigg(n_i^{\ell}
\sum_{k\in\gamma_0\setminus\gamma}c_{i,k}b_{i,k}^2\bigg)$ in
\eqref{appendix:thm:mle2:correct:eq1} is the 
% dominate bias term 
  dominant bias term   % AMH: check
for
$\hat{v}^2(\alpha, \gamma)$, which is contributed by the non-negligible random
effects missed by model~$\gamma$.
It is asymptotically positive with probability one when $\ell=1$. Hence
$\hat{v}^2(\alpha,\gamma)$ has a non-negligible positive bias when $\ell=1$.
On the other hand, for $\xi=\ell$ or nearly balanced data,
the following corollary shows that
$\hat{\sigma}_k^2(\alpha,\gamma)\xrightarrow{p}\sigma_k^2$; $k\in\gamma$, as $m\rightarrow \infty$,
even though
$(\alpha,\gamma)\in\mathcal{A}_0\times(\mathcal{G}\setminus\mathcal{G}_0)$ is
misspecified.

\begin{coro}
	Under the assumptions of Theorem~\ref{appendix:theorem:MLE 2}, with
	$\xi=\ell$ or $n_{\max}=O(n_{\min})$,
	\begin{align*}
	\hat{\sigma}_k^2(\alpha,\gamma)
	=&~ \left\{
	\begin{array}{ll}
	\displaystyle\frac{1}{m}\sum_{i=1}^m b_{i,k}^2+o_p(1);  &
	\mbox{if }k\in\gamma\cap\gamma_0, \\
	o_p(1); & \mbox{if }k\in\gamma\setminus\gamma_0.
	\end{array}
	\right.
	\end{align*}
	
	\noindent If $m\rightarrow\infty$, then
	\[
	\hat{\sigma}_k^2(\alpha,\gamma)\xrightarrow{p}\sigma_{k,0}^2; \quad k\in\gamma,\quad\mbox{as }N\rightarrow\infty.
	\]
	\label{coro:MLE 2}
\end{coro}

The following theorem presents the asymptotic properties of
$\hat{v}^2(\alpha,\gamma)$ and $\{\hat{\sigma}_k^2(\alpha,\gamma):k\in\gamma\}$
for $(\alpha,\gamma)\in(\mathcal{A}\setminus\mathcal{A}_0)\times\mathcal{G}$
under a misspecified fixed-effects model.

\begin{thm}
	Under the assumptions of Theorem~\ref{appendix:theorem:MLE} except that
	$(\alpha,\gamma)\in(\mathcal{A}\setminus\mathcal{A}_0)\times\mathcal{G}_0$,
	\begin{align}
	\begin{split}
	\hat{v}^2(\alpha,\gamma)
	=&~ v_0^2 +\frac{1}{N}\sum_{i=1}^m\bigg(n_i^{\xi} \sum_{j\in\alpha_0\setminus\alpha}d_{i,j}\beta_{j,0}^2\bigg)
	+o_p\bigg(\frac{1}{N}\sum_{i=1}^mn_i^{\xi}\bigg)\\
	&~   +O_p\Big(\frac{p+mq}{N}\Big)  + O_p(N^{-1/2})
	\end{split}
	\label{appendix:thm:mle3:correct:eq1}
	\end{align}
	\noindent and
	\begin{align}
	\hat{\sigma}_k^2(\alpha,\gamma)
	=&~\left\{
	\begin{array}{ll}
	\displaystyle\frac{1}{m}\sum_{i=1}^m b_{i,k}^2
	+o_p\bigg(\frac{1}{m}\sum_{i=1}^mn_i^{\xi-\ell}\bigg)+o_p(1);  &
	\mbox{if }k\in\gamma\cap\gamma_0, \\
	\displaystyle o_p\bigg(\frac{1}{m}\sum_{i=1}^mn_i^{\xi-\ell}\bigg)+o_p(1); & \mbox{if }k\in\gamma\setminus\gamma_0.
	\end{array}
	\right.
	\label{appendix:thm:mle3:correct:eq2}
	\end{align}
	\noindent In addition, if $\xi<1$, then
	\[
	\hat{v}^2(\alpha,\gamma)\xrightarrow{p}v_0^2,\quad\mbox{as }N\rightarrow\infty.
	\]
	Furthermore, if $\xi\in(0,1/2)$ and $p+mq=o\big(N^{1/2}\big)$, then
	\begin{align*}
	N^{1/2}(\hat{v}^2(\alpha,\gamma) - v_0^2) \xrightarrow{d} N(0,2v_0^4),\quad\mbox{as }N\rightarrow\infty.
	\end{align*}
	\label{appendix:theorem:MLE 3}
\end{thm}

Note that $\displaystyle\frac{1}{N}\sum_{i=1}^m\bigg(n_i^{\xi}
\sum_{j\in\alpha_0\setminus\alpha}d_{i,j}\beta_{j,0}^2\bigg)$ in
\eqref{appendix:thm:mle3:correct:eq1} is asymptotically positive with
probability one when $\xi=1$.
Therefore, under the assumptions of Theorem~\ref{appendix:theorem:MLE 3},
$\hat{v}^2(\alpha,\gamma)$ has a non-negligible positive bias when $\xi=1$.
Nevertheless, $\hat{\sigma}_k^2(\alpha,\gamma)$ is consistent for
$\gamma\in\mathcal{G}_0$ when $\xi\leq\ell$, as $m\rightarrow\infty$.

The following theorem establishes the asymptotic properties of
$\hat{v}^2(\alpha,\gamma)$ and $\{\hat{\sigma}_k^2(\alpha,\gamma):k\in\gamma\}$
for
$(\alpha,\gamma)\in(\mathcal{A}\setminus\mathcal{A}_0)\times(\mathcal{G}\setminus\mathcal{G}_0)$
when both the fixed-effects model and the random-effects model are
misspecified.

\begin{thm}
	Under the assumptions of Theorem~\ref{appendix:theorem:MLE} except that
	$(\alpha,\gamma)\in(\mathcal{A}\setminus\mathcal{A}_0)\times(\mathcal{G}\setminus\mathcal{G}_0)$,
	\begin{align}
	\begin{split}
	\hat{v}^2(\alpha,\gamma)
	=&~ v_0^2 +\frac{1}{N}\sum_{i=1}^m \Bigg(n_i^{\xi}\sum_{j\in\alpha_0\setminus\alpha}d_{i,j}\beta_{j,0}^2
	+n_i^{\ell}\sum_{k\in\gamma_0\setminus\gamma}c_{i,k}b_{i,k}^2\Bigg)\\
	&~  +o_p\bigg(\frac{1}{N}\sum_{i=1}^m(n_i^{\xi}+n_i^{\ell})\bigg)
	+O_p\Big(\frac{p+mq}{N}\Big)
	+O_p(N^{-1/2})
	\end{split}
	\label{appendix:thm:mle4:correct:eq1}
	\end{align}
	\noindent and 
	\begin{align}
	\hat{\sigma}_k^2(\alpha,\gamma)
	=&~ \left\{
	\begin{array}{ll}
	\displaystyle\frac{1}{m}\sum_{i=1}^m b_{i,k}^2
	+\displaystyle o_p(a_N^*(\xi,\ell))+o_p(1); &
	\mbox{if }k\in\gamma\cap\gamma_0, \\
	\displaystyle o_p(a_N^*(\xi,\ell))+o_p(1); &
	\mbox{if }k\in\gamma\setminus\gamma_0,
	\end{array}
	\right.
	\label{appendix:thm:mle4:correct:eq2}
	\end{align}
	\noindent where  $a_N^*(\xi,\ell)=\displaystyle\bigg(1+\frac{\sum_{i=1}^mn_i^{\ell}}{\sum_{i=1}^mn_i^{\xi}}\bigg)\bigg(\frac{\sum_{i=1}^mn_i^{\xi-\ell}}{m}\bigg)$. In addition, if $\max\{\xi,\ell\}<1$, then
	\[
	\hat{v}^2(\alpha,\gamma)\xrightarrow{p}v_0^2,\quad\mbox{as }N\rightarrow\infty.
	\]
	Furthermore, if $(\xi,\ell)\in(0,1/2)\times(0,1/2)$ and $p+mq=o\big(N^{1/2}\big)$,  then
	\begin{align*}
	N^{1/2}(\hat{v}^2(\alpha,\gamma) - v_0^2) \xrightarrow{d} N(0,2v_0^4),\quad\mbox{as }N\rightarrow\infty.
	\end{align*}
	\label{appendix:theorem:MLE 4}
\end{thm}

Note that $\displaystyle\frac{1}{N}\sum_{i=1}^m \Bigg(n_i^{\xi}\sum_{j\in\alpha_0\setminus\alpha}d_{i,j}\beta_{j,0}^2
+n_i^{\ell}\sum_{k\in\gamma_0\setminus\gamma}c_{i,k}b_{i,k}^2\Bigg)$
in \eqref{appendix:thm:mle4:correct:eq1} is asymptotically positive with
probability one when either $\xi=1$ or $\ell=1$.
Therefore, under the assumptions of Theorem~\ref{appendix:theorem:MLE 4},
$\hat{v}^2(\alpha,\gamma)$ has a non-negligible positive bias when either
$\xi=1$ or $\ell=1$.
%In addition, when $\xi\leq \ell$, $\hat{\sigma}_k^2(\alpha,\gamma)$ is consistent.
Also, we have the following corollary.

\begin{coro}
	Under the assumptions of Theorem~\ref{appendix:theorem:MLE 4}, with either
	(i) $\xi=\ell$ or (ii) $\xi<\ell$ and $n_{\max}=O(n_{\min})$,
	\begin{align*}
	\hat{\sigma}_k^2(\alpha,\gamma)
	=&~ \left\{
	\begin{array}{ll}
	\displaystyle\frac{1}{m}\sum_{i=1}^m b_{i,k}^2+o_p(1);  &
	\mbox{if }k\in\gamma\cap\gamma_0, \\
	o_p(1); & \mbox{if }k\in\gamma\setminus\gamma_0.
	\end{array}
	\right.
	\end{align*}
	\noindent If $m\rightarrow\infty$, then
	\[
	\hat{\sigma}_k^2(\alpha,\gamma)\xrightarrow{p}\sigma_{k,0}^2; \quad k\in\gamma,\quad\mbox{as }N\rightarrow\infty.
	\]
\end{coro}

\section{Simulations}
\label{section:simulation}

We conduct two simulation experiments for linear mixed-effects models.
The first one examines estimation of mixed-effects models,
and the second concerns confidence intervals.

\subsection{Experiment 1}
\label{section:increasing m}

We generate data according to \eqref{data:general} with $p=q=5$,
$(\sigma_{1,0}^2,\sigma_{2,0}^2,\sigma_{3,0}^2,\sigma_{4,0}^2,\sigma_{5,0}^2)'=(0,0.5,1,1.5,0)'$,
$\bm{\beta}_0=(1.2,-0.7,0.8,0,0)'$, and $v^2=1$, where $\bm{x}_{i,j}\sim
N(\bm{0},\bm{I}_{n_i})$ and $\bm{z}_{i,k}\sim N(\bm{0},\bm{I}_{n_i})$ are
independent,
for $i=1,\dots,m$ and $j,k=1,\dots,5$.
This setup satisfies (A1)--(A3) with $\xi=\ell=1$ and $d_{i,j}=c_{i,k}=1$, for
$i=1,\dots,m$ and $j,k=1,\dots,5$.
We consider parameter estimation under two scenarios corresponding to balanced
data and unbalanced data.
We also consider model selection under balanced data.

For parameter estimation, we first consider balanced data with
$m\in\{10,20,30\}$, $n_1=\cdots=n_m=m$, and hence $N=m^2$.
The ML estimators of $\sigma_1^2,\dots,\sigma_5^2$ and $v^2$ under the full model
$(\{1,\dots,5\},\{1,\dots,5\})\in\mathcal{A}_0\times\mathcal{G}_0$ based on
100 simulated replicates are summarized in Table~\ref{table:balanced true}.
The ML estimators of $\sigma_1^2,\dots,\sigma_3^2$ and $v^2$ under model
$(\{1,2,3\},\{1,2,3\})\in\mathcal{A}_0\times(\mathcal{G}\setminus\mathcal{G}_0)$
with correct fixed effects but misspecified random effects based on 100
simulated replicates are summarized in Table~\ref{table:balanced incorrect}.
The ML estimators of $\sigma_4^2,\sigma_5^2$, and $v^2$ under model
$(\{2,3,4,5\},\{4,5\})\in(\mathcal{A}\setminus\mathcal{A}_0)\times(\mathcal{G}\setminus\mathcal{G}_0)$
with both misspecified fixed and random effects based on 100 simulated
replicates are summarized in Table~\ref{table:balanced incorrect 2}.

\begin{table}[H]\centering
	\caption{Sample means and sample standard deviations (in parentheses) of 
	ML estimators of $\sigma_1^2,\dots,\sigma_5^2$ and $v^2$ for different
	values of $m$ obtained from full model in 
	Section~\ref{section:increasing m} with balanced data based on 100 simulated
	replicates. Values in row for $m=\infty$ are probability limits
	of ML estimators.}
	\smallskip
	\begin{tabular}{ccccccc}
		\hline
		$m$ & $\hat{\sigma}_1^2$ & $\hat{\sigma}_2^2$ & $\hat{\sigma}_3^2$ & $\hat{\sigma}_4^2$ & $\hat{\sigma}_5^2$ & $\hat{v}^2$ \\\hline
		10 & 0.033    & 0.467   & 0.994   & 1.453    & 0.041   & 0.850 \\
		& (0.062) & (0.281) & (0.564) & (0.615) & (0.073) & (0.165) \\
		20 & 0.008    & 0.512   & 1.028   & 1.470    & 0.008   & 0.983 \\
		& (0.017) & (0.194) & (0.362) & (0.490) & (0.014) & (0.087) \\
		30 & 0.003    & 0.490   & 0.994   & 1.534    & 0.004   & 0.989 \\
		& (0.006) & (0.116) & (0.260) & (0.396) & (0.007) & (0.049) \\
		$\infty$ & 0.000 & 0.500 & 1.000 & 1.500 & 0.000 & 1.000 \\
		True & 0.000 & 0.500 & 1.000 & 1.500 & 0.000 & 1.000 \\
		\hline
	\end{tabular}
	\label{table:balanced true}
\end{table}

\begin{table}[H]\centering
	\caption{Sample means and sample standard deviations (in parentheses) of 
	ML estimators of $\sigma_1^2,\sigma_2^2,\sigma_3^2$ and $v^2$ for different
	values of $m$ obtained from model $(\alpha,\gamma)=(\{1,2,3\},\{1,2,3\})$ in
	Section~\ref{section:increasing m} with balanced data based on 100 simulated
	replicates.
		Values in row for $m=\infty$ are probability limits of ML
		estimators.}
	\smallskip
	\begin{tabular}{cllll}
		\hline
		$m$ & ~~~~~~$\hat{\sigma}_1^2$ & ~~~~~~$\hat{\sigma}^2_2$ & ~~~~~~$\hat{\sigma}^2_3$ & ~~~~~~$\hat{v}^2$ \\\hline
		10 & 0.151 (0.541) & 0.583 (0.546) & 0.944 (0.572) & 2.363 (1.054) \\
		20 & 0.047 (0.091) & 0.555 (0.285) & 1.050 (0.392) & 2.415 (0.650) \\
		30 & 0.030 (0.070) & 0.521 (0.198) & 0.961 (0.268) & 2.442 (0.666) \\
		$\infty$ & 0.000 & 0.500 & 1.000 & 2.500 \\
		True & 0.000 & 0.500 & 1.000 & 1.000 \\
		\hline
	\end{tabular}
	\label{table:balanced incorrect}
\end{table}

\begin{table}[H]\centering
	\caption{Sample means and sample standard deviations (in parentheses) of 
	ML estimators of $\sigma_4^2,\sigma_5^2$ and $v^2$ for different values of
	$m$ obtained from model $(\alpha,\gamma)=(\{2,3,4,5\},\{4,5\})$ in 
	Section~\ref{section:increasing m} with balanced data based on 100 simulated
	replicates. 
		Values in row for $m=\infty$ are probability limits of ML
		estimators.}
	\smallskip
	\begin{tabular}{clll}
		\hline
		$m$ & ~~~~~~$\hat{\sigma}_4^2$ & ~~~~~~$\hat{\sigma}_5^2$ & ~~~~~~$\hat{v}^2$ \\\hline
		10 & 1.604 (1.073) & 0.176 (0.465) & 3.494 (0.788)\\
		20 & 1.353 (0.567) & 0.043 (0.077) & 3.915 (0.540)\\
		30 & 1.525 (0.427) & 0.030 (0.065) & 3.880 (0.436)\\
		$\infty$ & 1.500 & 0.000 & 3.690 \\
		True & 1.500 & 0.000 & 1.000\\
		\hline
	\end{tabular}
	\label{table:balanced incorrect 2}
\end{table}

As seen in Table~\ref{table:balanced true}, the ML estimators,
$\hat{\sigma}_1^2,\dots,\hat{\sigma}_5^2$ and $\hat{v}^2$, based on the full
model, have small biases except for $\hat{v}^2$ with $m=10$.
We note that their standard deviations tend to be smaller when $m$ is larger.
In particular, the standard deviations of $\hat{\sigma}_1^2$ and
$\hat{\sigma}_5^2$ are much smaller than the others, which echoes 
% Theorem~\ref{appendix:theorem:MLE}, that 
  Theorem~\ref{appendix:theorem:MLE}, that which shows that   % AMH: check
$\hat{v}^2_j$ has a faster convergence rate
when it converges to zero.
For model $(\alpha,\gamma)=(\{1,2,3\},\{1,2,3\})$ with misspecified random
effects, Table~\ref{table:balanced incorrect} shows that the ML estimator
$\hat{v}^2$ overestimates $v_0^2=1$ by about $\sigma_{4,0}^2=1.5$ on average,
particularly for larger values of $m$, which also complies with 
Theorem~\ref{appendix:theorem:MLE 2}.
Finally, for model $(\alpha,\gamma)=(\{2,3,4,5\},\{4,5\})$ with both fixed and
random effects misspecified, Table~\ref{table:balanced incorrect 2} confirms
that $\hat{v}^2$ is far from its true value and reasonably close to its
probability limit, $v_0^2+ \sigma_{2,0}^2+\sigma_{3,0}^2+\beta_{1,0}^2=3.69$,
derived in Theorem~\ref{appendix:theorem:MLE 4}.
In addition, $\hat{\sigma}_4^2$ tends to be closer to $\sigma_{4,0}^2$ when $m$
is larger, as expected from Theorem~\ref{appendix:theorem:MLE 4}.

Next, we consider unbalanced data with $m\in\{10,20,30\}$ and $N=m^2$.
We set $n_1=[N^{1/4}]$, $n_2=[N^{3/4}]$, $n_3=\cdots=n_{m-1}=[(N-n_1-n_2)/(m-2)]$,
and hence $n_m=N-\sum_{i=1}^{m-1}n_i$.
The ML estimators of $\sigma_1^2,\dots,\sigma_5^2$ and $v^2$ under the
full model $(\{1,\dots,5\},\{1,\dots,5\})\in\mathcal{A}_0\times\mathcal{G}_0$
based on 100 simulated replicates are summarized in Table~\ref{table:unbalanced true}.
The ML estimators of $\sigma_1^2,\dots,\sigma_3^2$ and $v^2$ under model
$(\{1,2,3\},\{1,2,3\})\in\mathcal{A}_0\times(\mathcal{G}\setminus\mathcal{G}_0)$
with correct fixed effects but misspecified random effects based on 100
simulated replicates are summarized in Table~\ref{table:unbalanced incorrect}.
The ML estimators of $\sigma_4^2$, $\sigma_5^2$, and $v^2$ under model
$(\{2,3,4,5\},\{4,5\})\in(\mathcal{A}\setminus\mathcal{A}_0)\times(\mathcal{G}\setminus\mathcal{G}_0)$
with both misspecified fixed and random effects based on 100 simulated
replicates are summarized in Table~\ref{table:unbalanced incorrect 2}.
The ML estimators of $\sigma_1^2,\dots,\sigma_5^2$ and $v^2$ based on
unbalanced data can be seen to perform similarly to those based on balanced
data.

\begin{table}[H]\centering
	\caption{Sample means and sample standard deviations (in parentheses) of 
	ML estimators of $\sigma_1^2,\dots,\sigma_5^2$ and $v^2$ for different
	values of $m$ obtained from full model in 
	Section~	\ref{section:increasing m} with unbalanced data based on 100 simulated
	replicates. Values in row for $m=\infty$ are probability limits
	of ML estimators.}
	\smallskip
	\begin{tabular}{ccrcccccc}
		\hline
		$m$ & $n_{\min}$ & \!\!$n_{\max}$\!\! & $\hat\sigma_1^2$ & $\hat\sigma_2^2$ & $\hat\sigma_3^2$ & $\hat\sigma_4^2$ &
		$\hat\sigma_5^2$ & $\hat{v}^2$\\\hline
		10 & 3 &  32  & 0.017   & 0.500    & 1.011   & 1.414    & 0.021   & 0.877 \\
		&    &        & (0.041) & (0.330) & (0.602) & (0.790) & (0.040) & (0.154) \\
		20 & 4 &  89  & 0.009   & 0.516    & 1.029   & 1.490    & 0.008   & 0.974 \\
		&    &        & (0.020) & (0.175) & (0.399) & (0.493) & (0.014) & (0.082) \\
		30 & 5 & 164 & 0.002   & 0.497    & 1.007   & 1.539   & 0.004    & 0.991 \\
		&    &        & (0.005) & (0.121) & (0.263) & (0.374) & (0.008) & (0.049) \\
		$\infty$ &&& 0.000 & 0.500 & 1.000 & 1.500 & 0.000 & 1.000 \\
		True &&& 0.000 & 0.500 & 1.000 & 1.500 & 0.000 & 1.000 \\
		\hline
	\end{tabular}
	\label{table:unbalanced true}
\end{table}

\begin{table}[H]\centering
	\caption{Sample means and sample standard deviations (in parentheses) of 
	ML estimators of $\sigma_1^2,\sigma_2^2,\sigma_3^2$, and $v^2$ for different
	values of $m$ obtained from model $(\alpha,\gamma)=(\{1,2,3\},\{1,2,3\})$ in
	Section~\ref{section:increasing m} with unbalanced data based on 100
	simulated replicates. Values in row for $m=\infty$ are
	probability limits of ML estimators.}
	\smallskip
	\begin{tabular}{ccrllll}
		\hline
		$m$ & $n_{\min}$ & \!\!$n_{\max}$\!\! & ~~~~~~$\hat\sigma_1^2$ & ~~~~~~$\hat\sigma_2^2$ & ~~~~~~$\hat\sigma_3^2$ & ~~~~~~$\hat{v}^2$ \\\hline
		10 & 3 &   32 & 0.213 (0.927) & 0.576 (0.696) & 1.175 (1.220) & 2.283 (1.127) \\
		20 & 4 &   89 & 0.044 (0.096) & 0.536 (0.260) & 1.091 (0.461) & 2.456 (0.792) \\
		30 & 5 & 165 &  0.028 (0.079) & 0.500 (0.208) & 0.959 (0.290) & 2.426 (0.645)\\
		$\infty$ &&& 0.000 & 0.500 & 1.000 & 2.500 \\
		True &&& 0.000 & 0.500 & 1.000 & 1.000 \\
		\hline
	\end{tabular}
	\label{table:unbalanced incorrect}
\end{table}

\begin{table}[H]\centering
	\caption{Sample means and sample standard deviations (in parentheses) of 
	ML estimators of $\sigma_4^2,\sigma_5^2$, and $v^2$ for different values of
	$m$ obtained from model $(\alpha,\gamma)=(\{2,3,4,5\},\{4,5\})$ in 
	Section~\ref{section:increasing m} with unbalanced data based on 100 simulated
	replicates. Values in row for $m=\infty$ are probability limits
	of ML estimators.}
	\smallskip
	\begin{tabular}{ccrlll}
		\hline
		$m$ & $n_{\min}$ & $n_{\max}$ & ~~~~~~$\hat{\sigma}_4^2$ & ~~~~~~$\hat{\sigma}_5^2$ & ~~~~~~$\hat{v}^2$ \\\hline
		10 & 3 &  32  & 1.522 (0.902) & 0.065 (0.180) & 3.535 (0.944) \\
		20 & 4 &  89  & 1.362 (0.540) & 0.057 (0.142) & 3.960 (0.716) \\
		30 & 5 & 164 &  1.494 (0.458) & 0.030 (0.068) & 3.892 (0.539) \\
		$\infty$ &&& 1.500 & 0.000 & 3.690 \\
		True &&& 1.500 & 0.000 & 1.000\\
		\hline
	\end{tabular}
	\label{table:unbalanced incorrect 2}
\end{table}

\subsection{Experiment 2}
\label{section:fixed m}

In the second experiment, we compare the conventional confidence interval given
by \eqref{CI:standard normal} with the proposed confidence interval given by
\eqref{CI:chisquare}.
Similar to Experiment~1, we generate data according to \eqref{data:general}
with $p=q=5$, $\bm{\beta}=(1.2,-0.7,0.8,0,0)'$, $v^2=1$, and
$(\sigma_{1,0}^2,\sigma_{2,0}^2,\sigma_{3,0}^2,\sigma_{4,0}^2,\sigma_{5,0}^2)'=(0,0.5,1,1.5,0)'$,
where $\bm{x}_{i,j}\sim N(\bm{0},\bm{\Sigma}_x)$ and $\bm{z}_{i,k}\sim
N(\bm{0},\bm{\Sigma}_z)$ are independent, for $i=1,\dots,m$ and
$j,k=1,\dots,5$.
Here we consider a more challenging situation of dependent covariates.
Specifically, we assume that $\bm{\Sigma}_x$ is a $5\times 5$ matrix with the
$(i,j)$-th entry $0.4^{|i-j|}$, and $\bm{\Sigma}_z$ is a $5\times 5$ matrix
with the $(i,j)$-th entry $0.6^{|i-j|}$.
We consider balanced data with $n=n_1=\cdots=n_m\in\{10,50,100\}$ and three
numbers of clusters, $m\in\{2,5,10\}$, resulting in a total of nine different
combinations.

We compare the 95$\%$ confidence intervals of \eqref{CI:standard normal} and
\eqref{CI:chisquare} for $\sigma_2^2$ and $\sigma_4^2$ based on model
$(\alpha,\gamma)=(\{1,2,3\},\{2,3,4\})$.
The coverage probabilities of both confidence intervals obtained from the two
methods for various cases based on 1,000 simulated replicates are shown in
Table~\ref{table:CI:sd and chisq}.
The proposed method has better coverage probabilities than the
conventional ones in almost all cases.
The coverage probabilities of our confidence interval tend to the nominal level
(i.e., $0.95$) as $n$ increases for all cases even when $m$ is very small.
In contrast, the conventional method tends to be too optimistic for both
$\sigma_2^2$ and $\sigma_4^2$.
For example, the coverage probabilities are less than $0.73$ when $m=2$
regardless of $n$.
Although the coverage probabilities are a bit closer to the nominal level when
$m$ is larger, they are still in the range of $(0.82,0.87)$ when $m=10$,
showing that the conventional confidence interval is not valid for small $m$.

\begin{table}[H]\centering
	\caption{ Coverage probabilities (denoted by $\hat{P}$) for $95\%$
	confidence intervals of $\sigma_2^2$ and $\sigma_{4}^2$ obtained from two
	methods in Section~\ref{section:fixed m} based on 1,000 simulated
	replicates. 
		Values given in parentheses are standard errors of coverage
		probabilities (evaluted by $\sqrt{\hat{P}(1-\hat{P})/1000}$).  }
	\smallskip
	\begin{tabular}{cccccc}
		\hline
		$m$ & $n$ & \multicolumn{2}{c}{Classical} & \multicolumn{2}{c}{Proposed}\\
		\cline{3-6}
		&& $\sigma_2^2$ & $\sigma_4^2$ & $\sigma_2^2$ & $\sigma_4^2$ \\
		\hline
		2 & 10 & 0.651 (0.015) & 0.649 (0.015) & 0.814 (0.012) & 0.763 (0.013) \\
		& 50 & 0.724 (0.014) & 0.703 (0.014) & 0.932 (0.008) & 0.935 (0.008) \\
		&100 & 0.725 (0.014) & 0.722 (0.014) & 0.942 (0.007) & 0.929 (0.008) \\
		\hline
		5 & 10 & 0.778 (0.013) & 0.738 (0.014) & 0.895 (0.010) & 0.871 (0.011) \\
		& 50 & 0.809 (0.012) & 0.818 (0.012) & 0.936 (0.008) & 0.937 (0.008) \\
		&100 & 0.811 (0.012) & 0.809 (0.012) & 0.940 (0.008) & 0.929 (0.008) \\
		\hline
		10& 10 & 0.838 (0.012) & 0.816 (0.012) & 0.900 (0.009) & 0.893 (0.010) \\
		& 50 & 0.874 (0.010) & 0.849 (0.011) & 0.952 (0.007) & 0.946 (0.007) \\
		&100 & 0.849 (0.011) & 0.867 (0.011) & 0.941 (0.007) & 0.956 (0.006) \\
		\hline
	\end{tabular}
	\label{table:CI:sd and chisq}
\end{table}

\section{Discussion}
\label{section:discussion}

In this article, we establish the asymptotic theory of the ML estimators of
random-effects parameters in linear mixed-effects models for unbalanced data,
without assuming that $m$ grows to infinity with $N$.
We not only allow the dimensions of both the fixed-effects and random-effects
models to go to infinity with $N$, but also allow both models to be
misspecified.
In addition, we provide an asymptotic valid confidence interval for the
random-effects parameters when $m$ is fixed.
These asymptotic results are essential for investigating the asymptotic
properties of model-selection methods for linear mixed-effects models, which to
the best of our knowledge have only been developed under the assumption of
$m\rightarrow\infty$.

Although it is common to assume the random effects to be uncorrelated
as done in model~\eqref{data:general},
it is also of interest to consider correlated random effects with no structure
imposed on $\bm{D}$. However, the technique developed in this article may not
be directly applicable to the latter situation; further research in this
direction is 
thus   % AMH: check
warranted.

Conditions (A1)~and~(A2) assume that the covariates are asymptotically
uncorrelated.
These restrictions can be relaxed. Here is a simple example.

\begin{lem}
	Consider the data generated from \eqref{data} with $m=1$, $n_1=N$, $p=q=2$,
	and the true parameters given in \eqref{data:linear mixed-effects}.
	Suppose that $(\alpha_0,\gamma_0)=(\{1,2\},\{1,2\})$ is the smallest true
	model and $(\alpha_1,\gamma_1)=(\{1\},\{1\})$ is a misspecified model
	defined in \eqref{model:linear mixed effects}.
	Let $\hat{\sigma}_k^2(\alpha,\gamma)$ and $\hat{v}^2(\alpha,\gamma)$ be
	the ML estimators of $\sigma_k^2$ and $v^2$ based on $(\alpha,\gamma)$.
	Assume that (A1)--(A3) hold except that $\bm{z}_{1,1}'\bm{z}_{1,2}=c_{1,12}N
	+ o(N)$ and $\bm{x}_{1,1}'\bm{x}_{1,2}=d_{1,12}N +o(N)$, for some constants
	$c_{1,12},d_{1,12}\in\mathbb{R}$.
	Then
	\begin{align*}
	\hat{v}^2(\alpha_0,\gamma_0)
	=&~ v_0^2 + O_p(N^{-1/2}),\\
	\hat{\sigma}_k^2(\alpha_0,\gamma_0)
	=&~ b_k^2 + O_p(N^{-1/2}); \quad k=1,2,\\
	\hat{v}^2(\alpha_1,\gamma_1)
	=&~ v_0^2 +
	\bigg(d_{1,2}-\frac{d_{1,12}^2}{d_{1,1}}\bigg)\beta_{2,0}^2+
	\bigg(c_{1,2}-\frac{c_{1,12}^2}{c_{1,1}}\bigg)b_2^2 + o_p(1),\\
	\hat{\sigma}_1^2(\alpha_1,\gamma_1)
	=&~ \bigg(b_{1}+\frac{c_{1,12}}{c_{1,1}}b_{2}\bigg)^2 + o_p(1),
	\end{align*}
	\noindent where $\beta_{2,0}\neq 0$ is the true parameter of $\beta_2$.
	\label{lemma:correlated covariates}
\end{lem}

From Lemma~\ref{lemma:correlated covariates}, it is not surprising to see that
$\hat{v}^2(\alpha_0,\gamma_0)\xrightarrow{p}v_0^2$.
On the other hand, $\hat{v}^2(\alpha_1,\gamma_1)$ tends to overestimate $v_0^2$ by
$(d_{1,2}-d_{1,12}^2/d_{1,1})\beta_{2,0}^2+(c_{1,2}-c_{1,12}^2/c_{1,1})b_2^2$.
Since $d_{1,2}-d_{1,12}^2/d_{1,1}\geq 0$ and $c_{1,2}-c_{1,12}^2/c_{1,1}\geq
0$, the amount of overestimation is smaller when either $c_{1,12}^2$ or
$d_{1,12}^2$ is larger.
In contrast, $\hat{\sigma}_1^2(\alpha_1,\gamma_1)$ tends to be more 
% biased upward 
  upward biased   % AMH: check
when $c_{1,12}^2$ is larger, since
$\mathrm{E}\big(b_{1}+(c_{1,12}/c_{1,1})b_2\big)^2=\sigma_1^2+(c_{1,12}/c_{1,1})^2\sigma_2^2$.
Lemma~\ref{lemma:correlated covariates} demonstrates how the correlations
between the two covariates affect the behavior of
$\hat{v}^2(\alpha_1,\gamma_1)$ and $\hat{\sigma}_1^2(\alpha_1,\gamma_1)$.
However, when the number of covariates is larger, the ML estimators of $v^2$
and $\{\sigma_k^2\}$ become much more complicated.
We leave this extension of Lemma~\ref{lemma:correlated covariates} to the
general case for future work.

%%%%%%%%%%%%%%%%%%%%%%%%%%%%%%%%%%%%%%%%%%%%%%
%% Support information (funding), if any,   %%
%% should be provided in the                %%
%% Acknowledgements section.                %%
%%%%%%%%%%%%%%%%%%%%%%%%%%%%%%%%%%%%%%%%%%%%%%
\section*{Acknowledgements}
The research of Chih-Hao Chang is supported by ROC Ministry of Science and
Technology grant MOST 107-2118-M-390-001.

The research of Hsin-Cheng Huang is supported by  ROC Ministry of Science and
Technology grant MOST 106-2118-M-001-002-MY3.

The research of Ching-Kang Ing is supported by the Science Vanguard Research
Program under the Ministry of Science and Technology, Taiwan, ROC.

%%%%%%%%%%%%%%%%%%%%%%%%%%%%%%%%%%%%%%%%%%%%%%
%% Supplementary Material, if any, should   %%
%% be provided in {supplement} environment  %%
%% with title and short description.        %%
%%%%%%%%%%%%%%%%%%%%%%%%%%%%%%%%%%%%%%%%%%%%%%

%%  imsart-number.bst  will be used to                     %%
%%  MR numbers will be added by VTeX.                      %%
%%                                                         %%
%%  Use \cite{...} to cite references in text.             %%
%%                                                         %%
%%%%%%%%%%%%%%%%%%%%%%%%%%%%%%%%%%%%%%%%%%%%%%%%%%%%%%%%%%%%%

%% if your bibliography is in bibtex format, uncomment commands:
\bibliographystyle{imsart-number} % Style BST file
%\bibliography{bibliography}       % Bibliography file (usually '*.bib')

%% or include bibliography directly:

\bibliography{mybibfilemle}

\newpage

\renewcommand{\theequation}{\thesection.\arabic{equation}}
\renewcommand\thesection{\Alph{section}}

\setcounter{section}{0}
\setcounter{equation}{0}
\begin{frontmatter}
	
	\title{}
\end{frontmatter}

\section*{Supplementary Material}
The supplementary materials consist of three appendices that prove all the
theoretical results except for Theorem~\ref{theorem:CI}, whose proof is
straightforward and is hence omitted. Appendix~\ref{appendix:lemma} contains
auxiliary lemmas that are required in the proofs. 
Appendix~\ref{appendix:proofs} provides proofs for Example~\ref{prop:prediction} and
Theorems~\ref{appendix:theorem:MLE} and 
\ref{appendix:theorem:MLE 2}--\ref{appendix:theorem:MLE 4}. Appendix~\ref{appendix:lemma proofs} gives
proofs for all the lemmas.

\begin{appendix}
	\section{Auxiliary Lemmas}
	\label{appendix:lemma}
	We start with the following matrix identities, which will be repeated
	applied:
	\begin{align}
	\det(\bm{A}+ \bm{c}\bm{d}')
	=&~ \det(\bm{A})(1+ \bm{d}'\bm{A}^{-1}\bm{c}),
	\label{eq:det decompose}\\
	(\bm{A}+ \bm{c}\bm{d}')^{-1}
	=&~ \bm{A}^{-1} - \frac{\bm{A}^{-1}\bm{c}\bm{d}'\bm{A}^{-1}}{1+\bm{d}'\bm{A}^{-1}\bm{c}},
	\label{eq:inv decompose}
	\end{align}
	
	\noindent where $\bm{A}$ is an $n\times n$ nonsingular matrix, and $\bm{c}$
	and $\bm{d}$ are $n\times 1$ column vectors.
	Note that (\ref{eq:inv decompose}) is applied iteratively to establish
	the decomposition of the precision matrix $\bm{H}_i^{-1}(\gamma,\bm\theta)$,
	where
	\begin{align}
	\bm{H}_i(\gamma,\bm\theta)
	\equiv &~ \sum_{k\in\gamma}\theta_k\bm{z}_{i,k}\bm{z}_{i,k}'+\bm{I}_{n_i}.
	\label{matrix:H:lmm}
	\end{align}
	\noindent Heuristically speaking, let $\bm{z}_{i,(s)}$;
	$s=1,\dots,q(\gamma)$ be the $s$-th column of $\bm{Z}_i(\gamma)$ and
	\begin{align}
	\bm{H}_{i,t}(\gamma,\bm\theta)
	= \sum_{s=1}^t\theta_{(s)}\bm{z}_{i,(s)}\bm{z}_{i,(s)}'+\bm{I}_{n_i};\quad t=1,\dots,q(\gamma),
	\label{proof:lem:z:eq0}
	\end{align}
	
	\noindent where $\theta_{(s)}$ denotes the $s$-th element of $\bm\theta$;
	$s=1,\dots,q(\gamma)$.
	Suppose that $q(\gamma)=q$. Then by (\ref{eq:inv decompose}),
	\begin{align}
	\bm{H}_{i,q}^{-1}(\gamma,\bm\theta)
	=&~ \bm{H}_{i,q-1}^{-1}(\gamma,\bm\theta) - \frac{\theta_{q}\bm{H}_{i,q-1}^{-1}(\gamma,\bm\theta)\bm{z}_{i,q}\bm{z}_{i,q}'\bm{H}_{i,q-1}^{-1}(\gamma,\bm\theta)}
	{1+\theta_{q}\bm{z}_{i,q}'\bm{H}_{i,q-1}^{-1}(\gamma,\bm\theta)\bm{z}_{i,q}}.
	\label{matrix:H first step}
	\end{align}
	
	\noindent Applying (\ref{eq:inv decompose}) iteratively, we obtain
	the decomposition
	\begin{align}
	\bm{H}_{i,q}^{-1}(\gamma,\bm\theta)
	=&~ \bm{I}_{n_i} - \sum_{k=1}^q \frac{\theta_{k}\bm{H}_{i,k-1}^{-1}(\gamma,\bm\theta)\bm{z}_{i,k}\bm{z}_{i,k}'
		\bm{H}_{i,k-1}^{-1}(\gamma,\bm\theta)}
	{1+\theta_{k}\bm{z}_{i,k}'\bm{H}_{i,k-1}^{-1}(\gamma,\bm\theta)\bm{z}_{i,k}};
	\label{fn:expand inv H}
	\end{align}
	\noindent note that $\bm{H}_{i,0}(\gamma,\bm\theta)=\bm{I}_{n_i}$.
	The proofs of Lemmas \ref{appendix:lemma:z x},
	\ref{appendix:lemma:z}, and
	\ref{appendix:lemma:epsilon} are then based on the induction and the
	decomposition of (\ref{fn:expand inv H}).
	
	The proofs of theorems in Section~\ref{section:MLE} heavily rely on the asymptotic properties of the quadratic forms, $\bm{x}'_{i,j}\bm{H}_i^{-1}(\gamma,\bm\theta)\bm{x}_{i,j^*}$,
	$\bm{z}'_{i,k}\bm{H}_i^{-1}(\gamma,\bm\theta)\bm{z}_{i,k^*}$, $\bm{\epsilon}'_i\bm{H}_i^{-1}(\gamma,\bm\theta)\bm{\epsilon}_i$,
	$\bm{x}'_{i,j}\bm{H}_i^{-1}(\gamma,\bm\theta)\bm{z}_{i,k}$,  $\bm{x}'_{i,j}\bm{H}_i^{-1}(\gamma,\bm\theta)\bm{\epsilon}_i$,
	and  $\bm{z}'_{i,k}\bm{H}_i^{-1}(\gamma,\bm\theta)\bm{\epsilon}_i$,
	with $\bm{H}_i(\gamma,\bm{\theta})$ defined in \eqref{matrix:H:lmm},
	for $i=1,\dots,m$; $j,j^*=1,\dots,p$ and $k,k^*=1,\dots,q$.
	The following lemmas give their convergence rates.
	
	\begin{lem}
		Consider the linear mixed-effects model $(\alpha,\gamma)$ of
		\eqref{model:linear mixed effects}.
		Suppose that (A0)--(A3) hold. Then for $\bm{H}_i(\gamma,\bm\theta)$
		defined in (\ref{matrix:H:lmm}), we have
		\begin{enumerate}
			\item[(i)] For $i=1,\dots,m$ and $j,j^*=1,\dots,p$,
			\begin{align*}
			\sup_{\bm\theta\in[0,\infty)^{q(\gamma)}}\big|\bm{x}_{i,j}'\bm{H}_i^{-1}(\gamma,\bm\theta)\bm{x}_{i,j^*}\big|
			=&~ \left\{
			\begin{array}{ll} d_{i,j}n_i^\xi + o(n_i^\xi);
			& \mbox{if }j=j^*,\\
			o(n_i^{\xi-\tau}); &\mbox{if } j\neq j^*.
			\end{array}\right.
			%\label{appendix:lemma:z x:eq1}
			\end{align*}
			\item[(ii)] For $i=1,\dots,m$, $j=1,\dots,p$ and $k\notin\gamma$,
			\begin{align*}
			\sup_{\bm\theta\in[0,\infty)^{q(\gamma)}}\big|\bm{x}_{i,j}'\bm{H}_i^{-1}(\gamma,\bm\theta)\bm{z}_{i,k}\big|
			=&~ o(n_i^{(\xi+\ell)/2-\tau}).
			%\label{appendix:lemma:z x:eq3}
			\end{align*}
			\item[(iii)] For $i=1,\dots,m$, $j=1,\dots,p$ and $k\in\gamma$,
			\begin{align*}
			\begin{split}
			\sup_{\bm\theta\in[0,\infty)^{q(\gamma)}}\theta_k\big|\bm{x}_{i,j}'\bm{H}_i^{-1}(\gamma,\bm\theta)\bm{z}_{i,k}\big|
			=&~ o_p(n_i^{(\xi-\ell)/2-\tau}),\\
			\sup_{\bm\theta\in[0,\infty)^{q(\gamma)}}\big|\bm{x}_{i,j}'\bm{H}_i^{-1}(\gamma,\bm\theta)\bm{z}_{i,k}\big|
			=&~ o(n_i^{(\xi+\ell)/2-\tau}).
			\end{split}
			%\label{appendix:lemma:z x:eq4}
			\end{align*}
		\end{enumerate}
		\label{appendix:lemma:z x}
	\end{lem}
	
	\begin{lem}
		Consider the linear mixed-effects model $(\alpha,\gamma)$ of
		\eqref{model:linear mixed effects}.
		Suppose that (A0) and (A2) hold. Then for $\bm{H}_i(\gamma,\bm\theta)$
		defined in (\ref{matrix:H:lmm}), we have
		\begin{enumerate}
			\item[(i)] For $i=1,\dots,m$ and $k,k^*\notin\gamma$,
			\begin{align*}
			\sup_{\bm\theta\in[0,\infty)^{q(\gamma)}}\big|\bm{z}_{i,k}'\bm{H}_i^{-1}(\gamma,\bm\theta)\bm{z}_{i,k^*}\big|
			=&~ \left\{
			\begin{array}{ll}
			c_{i,k}n_i^\ell + o(n_i^\ell); & \mbox{if }k=k^*,\\
			o(n_i^{\ell-\tau}); & \mbox{if }k\neq k^*.
			\end{array}
			\right.
			%\label{appendix:lemma:z:eq1}
			\end{align*}
			\item[(ii)] For $i=1,\dots,m$ and $k\in \gamma$,
			\begin{align*}
			\begin{split}
			\sup_{\bm\theta\in[0,\infty)^{q(\gamma)}}\big|\theta_k^2\bm{z}_{i,k}'\bm{H}_i^{-1}(\gamma,\bm\theta)\bm{z}_{i,k}-\theta_k\big|
			=&~ O(n_i^{-\ell}),\\
			\sup_{\bm\theta\in[0,\infty)^{q(\gamma)}}\big|\bm{z}_{i,k}'\bm{H}_i^{-1}(\gamma,\bm\theta)\bm{z}_{i,k}\big|
			=&~ O(n_i^\ell).
			\end{split}
			%\label{appendix:lemma:z:eq3}
			\end{align*}
			\item[(iii)] For $i=1,\dots,m$ and $k,k^*\in\gamma$ with $k\neq k^*$,
			\begin{align*}
			\begin{split}
			\sup_{\bm\theta\in[0,\infty)^{q(\gamma)}}\theta_k\theta_{k^*}\big|\bm{z}_{i,k}'\bm{H}_i^{-1}(\gamma,\bm\theta)\bm{z}_{i,k^*}\big|
			=&~ o(n_i^{-\ell-\tau}),\\
			\sup_{\bm\theta\in[0,\infty)^{q(\gamma)}}\theta_{k}\big|\bm{z}_{i,k}'\bm{H}_i^{-1}(\gamma,\bm\theta)\bm{z}_{i,k^*}\big|
			=&~ o( n_i^{-\tau}),\\
			\sup_{\bm\theta\in[0,\infty)^{q(\gamma)}}\big|\bm{z}_{i,k}'\bm{H}_i^{-1}(\gamma,\bm\theta)\bm{z}_{i,k^*}\big|
			=&~ o(n_i^{\ell-\tau}).
			\end{split}
			%\label{appendix:lemma:z:eq4}
			\end{align*}
			\item[(iv)] For $i=1,\dots,m$, $k\in\gamma$ and $k^*\notin\gamma$,
			\begin{align*}
			\begin{split}
			\sup_{\bm\theta\in[0,\infty)^{q(\gamma)}}\theta_k\big|\bm{z}_{i,k}'\bm{H}_i^{-1}(\gamma,\bm\theta)\bm{z}_{i,k^*}\big|
			=&~ o(n_i^{-\tau}),\\
			\sup_{\bm\theta\in[0,\infty)^{q(\gamma)}}\big|\bm{z}_{i,k}'\bm{H}_i^{-1}(\gamma,\bm\theta)\bm{z}_{i,k^*}\big|
			=&~ o(n_i^{\ell-\tau}).
			\end{split}
			%\label{appendix:lemma:z:eq5}
			\end{align*}
		\end{enumerate}
		\label{appendix:lemma:z}
	\end{lem}
	
	\begin{lem}
		Consider the linear mixed-effects model $(\alpha,\gamma)$ of
		\eqref{model:linear mixed effects}.
		Suppose that (A0)--(A3) hold. Then for $\bm{H}_i(\gamma,\bm\theta)$
		defined in (\ref{matrix:H:lmm}), we have
		\begin{enumerate}
			\item[(i)] For $i=1,\dots,m$ and $k\in\gamma$,
			\begin{align*}
			\begin{split}
			\sup_{\bm\theta\in[0,\infty)^{q(\gamma)}}\theta_k\big|\bm{z}_{i,k}'\bm{H}_i^{-1}(\gamma,\bm\theta)\bm\epsilon_i\big|
			=&~ O_p(n_i^{-\ell/2}),\\
			\sup_{\bm\theta\in[0,\infty)^{q(\gamma)}}\big|\bm{z}_{i,k}'\bm{H}_i^{-1}(\gamma,\bm\theta)\bm\epsilon_i\big|
			=&~ O_p(n_i^{\ell/2}).
			\end{split}
			%\label{appendix:lemma:epsilon:eq1}
			\end{align*}
			\item[(ii)] For $i=1,\dots,m$ and $k\notin\gamma$,
			\begin{align*}
			\sup_{\bm\theta\in[0,\infty)^{q(\gamma)}}\big|\bm{z}_{i,k}'\bm{H}_i^{-1}(\gamma,\bm\theta)\bm\epsilon_i\big|
			=&~ O_p(n_i^{\ell/2}).
			%\label{appendix:lemma:epsilon:eq2}
			\end{align*}
			\item[(iii)] For $i=1,\dots,m$ and $j=1,\dots,p$,
			\begin{align*}
			\sup_{\bm\theta\in[0,\infty)^{q(\gamma)}}\big|\bm{x}_{i,j}'\bm{H}_i^{-1}(\gamma,\bm\theta)\bm\epsilon_i\big|
			=&~ O_p(n_i^{\xi/2}).
			%\label{appendix:lemma:epsilon:eq3}
			\end{align*}
			\noindent In addition,
			\begin{align*}
			\sup_{\bm\theta\in[0,\infty)^{q(\gamma)}}\bigg|\sum_{i=1}^m\bm{x}_{i,j}'\bm{H}_i^{-1}(\gamma,\bm\theta)\bm\epsilon_i\bigg|
			=&~ O_p\bigg(\bigg(\sum_{i=1}^mn_i^{\xi}\bigg)^{1/2}\bigg).
			\end{align*}
			\item[(iv)] For $i=1,\dots,m$,
			\begin{align*}
			\sup_{\bm\theta\in[0,\infty)^{q(\gamma)}}\bm\epsilon_i'\bm{H}_i^{-1}(\gamma,\bm\theta)\bm\epsilon_i
			=&~ \bm\epsilon_i'\bm\epsilon_i + O_p(q).
			%\label{appendix:lemma:epsilon:eq4}
			\end{align*}
		\end{enumerate}
		\label{appendix:lemma:epsilon}
	\end{lem}
	
	Note that Lemma \ref{appendix:lemma:z x} (i) implies that, for
	$(\alpha,\gamma)\in\mathcal{A}\times\mathcal{G}$,
	\begin{align*}
	\sum_{i=1}^m\bm{X}_{i}(\alpha)'\bm{H}_i^{-1}(\gamma,\bm\theta)\bm{X}_{i}(\alpha)
	=&~ \bigg(\sum_{i=1}^m n_i^{\xi}\bigg)\bm{T}(\alpha)+
	\bigg\{o\bigg(\sum_{i=1}^m n_i^{\xi-\tau}\bigg)\bigg\}_{p(\alpha)\times p(\alpha)}\\
	=&~	\bigg(\sum_{i=1}^m n_i^{\xi}\bigg)\bm{T}(\alpha)+
	\bigg\{o\bigg(n_{\min}^{-\tau}\sum_{i=1}^m n_i^{\xi}\bigg)\bigg\}_{p(\alpha)\times p(\alpha)}
	\end{align*}
	\noindent uniformly over $\bm\theta\in[0,\infty)^{q(\gamma)}$, where
	$\{a\}_{k\times j}$ denotes a $k\times j$ matrix with elements equal
	to $a$ and $\bm{T}(\alpha)$ is a diagonal matrix with diagonal
	elements bounded away from $0$ and $\infty$. Hence by (\ref{eq:inv
		decompose}) with $\bm{c},\bm{d} =
	\{o(n_{\min}^{-\tau/2})\}_{p(\alpha)\times1}$ and
	$\bm{A}=\bm{T}(\alpha)$, we have, for
	$(\alpha,\gamma)\in\mathcal{A}\times\mathcal{G}$,
	\begin{align}
	\bigg(\frac{\sum_{i=1}^m\bm{X}_{i}(\alpha)'\bm{H}_i^{-1}(\gamma,\bm\theta)\bm{X}_{i}(\alpha)}{\sum_{i=1}^m n_i^{\xi}}
	\bigg)^{-1}
	=&~ \bm{T}^{-1}(\alpha) + \{o(n_{\min}^{-\tau})\}_{p(\alpha)\times p(\alpha)}
	\label{proof:lemma:xze:eq0}
	\end{align}
	\noindent uniformly over $\bm\theta\in[0,\infty)^{q(\gamma)}$, which plays a
	key role in proving lemmas for theorems.

	The following lemma shows that $\hat\theta_k$ does not converge to $0$ in
	probability for $k\in\gamma\cap\gamma_0$, which allows us to restrict the
	parameter space of $\bm{\theta}$ from $[0,\infty)^{q(\gamma)}$ to
	\begin{align}
	\Theta_\gamma = \{\bm\theta\in[0,\infty)^{q(\gamma)}:\bm\theta(\gamma\cap\gamma_0)\in(0,\infty)^{q(\gamma\cap\gamma_0)}\}.
	\label{space:theta restricted}
	\end{align}
	\begin{lem}
		Under the assumptions of Theorem~\ref{appendix:theorem:MLE},
		let $\bm\theta_0^\dag$ be $\bm\theta$ except that
		$\{\theta_k:k\in\gamma\cap\gamma_0\}$ are replaced by $\{\theta_{k,0}:k\in\gamma\cap\gamma_0\}$.
		Then for any $(\alpha,\gamma)\in\mathcal{A}\times\mathcal{G}$, $v^2>0$,
		and $\bm\theta\in[0,\infty)^{q(\gamma)}$ with $\theta_k\rightarrow0$ for
		some $k\in\gamma\cap\gamma_0$, we have
		\begin{align*}
		-2\log L(\bm\theta,v^2;\alpha,\gamma)
		-\{-2\log L(\bm\theta_0^\dag,v^2,\alpha,\gamma)\}
		\xrightarrow{p}\infty
		\end{align*}
		\noindent as $N\rightarrow\infty$, where $-2\log
		L(\bm\theta,v^2;\alpha,\gamma)$ is given in (\ref{fn:likelihood}).
		\label{appendix:prop:compact space}
	\end{lem}

	Based on Lemma \ref{appendix:prop:compact space}, the following lemma is
	needed to develop the convergence rates of components of the likelihood
	equations given in \eqref{partial:v} and \eqref{partial:theta2}, uniformly
	over $\Theta_{\gamma}$ defined in \eqref{space:theta restricted}. 
	\begin{lem}
		Consider a mixed-effects model
		$(\alpha,\gamma)\in\mathcal{A}\times\mathcal{G}$ with
		$\bm{H}(\gamma,\bm\theta)$ defined in (\ref{matrix:H}) and
		$\Theta_\gamma$ defined in (\ref{space:theta restricted}). Suppose that
		(A0)--(A3) hold. Then
		\begin{enumerate}
			\item[(i)] For $i,i^*=1,\dots,m$, $(\alpha,\gamma)\in\mathcal{A}\times\mathcal{G}$ and $k, k^*\in\gamma$,
			\begin{align*}
			\begin{split}
			\sup_{\bm\theta\in\Theta_\gamma}\theta_k\theta_{k^*}
			\big|\bm{h}_{i,k}'\bm{H}^{-1}(\gamma,\bm\theta)\bm{M}(\alpha,\gamma;\bm\theta)\bm{h}_{i^*,k^*}\big|
			=&~ o\Bigg(\frac{n_i^{(\xi-\ell)/2}n_{i^*}^{(\xi-\ell)/2-\tau}}{\sum_{i=1}^mn_i^{\xi}}\Bigg),\\
			\sup_{\bm\theta\in\Theta_\gamma}\theta_k
			\big|\bm{h}_{i,k}'\bm{H}^{-1}(\gamma,\bm\theta)\bm{M}(\alpha,\gamma;\bm\theta)\bm{h}_{i^*,k^*}\big|
			=&~ o\Bigg(\frac{n_i^{(\xi-\ell)/2}n_{i^*}^{(\xi+\ell)/2-\tau}}{\sum_{i=1}^mn_i^{\xi}}\Bigg),\\
			\sup_{\bm\theta\in\Theta_\gamma}
			\big|\bm{h}_{i,k}'\bm{H}^{-1}(\gamma,\bm\theta)\bm{M}(\alpha,\gamma;\bm\theta)\bm{h}_{i^*,k^*}\big|
			=&~ o\Bigg(\frac{n_i^{(\xi+\ell)/2}n_{i^*}^{(\xi+\ell)/2-\tau}}{\sum_{i=1}^mn_i^{\xi}}\Bigg).
			\end{split}
			%\label{appendix:lemma:ezx:eq1}
			\end{align*}
			\item[(ii)] For $i,i^*=1,\dots,m$, $(\alpha,\gamma)\in\mathcal{A}\times\mathcal{G}$, $k\in\gamma$ and $k^*\notin\gamma$,
			\begin{align*}
			\begin{split}
			\sup_{\bm\theta\in\Theta_\gamma}\theta_k
			\big|\bm{h}_{i,k}'\bm{H}^{-1}(\gamma,\bm\theta)\bm{M}(\alpha,\gamma;\bm\theta)\bm{h}_{i^*,k^*}\big|
			=&~  o\Bigg(\frac{n_i^{(\xi-\ell)/2}n_{i^*}^{(\xi+\ell)/2-\tau}}{\sum_{i=1}^mn_i^{\xi}}\Bigg),\\
			\sup_{\bm\theta\in\Theta_\gamma}
			\big|\bm{h}_{i,k}'\bm{H}^{-1}(\gamma,\bm\theta)\bm{M}(\alpha,\gamma;\bm\theta)\bm{h}_{i^*,k^*}\big|
			=&~  o\Bigg(\frac{n_i^{(\xi+\ell)/2}n_{i^*}^{(\xi+\ell)/2-\tau}}{\sum_{i=1}^mn_i^{\xi}}\Bigg).
			\end{split}
			%\label{appendix:lemma:ezx:eq2}
			\end{align*}
			\item[(iii)] For $i=1,\dots,m$, $(\alpha,\gamma)\in\mathcal{A}\times\mathcal{G}$ and $k\in\gamma$,
			\begin{align*}
			\begin{split}
			\sup_{\bm\theta\in\Theta_\gamma}\theta_k\big|\bm{h}_{i,k}'\bm{H}^{-1}(\gamma,\bm\theta)\bm{M}(\alpha,\gamma;\bm\theta)\bm\epsilon\big|
			=&~ o_p(n_i^{-\ell/2}),\\
			\sup_{\bm\theta\in\Theta_\gamma}\big|\bm{h}_{i,k}'\bm{H}^{-1}(\gamma,\bm\theta)\bm{M}(\alpha,\gamma;\bm\theta)\bm\epsilon\big|
			=&~ o_p(n_i^{\ell/2}).
			\end{split}
			%\label{appendix:lemma:ezx:eq3}
			\end{align*}
			\item[(iv)] For $i=1,\dots,m$, $(\alpha,\gamma)\in(\mathcal{A}\setminus\mathcal{A}_0)\times\mathcal{G}$ and $k\in\gamma$,
			\begin{align*}
			\begin{split}
			\sup_{\bm\theta\in\Theta_\gamma}\theta_k\big|\bm{h}_{i,k}'\bm{H}^{-1}(\gamma,\bm\theta)\bm{M}(\alpha,\gamma;\bm\theta)\bm{X}(\alpha_0\setminus\alpha)\bm\beta(\alpha_0\setminus\alpha)\big|
			=&~ o(n_i^{(\xi-\ell)/2-\tau}),\\
			\sup_{\bm\theta\in\Theta_\gamma}\big|\bm{h}_{i,k}'\bm{H}^{-1}(\gamma,\bm\theta)\bm{M}(\alpha,\gamma;\bm\theta)\bm{X}(\alpha_0\setminus\alpha)\bm\beta(\alpha_0\setminus\alpha)\big|
			=&~ o(n_i^{(\xi+\ell)/2-\tau}).
			\end{split}
			%\label{appendix:lemma:ezx:eq4}
			\end{align*}
			\item[(v)] For $i=1,\dots,m$ and $(\alpha,\gamma)\in\mathcal{A}\times\mathcal{G}$,
			\begin{align*}
			\sup_{\bm\theta\in\Theta_\gamma}\bm\epsilon'\bm{H}^{-1}(\gamma,\bm\theta)\bm{M}(\alpha,\gamma;\bm\theta)\bm\epsilon
			=&~ O_p(p(\alpha)).
			%\label{appendix:lemma:ezx:eq5}
			\end{align*}
			\item[(vi)] For $i=1,\dots,m$, $(\alpha,\gamma)\in\mathcal{A}\times\mathcal{G}$ and $k\notin\gamma$,
			\begin{align*}
			\sup_{\bm\theta\in\Theta_\gamma}\big|\bm{h}_{i,k}'\bm{H}^{-1}(\gamma,\bm\theta)\bm{M}(\alpha,\gamma;\bm\theta)\bm\epsilon\big|
			=&~ o_p(n_i^{\ell/2}).
			%\label{appendix:lemma:ezx:eq6}
			\end{align*}
			\item[(vii)] For $(\alpha,\gamma)\in(\mathcal{A}\setminus\mathcal{A}_0)\times\mathcal{G}$,
			\begin{align*}
			\begin{split}
			\sup_{\bm\theta\in\Theta_\gamma}\big|\bm\epsilon'\bm{H}^{-1}(\gamma,\bm\theta)\bm{M}(\alpha,\gamma;\bm\theta)\bm{X}(\alpha_0\setminus\alpha)\bm\beta(\alpha_0\setminus\alpha)\big|
			= o_p\bigg(\bigg(\sum_{i=1}^m n_i^{\xi}\bigg)^{1/2}\bigg).
			\end{split}
			%\label{appendix:lemma:ezx:eq6-1}
			\end{align*}
			\item[(viii)] For $i,i^*=1,\dots,m$, $(\alpha,\gamma)\in\mathcal{A}\times\mathcal{G}$ and $k,k^*\notin\gamma$,
			\begin{align*}
			\sup_{\bm\theta\in\Theta_\gamma}
			\big|\bm{h}_{i,k}'\bm{H}^{-1}(\gamma,\bm\theta)\bm{M}(\alpha,\gamma;\bm\theta)\bm{h}_{i^*,k^*}\big|
			=&~  o_p\Bigg(\frac{n_i^{(\xi+\ell)/2}n_{i^*}^{(\xi+\ell)/2-\tau}}{\sum_{i=1}^mn_i^{\xi}}\Bigg).
			%\label{appendix:lemma:ezx:eq7}
			\end{align*}
			\item[(ix)] For $i=1,\dots,m$, $(\alpha,\gamma)\in(\mathcal{A}\setminus\mathcal{A}_0)\times\mathcal{G}$ and $k\notin\gamma$,
			\begin{align*}
			\begin{split}
			\sup_{\bm\theta\in\Theta_\gamma}\big|\bm{h}_{i,k}'\bm{H}^{-1}(\gamma,\bm\theta)
			\bm{M}(\alpha,\gamma;\bm\theta)\bm{X}(\alpha_0\setminus\alpha)\bm\beta(\alpha_0\setminus\alpha)\big|
			=o(n_i^{(\xi+\ell)/2-\tau}).
			\end{split}
			%\label{appendix:lemma:ezx:eq8}
			\end{align*}
			\item[(x)] For $(\alpha,\gamma)\in(\mathcal{A}\setminus\mathcal{A}_0)\times\mathcal{G}$,
			\begin{align*}
			\begin{split}
			\sup_{\bm\theta\in\Theta_\gamma}\big|
			&\bm\beta(\alpha_0\setminus\alpha)'\bm{X}(\alpha_0\setminus\alpha)'\bm{H}^{-1}(\gamma,\bm\theta)\\
			&~  \times\bm{M}(\alpha,\gamma;\bm\theta)\bm{X}(\alpha_0\setminus\alpha)\bm\beta(\alpha_0\setminus\alpha)\big|
			=o\bigg(\sum_{i=1}^mn_i^{\xi-\tau}\bigg).
			\end{split}
			%\label{appendix:lemma:ezx:eq9}
			\end{align*}
		\end{enumerate}
		\label{appendix:lemma:xze}
	\end{lem}

	\setcounter{equation}{0}
	\section{Theoretical Proofs}
	\label{appendix:proofs}

	\subsection{Proof of Theorem \ref{appendix:theorem:MLE}}
	\label{section:proofs}
	
	We shall focus on the asymptotic properties of $\hat{v}^2(\alpha,\gamma)$
	and $\big\{\hat{\theta}_k(\alpha,\gamma):k\in\gamma\big\}$, and derive the
	asymptotic properties of $\{\hat{\sigma}_k^2(\alpha,\gamma):k\in\gamma\}$
	via
	$\hat{\sigma}_k^2(\alpha,\gamma)=\hat{v}^2(\alpha,\gamma)\hat{\theta}_k(\alpha,\gamma)$;
	$k\in\gamma$. 
	If $\hat{v}^2(\alpha,\gamma)>0$ and $\hat{\theta}_k(\alpha,\gamma)>0$; $k\in\gamma$,
	then we can derive them using the likelihood equations.
	Differentiating the profile log-likelihood function of (\ref{fn:likelihood})
	with respect to $v^2$ and $\{\theta_k:k\in\gamma\}$, we obtain
	\begin{equation}
	\frac{\partial}{\partial v^2}\{-2\log L(\bm\theta,v^2;\alpha,\gamma)\}
	=\frac{N}{v^2} - \frac{
		\bm{y}'\bm{H}^{-1}(\gamma,\bm\theta)(\bm{I}_{N}-\bm{M}(\alpha,\gamma;\bm\theta))\bm{y}}{v^4}
	\label{partial:v}
	\end{equation}
	\noindent and
	\begin{align}
	\begin{split}
	\frac{\partial}{\partial \theta_k}\{-2\log L(\bm\theta,v^2;\alpha,&\gamma)\}
	=\sum_{i=1}^m\bigg\{\bm{z}_{i,k}'\bm{H}_i^{-1}(\gamma,\bm\theta)\bm{z}_{i,k}\\
	&~   -\frac{\{\bm{h}_{i,k}'\bm{H}^{-1}(\gamma,\bm\theta)(\bm{I}_{N}-\bm{M}(\alpha,\gamma;\bm\theta))\bm{y}\}^2}{v^2}\bigg\}.
	\end{split}
	\label{partial:theta2}
	\end{align}
	
	\noindent
	To derive $\hat{v}^2(\alpha,\gamma)$ and
	$\big\{\hat{\theta}_k(\alpha,\gamma):k\in\gamma\big\}$,
	we must study the convergence rate of each term on the right-hand sides of
	both \eqref{partial:v} and \eqref{partial:theta2} by Lemmas
	\ref{appendix:lemma:z x}--\ref{appendix:lemma:epsilon} and Lemma
	\ref{appendix:lemma:xze}.
	
	We first prove (\ref{appendix:thm:mle:correct:eq1}) using (\ref{partial:v}).
	Consider the following decomposition of
	$\bm{y}'\bm{H}^{-1}(\gamma,\bm\theta)(\bm{I}_{N}-\bm{M}(\alpha,\gamma;\bm\theta))\bm{y}$
	in (\ref{partial:v}):
	\begin{align}
	\begin{split}
	\bm{y}'\bm{H}^{-1}&(\gamma,\bm\theta)(\bm{I}_{N}-\bm{M}(\alpha,\gamma;\bm\theta))\bm{y}\\
	=&~ \bm{\mu}'_0\bm{H}^{-1}(\gamma,\bm\theta)(\bm{I}_{N}-\bm{M}(\alpha,\gamma;\bm\theta))\bm{\mu}_0\\
	&~ +2\bm{\mu}'_0\bm{H}^{-1}(\gamma,\bm\theta)(\bm{I}_{N}-\bm{M}(\alpha,\gamma;\bm\theta))(\bm{Z}(\gamma_0)\bm{b}(\gamma_0)+\bm\epsilon)\\
	&~ +(\bm{Z}(\gamma_0)\bm{b}(\gamma_0)+\bm\epsilon)'\bm{H}^{-1}(\gamma,\bm\theta)(\bm{Z}(\gamma_0)\bm{b}(\gamma_0)+\bm\epsilon)\\
	&~ -(\bm{Z}(\gamma_0)\bm{b}(\gamma_0)+\bm\epsilon)'\bm{H}^{-1}(\gamma,\bm\theta)\bm{M}(\alpha,\gamma;\bm\theta)(\bm{Z}(\gamma_0)\bm{b}(\gamma_0)+\bm\epsilon).
	\end{split}
	\label{eq:decomposition of v}
	\end{align}
	
	\noindent The first two terms of \eqref{eq:decomposition of v} are zeros
	because
	\begin{align}
	(\bm{I}_{N}-\bm{M}(\alpha,\gamma;\bm\theta))\bm\mu_0 = \bm{0};\quad\alpha\in\mathcal{A}_0.
	\label{proof:prop:compact space:eq2}
	\end{align}
	
	\noindent By Lemma \ref{appendix:lemma:z} (ii)--(iii), Lemma
	\ref{appendix:lemma:epsilon} (i), and Lemma \ref{appendix:lemma:epsilon}
	(iv), the third term of \eqref{eq:decomposition of v} can be written as
	\begin{align*}
	\sum_{i=1}^m(\bm{Z}_i(\gamma_0)&\bm{b}_i(\gamma_0)+\bm\epsilon_i)'\bm{H}_i^{-1}(\gamma,\bm\theta)(\bm{Z}_i(\gamma_0)\bm{b}_i(\gamma_0)+\bm\epsilon_i)\\
	=&~ \sum_{i=1}^m\bm\epsilon_i'\bm\epsilon_i +O_p\bigg(\sum_{k\in\gamma_0}\frac{m}{\theta_k}\bigg)
	+o_p\bigg(\sum_{k\in\gamma_0}\frac{m}{\theta_k^2}\bigg)+O_p(mq)
	\end{align*}
	
	\noindent uniformly over $\bm\theta\in\Theta_\gamma$.
	Note that by the Cauchy--Schwarz inequality, 
	\begin{align}
	\bigg(\sum_{i=1}^m n_i^{(\xi-\ell)/2}\bigg)^2=O\bigg(\sum_{i=1}^mn_i^{\xi}\sum_{i^*=1}^m n_{i^*}^{-\ell}\bigg).
	\label{cauchy}
	\end{align}
	\noindent Hence, by Lemma~\ref{appendix:lemma:xze}~(i), 
	Lemma~\ref{appendix:lemma:xze}~(iii), and Lemma~\ref{appendix:lemma:xze}~(v), the
	last term of \eqref{eq:decomposition of v} can be written as
	\begin{align*}
	\bigg\{\bigg(\sum_{i=1}^m&\sum_{k\in\gamma_0}b_{i,k}\bm{h}_{i,k}\bigg)+\bm\epsilon\bigg\}'
	\bm{H}^{-1}(\gamma,\bm\theta)\bm{M}(\alpha,\gamma;\bm\theta)
	\bigg\{\bigg(\sum_{i=1}^m\sum_{k\in\gamma_0}b_{i,k}\bm{h}_{i,k}\bigg)+\bm\epsilon\bigg\}\\
	=&~ o_p\bigg(\sum_{k,k^*\in\gamma_0}\frac{m}{\theta_k\theta_{k^*}}\bigg)
	+o_p\bigg(\sum_{k\in\gamma_0}\frac{m}{\theta_k}\bigg)
	+O_p(p+mq)
	\end{align*}
	\noindent uniformly over $\bm\theta\in\Theta_\gamma$.
	Therefore, we can rewrite \eqref{eq:decomposition of v} as
	\begin{align*}
	\bm{y}'\bm{H}^{-1}(\gamma,\bm\theta)&(\bm{I}_{N}-\bm{M}(\alpha,\gamma;\bm\theta))\bm{y}\\
	=&~ \bm\epsilon'\bm\epsilon
	+o_p\bigg(\sum_{k,k^*\in\gamma_0}\frac{m}{\theta_k\theta_{k^*}}\bigg)
	+O_p\bigg(\sum_{k\in\gamma_0}\frac{m}{\theta_k}\bigg)
	+O_p(p+mq).
	\end{align*}
	
	\noindent It follows from (\ref{partial:v}) that for $v^2\in(0,\infty)$,
	\begin{align*}
	\begin{split}
	v^4\bigg\{\frac{\partial}{\partial v^2}\{-2\log L(\bm\theta,v^2;\alpha,\gamma)\}\bigg\}
	=&~ N\bigg(v^2-\frac{\bm\epsilon'\bm\epsilon}{N}\bigg)+o_p\bigg(\sum_{k,k^*\in\gamma_0}\frac{m}{\theta_k\theta_{k^*}}\bigg)\\
	&~  +O_p\bigg(\sum_{k\in\gamma_0}\frac{m}{\theta_k}\bigg)
	+O_p(p+mq)
	\end{split}
	\end{align*}
	
	\noindent uniformly over $\bm\theta\in\Theta_\gamma$.
	This and Lemma \ref{appendix:prop:compact space} imply that
	\begin{align}
	\begin{split}
	\hat{v}^2(\alpha,\gamma)
	=&~ \frac{\bm\epsilon'\bm\epsilon}{N}
	+O_p\Big(\frac{p+mq}{N}\Big).
	\end{split}
	\label{eq:ML of v}
	\end{align}
	
	\noindent Thus (\ref{appendix:thm:mle:correct:eq1})
	follows by applying the law of large numbers to
	$\bm{\epsilon}'\bm{\epsilon}/N$.
	In addition, the asymptotic normality of $\hat{v}^2(\alpha,\gamma)$ follows
	by $p+mq=o(N^{1/2})$ and an application of the central limit theorem to
	$\bm{\epsilon}'\bm{\epsilon}/N$ in \eqref{eq:ML of v}.
	
	Next, we prove \eqref{appendix:thm:mle:correct:eq2}, for
	$k\in\gamma\cap\gamma_0$, using \eqref{partial:theta2}.
	By Lemma \ref{appendix:lemma:xze} (i) and Lemma \ref{appendix:lemma:xze}
	(iii), we have, for $k\in\gamma\cap\gamma_0$,
	\begin{align}
	\begin{split}
	\theta_{k}\bm{h}_{i,k}'&\bm{H}^{-1}(\gamma,\bm\theta)\bm{M}(\alpha,\gamma;\bm\theta)
	(\bm{Z}(\gamma_0)\bm{b}(\gamma_0)+\bm\epsilon)\\
	=&~
	\theta_{k}\bm{h}_{i,k}'\bm{H}^{-1}(\gamma,\bm\theta)\bm{M}(\alpha,\gamma;\bm\theta)
	\bigg\{\bigg(\sum_{i^*=1}^m\sum_{k^*\in\gamma_0}b_{i^*,k^*}\bm{h}_{i^*,k^*}\bigg)+\bm\epsilon\bigg\}\\
	=&~ o_p\bigg(\sum_{k^*\in\gamma_0}\frac{n_i^{(\xi-\ell)/2}\sum_{i^*=1}^mn_{i^*}^{(\xi-\ell)/2}}{\theta_{k^*}\sum_{i=1}^m n_i^{\xi}}\bigg)+o_p(n_i^{-\ell/2})
	\end{split}
	\label{proof:thm2:eq01}
	\end{align}
	
	\noindent uniformly over $\bm\theta\in\Theta_\gamma$.
	This and (\ref{proof:prop:compact space:eq2}) imply that for
	$k\in\gamma\cap\gamma_0$,
	\begin{align*}
	\theta_k&\bm{h}_{i,k}'\bm{H}^{-1}(\gamma,\bm\theta)(\bm{I}_{N}-\bm{M}(\alpha,\gamma;\bm\theta))\bm{y}\\
	=&~\theta_k\bm{h}_{i,k}'\bm{H}^{-1}(\gamma,\bm\theta)(\bm{I}_{N}-\bm{M}(\alpha,\gamma;\bm\theta))
	(\bm{Z}(\gamma_0)\bm{b}(\gamma_0)+\bm\epsilon)\\
	=&~ \theta_k\bm{h}_{i,k}'\bm{H}^{-1}(\gamma,\bm\theta)
	(\bm{Z}(\gamma_0)\bm{b}(\gamma_0)+\bm\epsilon)\\
	&~    +o_p\bigg(\sum_{k^*\in\gamma_0}\frac{n_i^{(\xi-\ell)/2}\sum_{i^*=1}^mn_{i^*}^{(\xi-\ell)/2}}{\theta_{k^*}\sum_{i=1}^m n_i^{\xi}}\bigg)
	+o_p(n_i^{-\ell/2})\\
	=&~ \theta_k\bm{z}_{i,k}'\bm{H}_i^{-1}(\gamma,\bm\theta)(\bm{Z}_i(\gamma_0)\bm{b}_i(\gamma_0)+\bm{\epsilon}_i)\\
	&~    +o_p\bigg(\sum_{k^*\in\gamma_0}\frac{n_i^{(\xi-\ell)/2}\sum_{i^*=1}^mn_{i^*}^{(\xi-\ell)/2}}{\theta_{k^*}\sum_{i=1}^m n_i^{\xi}}\bigg)+o_p(n_i^{-\ell/2})\\
	=&~ b_{i,k}+O_p(n_i^{-\ell/2})
	+o_p\bigg(\sum_{k^*\in\gamma_0}\frac{n_i^{(\xi-\ell)/2}\sum_{i^*=1}^mn_{i^*}^{(\xi-\ell)/2}}{\theta_{k^*}\sum_{i=1}^m n_i^{\xi}}\bigg)
	\end{align*}
	
	\noindent uniformly over $\bm\theta\in\Theta_\gamma$, where the last
	equality follows from Lemma \ref{appendix:lemma:z} (ii)--(iii) and Lemma
	\ref{appendix:lemma:epsilon} (i).
	Hence, for $k\in\gamma\cap\gamma_0$,
	\begin{align*}
	\theta_k^2&\{\bm{h}_{i,k}'\bm{H}^{-1}(\gamma,\bm\theta)(\bm{I}_{N}-\bm{M}(\alpha,\gamma;\bm\theta))\bm{y}\}^2\\
	=&~ b_{i,k}^2
	+O_p(n_i^{-\ell/2})
	+o_p\bigg(\sum_{k^*\in\gamma_0}\frac{n_i^{(\xi-\ell)/2}\sum_{i^*=1}^mn_{i^*}^{(\xi-\ell)/2}}{\theta_{k^*}\sum_{i=1}^m n_i^{\xi}}\bigg)
	\end{align*}
	
	\noindent uniformly over $\bm\theta\in\Theta_\gamma$.
	This together with Lemma \ref{appendix:lemma:z} (ii) and
	(\ref{partial:theta2}) imply that for $k\in\gamma\cap\gamma_0$,
	\begin{align}
	\begin{split}
	\theta_k^2& \bigg\{\frac{\partial}{\partial \theta_k}
	\{-2\log L(\bm\theta,v^2;\alpha,\gamma)\}\bigg\}\\
	=&~ m\bigg(\theta_k -\frac{1}{m}\sum_{i=1}^m\frac{b_{i,k}^2}{v^2} \bigg)+O_p\bigg(\sum_{i=1}^mn_i^{-\ell/2}\bigg)\\
	&~  +o_p\bigg(\sum_{k^*\in\gamma_0}\frac{\sum_{i=1}^mn_i^{(\xi-\ell)/2}\sum_{i^*=1}^mn_{i^*}^{(\xi-\ell)/2-\tau}}{\theta_{k^*}\sum_{i=1}^m n_i^{\xi}}\bigg)
	\end{split}
	\label{exist of theta k 0}
	\end{align}
	
	\noindent uniformly over $\bm\theta\in\Theta_\gamma$.
	By \eqref{cauchy}, Lemma \ref{appendix:prop:compact space} and setting
	\eqref{exist of theta k 0} to $0$, we obtain
	\[
	\hat\theta_k(\alpha,\gamma)
	=\frac{1}{m}\sum_{i=1}^m\frac{b_{i,k}^2}{\hat{v}^2(\alpha,\gamma)}
	+O_p\bigg(\frac{1}{m}\sum_{i=1}^mn_i^{-\ell/2}\bigg);\quad k\in\gamma\cap\gamma_0.
	\]
	This proves (\ref{appendix:thm:mle:correct:eq2}), for $k\in\gamma\cap\gamma_0$.
	
	It remains to prove (\ref{appendix:thm:mle:correct:eq2}), for
	$k\in\gamma\setminus\gamma_0$.
	We prove by showing that \eqref{partial:theta2} is asymptotically
	nonnegative, for $\theta_k\in\big(n_{\max}^{-\ell},\infty\big)$;
	$k\in\gamma\setminus\gamma_0$
	using a recursive argument. Let $\bm\theta^\dag$ be $\bm\theta$ except that
	$\{\theta_k: k\in\gamma\cap\gamma_0\}$ are replaced by
	$\{\hat{\theta}_k(\alpha,\gamma):k\in\gamma\cap\gamma_0\}$. By Lemma
	\ref{appendix:lemma:xze} (i) and Lemma \ref{appendix:lemma:xze} (iii), we
	have, for $k\in\gamma\setminus\gamma_0$,
	\begin{align*}
	\theta_k\bm{h}_{i,k}'
	& \bm{H}^{-1}(\gamma,\bm\theta^\dag)\bm{M}(\alpha,\gamma;\bm\theta^\dag)(\bm{Z}(\gamma_0)\bm{b}(\gamma_0)+\bm\epsilon)\\
	=&~\theta_k\bm{h}_{i,k}'\bm{H}^{-1}(\gamma,\bm\theta^\dag)\bm{M}(\alpha,\gamma;\bm\theta^\dag)
	\bigg(\sum_{i^*=1}^m\sum_{k^*\in\gamma_0}b_{i^*,k^*}\bm{h}_{i^*,k^*}+\bm\epsilon\bigg)\\
	=&~o_p\bigg(\frac{n_i^{(\xi-\ell)/2}\sum_{i^*=1}^mn_{i^*}^{(\xi-\ell)/2}}{\sum_{i=1}^m n_i^{\xi}}\bigg)+o_p(n_i^{-\ell/2})
	\end{align*}
	
	\noindent uniformly over $\bm{\theta}(\gamma\setminus\gamma_0) \in[0,\infty)^{q(\gamma\setminus\gamma_0)}$.
	This and (\ref{proof:prop:compact space:eq2}) imply that for
	$k\in\gamma\setminus\gamma_0$,
	\begin{align}
	\begin{split}
	\theta_{k}&\bm{h}_{i,k}'\bm{H}^{-1}(\gamma,\bm\theta^\dag)(\bm{I}_{N}-\bm{M}(\alpha,\gamma;\bm{\theta}^\dag ))\bm{y}\\
	=&~ \theta_{k}\bm{h}_{i,k}'\bm{H}^{-1}(\gamma,\bm{\theta}^\dag )(\bm{I}_{N}-\bm{M}(\alpha,\gamma;\bm{\theta}^\dag ))
	(\bm{Z}(\gamma_0)\bm{b}(\gamma_0)+\bm\epsilon)\\
	=&~\theta_{k}\bm{h}_{i,k}'\bm{H}^{-1}(\gamma,\bm{\theta}^\dag )(\bm{Z}(\gamma_0)\bm{b}(\gamma_0)+\bm\epsilon)\\
	&~
	+o_p\bigg(\frac{n_i^{(\xi-\ell)/2}\sum_{i^*=1}^mn_{i^*}^{(\xi-\ell)/2}}{\sum_{i=1}^m n_i^{\xi}}\bigg)+o_p(n_i^{-\ell/2})\\
	=&~ \theta_k\bm{z}_{i,k}'\bm{H}_i^{-1}(\gamma,\bm{\theta}^\dag )
	\bigg(\sum_{k^*\in\gamma_0}\bm{z}_{i,k^*}b_{i,k^*}+\bm\epsilon_i\bigg)\\
	&~+o_p\bigg(\frac{n_i^{(\xi-\ell)/2}\sum_{i^*=1}^mn_{i^*}^{(\xi-\ell)/2}}{\sum_{i=1}^m n_i^{\xi}}\bigg)+o_p(n_i^{-\ell/2})\\
	=&~	O_p(n_i^{-\ell/2})+o_p\bigg(\frac{n_i^{(\xi-\ell)/2}\sum_{i^*=1}^mn_{i^*}^{(\xi-\ell)/2}}{\sum_{i=1}^m n_i^{\xi}}\bigg)
	\end{split}
	\label{proof:thm:mle part1:eq15-0}
	\end{align}
	
	\noindent uniformly over $\bm{\theta}(\gamma\setminus\gamma_0)
	\in[0,\infty)^{q(\gamma\setminus\gamma_0)}$, where the last equality follows
	from Lemma \ref{appendix:lemma:z} (iii) and Lemma
	\ref{appendix:lemma:epsilon} (i).
	Hence by \eqref{cauchy}, Lemma \ref{appendix:lemma:z} (ii), and
	(\ref{partial:theta2}), we have, for
	$\bm\theta(\gamma\setminus\gamma_0)\in[0,\infty)^{q(\gamma\setminus\gamma_0)}$
	and $k\in\gamma\setminus\gamma_0$,
	\begin{align*}
	\begin{split}
	\theta_k^2&\bigg\{\frac{\partial}{\partial \theta_k}\{-2\log L(\bm\theta^\dag,v^2;\alpha,\gamma)\}\bigg\}\\
	=&~ m\theta_k +O_p\bigg(\sum_{i=1}^m n_i^{-\ell}\bigg)
	+o_p\bigg(\frac{\sum_{i=1}^mn_i^{(\xi-\ell)/2}\sum_{i^*=1}^mn_{i^*}^{(\xi-\ell)/2-\tau}}{\sum_{i=1}^m n_i^{\xi}}\bigg)\\
	=&~	m\theta_k +O_p\bigg(\sum_{i=1}^m n_i^{-\ell}\bigg)\\
	=&~ m\theta_k+ o_p(m\log(n_{\min})n_{\min}^{-\ell}).
	\end{split}
	\end{align*}
	\noindent This implies that $-2\log L(\bm\theta^\dag,v^2;\alpha,\gamma)$
	is an asymptotically nondecreasing function on
	$\theta_k\in(\log(n_{\min})n_{\min}^{-\ell},\infty)$, for
	$k\in\gamma\setminus\gamma_0$ given other
	$\bm\theta(\gamma\setminus\gamma_0)\in[0,\infty)^{q(\gamma\setminus\gamma_0)}$.
	It follows that $\hat\theta_k(\alpha,\gamma)
	\in\big[0,\log(n_{\min})n_{\min}^{-\ell}\big)$;
	$k\in\gamma\setminus\gamma_0$. The above convergence rate can be recursively
	improved.
	Without loss of generality, assume that $n_{\min}=n_{1}\leq n_{2}\leq
	\cdots\leq n_{m}=n_{\max}$.
	We can restrict the parameter space of $\theta_k$ in the next step to
	\begin{align}
	\Theta_{\gamma,k,i}=\big\{\bm{\theta}(\gamma\setminus\gamma_0)
	\in[0,\infty)^{q(\gamma\setminus\gamma_0)}:\theta_k\leq \log(n_{\min})n_i^{-\ell}\big\}
	\label{space:theta dag k}
	\end{align}
	
	\noindent with $i=1$.
	Then, by Lemma \ref{appendix:lemma:xze} (i) and Lemma
	\ref{appendix:lemma:xze} (iii), we have, for $k\in\gamma\setminus\gamma_0$,
	\begin{align*}
	\theta_k&\bm{h}_{1,k}'\bm{H}^{-1}(\gamma,\bm{\theta}^\dag )\bm{M}(\alpha,\gamma;\bm{\theta}^\dag )(\bm{Z}(\gamma_0)\bm{b}(\gamma_0)+\bm\epsilon)\\
	=&~ o_p\bigg(\frac{n_1^{(\xi-\ell)/2}\sum_{i^*=1}^mn_{i^*}^{(\xi-\ell)/2}}{\sum_{i=1}^m n_i^{\xi}}\bigg)+
	o_p(\theta_kn_1^{\ell/2})
	\end{align*}
	uniformly over $\bm{\theta}(\gamma\setminus\gamma_0)\in\Theta_{\gamma,k,1}$.
	This and (\ref{proof:prop:compact space:eq2}) imply that for
	$k\in\gamma\setminus\gamma_0$,
	\begin{align*}
	\begin{split}
	\theta_{k}\bm{h}_{1,k}'\bm{H}^{-1}&(\gamma,\bm{\theta}^\dag )(\bm{I}_{N}-\bm{M}(\alpha,\gamma;\bm{\theta}^\dag ))\bm{y}\\
	=&~ \theta_{k}\bm{h}_{1,k}'\bm{H}^{-1}(\gamma,\bm{\theta}^\dag )(\bm{I}_{N}-\bm{M}(\alpha,\gamma;\bm{\theta}^\dag ))
	(\bm{Z}(\gamma_0)\bm{b}(\gamma_0)+\bm\epsilon)\\
	=&~\theta_{k}\bm{h}_{1,k}'\bm{H}^{-1}(\gamma,\bm{\theta}^\dag )(\bm{Z}(\gamma_0)\bm{b}(\gamma_0)+\bm\epsilon)\\
	&~+o_p\bigg(\frac{n_1^{(\xi-\ell)/2}\sum_{i^*=1}^mn_{i^*}^{(\xi-\ell)/2}}{\sum_{i=1}^m n_i^{\xi}}\bigg)+
	o_p(\theta_kn_1^{\ell/2})\\
	=&~\theta_k\bm{z}_{1,k}'\bm{H}_1^{-1}(\gamma,\bm{\theta}^\dag )
	(\bm{Z}_1(\gamma_0)\bm{b}_1(\gamma_0)+\bm\epsilon_1)\\
	&~+o_p\bigg(\frac{n_1^{(\xi-\ell)/2}\sum_{i^*=1}^mn_{i^*}^{(\xi-\ell)/2}}{\sum_{i=1}^m n_i^{\xi}}\bigg)+
	o_p(\theta_kn_1^{\ell/2})\\
	=&~O_p(\theta_k n_1^{\ell/2})+o_p\bigg(\frac{n_1^{(\xi-\ell)/2}\sum_{i^*=1}^mn_{i^*}^{(\xi-\ell)/2}}{\sum_{i=1}^m n_i^{\xi}}\bigg)
	\end{split}
	\end{align*}
	
	\noindent uniformly over
	$\bm{\theta}(\gamma\setminus\gamma_0)\in\Theta_{\gamma,k,1}$, where the last
	equality follows from Lemma~\ref{appendix:lemma:z}~(iii) and Lemma~\ref{appendix:lemma:epsilon}~(i).
	Hence by \eqref{cauchy}, Lemma~\ref{appendix:lemma:z}~(ii), \eqref{partial:theta2},
	and \eqref{proof:thm:mle part1:eq15-0}, we have 
	\begin{align*}
	\theta_k^2&\bigg\{\frac{\partial}{\partial \theta_k}\{-2\log L(\bm\theta^\dag,v^2;\alpha,\gamma)\}\bigg\}\\
	=&~	(m-1)\theta_k+O_p(\theta_k^2n_1^{\ell})+
	O_p\bigg(\sum_{i=2}^mn_i^{-\ell}\bigg)+o_p\bigg(\sum_{i=1}^m n_i^{-\ell}\bigg)\\
	=&~ (m-1)\theta_k+O_p(\log(n_{\min})\theta_k)+
	O_p\bigg(\sum_{i=2}^mn_i^{-\ell}\bigg)
	\end{align*}
	
	\noindent uniformly over
	$\bm{\theta}(\gamma\setminus\gamma_0)\in\Theta_{\gamma,k,1}$. Hence, setting
	the above equation equal to $0$, we have
	\[
	\hat\theta_k(\alpha,\gamma)=
	\frac{1}{m-1+O_p(\log(n_{\min}))}O_p\bigg(\sum_{i=2}^m n_i^{-\ell}\bigg)
	=O_p(n_{2}^{-\ell}).
	\]
	Now we can further restrict the parameter space of $\theta_k$ to
	$\Theta_{\gamma,k,2}$ in \eqref{space:theta dag k}.
	Continuing this procedure, we can recursively obtain
	$\hat\theta_k(\alpha,\gamma)=O_p(n_i^{-\ell})$;
	$k\in\gamma\setminus\gamma_0$, for $i=3,\dots,m$. This completes the proof
	of  \eqref{appendix:thm:mle:correct:eq2}, for $k\in\gamma\setminus\gamma_0$.
	Hence the proof of Theorem~\ref{appendix:theorem:MLE} is complete.
	
		\subsection{Proof of Example \ref{prop:prediction}}
			Note that for $q=1$, $\bm{Z}_i=\bm{z}_{i,1}$ and $\bm{b}_i=b_{i,1}$.
			Note that by Lemma \ref{appendix:prop:compact space}, we consider the
			sample space $(\sigma_1^2,v^2)\in(0,\infty)^2$.   
			We first derive the explicit forms of the ML estimators $\hat\theta_1$
			and $\hat{v}^2$. 
		By \eqref{partial:theta2}, we have
		\begin{align*}
		\frac{\partial}{\partial \theta_1}\{-2\log L(\theta_1,v^2)\}
		=&~	\sum_{i=1}^m 
		\frac{\bm{z}_{i,1}'\bm{z}_{i,1}}{1+\theta_1\bm{z}_{i,1}'\bm{z}_{i,1}}
		-\frac{1}{v^2}\sum_{i=1}^m
		\bigg\{\bm{z}_{i,1}'\bigg(
		\bm{I}_n - \frac{\theta_1\bm{z}_{i,1}\bm{z}_{i,1}'}{1+\theta_1\bm{z}_{i,1}'\bm{z}_{i,1}}
		\bigg)\bm{y}_i\bigg\} ^2\\
		=&~	\sum_{i=1}^m \frac{\bm{z}_{i,1}'\bm{z}_{i,1}}{1+\theta_1\bm{z}_{i,1}'\bm{z}_{i,1}}
		-\frac{1}{v^2}
		\sum_{i=1}^m\bigg\{
		\frac{\bm{z}_{i,1}'\bm{z}_{i,1}b_{i,1}}{1+\theta_1\bm{z}_{i,1}'\bm{z}_{i,1}}
		+\frac{\bm{z}_{i,1}'\bm\epsilon_i}{1+\theta_1\bm{z}_{i,1}'\bm{z}_{i,1}}
		\bigg\}^2\\
		=&~ \sum_{i=1}^m\bigg(\frac{1}{\theta_1}-\frac{1}{\theta_1(1+\theta_1\bm{z}_{i,1}'\bm{z}_{i,1})}\bigg)\\
		&~
		-\frac{1}{v^2}
		\sum_{i=1}^m
		\bigg\{
		\frac{b_{i,1}}{\theta_1}
		-\frac{b_{i,1}}{\theta_1(1+\theta_1\bm{z}_{i,1}'\bm{z}_{i,1})}
		+\frac{\bm{z}_{i,1}'\bm\epsilon_i}{1+\theta_1\bm{z}_{i,1}'\bm{z}_{i,1}}
		\bigg\}^2\\
		=&~	\frac{m}{\theta_1}
		-\frac{\sum_{i=1}^mb_{i,1}^2}{v^2\theta_1^2}
		+2\sum_{i=1}^m\frac{b_{i,1}\bm{z}_{i,1}'\bm\epsilon_i}
		{v^2\theta_1(1+\theta_1\bm{z}_{i,1}'\bm{z}_{i,1})}
		+R(\sigma_1^2,v^2),
		\end{align*}
		\noindent where $\sigma_1^2=\theta_1v^2$ and 
		\begin{align}
		\begin{split}
		R(\sigma_1^2,v^2)
		=&~	-\sum_{i=1}^m \frac{1}{\theta_1(1+\theta_1\bm{z}_{i,1}'\bm{z}_{i,1})}-\sum_{i=1}^m\frac{(\bm{z}_{i,1}'\bm\epsilon_i)^2}{v^2\{1+\theta_1\bm{z}_{i,1}'\bm{z}_{i,1}\}^2}
		+\sum_{i=1}^m
		\frac{2b_{i,1}\bm{z}_{i,1}'\bm\epsilon_i}{v^2\theta_1\{1+\theta_1\bm{z}_{i,1}'\bm{z}_{i,1}\}^2}\\
		&~+\sum_{i=1}^m\frac{2b_{i,1}^2}{v^2\theta_1^2(1+\theta_1\bm{z}_{i,1}'\bm{z}_{i,1})}
		-\sum_{i=1}^m \frac{b_{i,1}^2}{v^2\theta_1^2\{1+\theta_1\bm{z}_{i,1}'\bm{z}_{i,1}\}^2}.
		\end{split}
		\label{eq:R}
		\end{align}
		\noindent Note that ML estimators $\hat\sigma_1^2=\hat\theta_1\hat{v}^2$
		and $\hat{v}^2$ satisfy
		\begin{align*}
		0=&~	\frac{m}{\hat\theta_1}
		-\frac{\sum_{i=1}^mb_{i,1}^2}{\hat{v}^2\hat\theta_1^2}
		+\sum_{i=1}^m\frac{2b_{i,1}\bm{z}_{i,1}'\bm\epsilon_i}{\hat{v}^2\hat\theta_1(1+\hat\theta_1(\bm{z}_{i,1}'\bm{z}_{i,1}))}+R(\hat\sigma_1^2,\hat{v}^2),
		\end{align*}
		\noindent which implies that
		\begin{align}
		\begin{split}
		\hat\sigma_1^2
		=&~\hat\theta_1\hat{v}^2
		=\frac{1}{m}\sum_{i=1}^mb_{i,1}^2
		+\frac{1}{m}\sum_{i=1}^m\frac{2\hat\theta_1b_{i,1}\bm{z}_{i,1}'\bm\epsilon_i}{1+\hat\theta_1\bm{z}_{i,1}'\bm{z}_{i,1}}
		+\frac{\hat\theta_1^2}{m}R(\hat\sigma_1^2,\hat{v}^2)\\
		=&~
		\frac{1}{m}\sum_{i=1}^mb_{i,1}^2
		+\frac{1}{m}\sum_{i=1}^m\frac{2b_{i,1}\bm{z}_{i,1}'\bm\epsilon_i}{\bm{z}_{i,1}'\bm{z}_{i,1}}
		-\frac{1}{m}\sum_{i=1}^m\frac{2b_{i,1}\bm{z}_{i,1}'\bm\epsilon_i}{(1+\hat\theta_1\bm{z}_{i,1}'\bm{z}_{i,1})\bm{z}_{i,1}'\bm{z}_{i,1}}
		+\frac{\hat\theta_1^2}{m}R(\hat\sigma_1^2,\hat{v}^2)\\	
		=&~	\frac{1}{m}\sum_{i=1}^mb_{i,1}^2
		+\frac{1}{m}\sum_{i=1}^m\frac{2b_{i,1}\bm{z}_{i,1}'\bm\epsilon_i}{\bm{z}_{i,1}'\bm{z}_{i,1}}
		+R^*(\hat\sigma_1^2,\hat{v^2}),
		\end{split}
		\label{proof:exm:mle1}
		\end{align}
		\noindent where 
		\begin{align}
		R^*(\hat\sigma_1^2,\hat{v}^2)
		=&~	-\frac{1}{m}\sum_{i=1}^m\frac{2b_{i,1}\bm{z}_{i,1}'\bm\epsilon_i}{(1+\hat\theta_1\bm{z}_{i,1}'\bm{z}_{i,1})\bm{z}_{i,1}'\bm{z}_{i,1}}
		+\frac{\hat\theta_1^2}{m}R(\hat\theta_1,\hat{v}^2)
		\label{eq:R*}
		\end{align}
		\noindent with $R(\sigma_1^2,v^2)$ defined in \eqref{eq:R}. By
		\eqref{proof:exm:mle1}, we have 
		\begin{align}
		\begin{split}
		\sum_{i=1}^m b_{i,1}^2
		=&~ O_p(\hat\sigma_1^2),\\
		b_{i,1}^2 
		=&~	O_p(\hat\sigma_1^2),\\
		b_{i,1}\bm{z}_{i,1}'\bm\epsilon_i
		=&~O_p(1+\hat\theta_1\bm{z}_{i,1}'\bm{z}_{i,1}).
		\end{split}
		\label{eq:dominates}
		\end{align}	
		\noindent By \eqref{eq:R*} and \eqref{eq:dominates}, we have 
		\begin{align}
		\begin{split}
		R^*(\hat\sigma_1^2,\hat{v}^2)
		=&~o_p(n^{-1}).
		\end{split}
		\label{eq:uniform bound R*}		
		\end{align}
		\noindent Similarly, by \eqref{partial:v}, we have
		\begin{align*}
		\frac{\partial}{\partial v^2}\{-2\log L(\theta_1,v^2)\}
		=&~	\frac{N}{v^2} - \frac{1}{v^4}
		\sum_{i=1}^m \bm{y}_i'\bigg(\bm{I}_n - \frac{\theta_1\bm{z}_{i,1}\bm{z}_{i,1}'}{1+\theta_1\bm{z}_{i,1}'\bm{z}_{i,1}}\bigg)\bm{y}_i
		\\
		=&~	\frac{N}{v^2}-\frac{1}{v^4}
		\sum_{i=1}^m
		(\bm{z}_{i,1}b_{i,1}+\bm\epsilon_i)'\bigg(\bm{I}_n - \frac{\theta_1\bm{z}_{i,1}\bm{z}_{i,1}'}{1+\theta_1\bm{z}_{i,1}'\bm{z}_{i,1}}\bigg)(\bm{z}_{i,1}b_{i,1}+\bm\epsilon_i)\\
		=&~	\frac{N}{v^2}-\frac{1}{v^4}
		\sum_{i=1}^m
		\bigg\{
		\frac{b_{i,1}^2\bm{z}_{i,1}'\bm{z}_{i,1}}{1+\theta_1\bm{z}_{i,1}'\bm{z}_{i,1}}
		+\frac{2b_{i,1}\bm{z}_{i,1}'\bm\epsilon_i}
		{1+\theta_1\bm{z}_{i,1}'\bm{z}_{i,1}}
		+\bm\epsilon_i'\bm\epsilon_i
		-\frac{\theta_1(\bm{z}_{i,1}'\bm\epsilon_i)^2}{1+\theta_1\bm{z}_{i,1}'\bm{z}_{i,1}}
		\bigg\}.
		\end{align*}
		\noindent The ML estimators $\hat\theta_1$ and $\hat{v}^2$ satisfy
		\begin{align*}
		0=&~ \frac{N}{\hat{v}^2} - \frac{N}{\hat{v}^4}
		\sum_{i=1}^m \bigg\{
		\frac{b_{i,1}^2\bm{z}_{i,1}'\bm{z}_{i,1}}{1+\hat\theta_1\bm{z}_{i,1}'\bm{z}_{i,1}}
		+\frac{2b_{i,1}\bm{z}_{i,1}'\bm\epsilon_i}
		{1+\hat\theta_1\bm{z}_{i,1}'\bm{z}_{i,1}}
		+\bm\epsilon_i'\bm\epsilon_i
		-\frac{\hat\theta_1(\bm{z}_{i,1}'\bm\epsilon_i)^2}{1+\hat\theta_1\bm{z}_{i,1}'\bm{z}_{i,1}}
		\bigg\},
		\end{align*} 
		\noindent which implies that 
		\begin{align}
		\hat{v}^2
		=&~	\frac{1}{N}\sum_{i=1}^m\bm\epsilon_i'\bm\epsilon_i
		+R^{\dag}(\hat\sigma_1^2,\hat{v}^2),
		\label{proof:exm:mle2}
		\end{align}
		\noindent with
		\begin{align*}
		R^{\dag}(\hat\sigma_1^2,\hat{v}^2)
		=&~	
		\frac{1}{N}\sum_{i=1}^m \bigg\{
		\frac{b_{i,1}^2\bm{z}_{i,1}'\bm{z}_{i,1}}{1+\hat\theta_1\bm{z}_{i,1}'\bm{z}_{i,1}}
		+\frac{2b_{i,1}\bm{z}_{i,1}'\bm\epsilon_i}
		{1+\hat\theta_1\bm{z}_{i,1}'\bm{z}_{i,1}}
		-\frac{\hat\theta_1(\bm{z}_{i,1}'\bm\epsilon_i)^2}{1+\hat\theta_1\bm{z}_{i,1}'\bm{z}_{i,1}}
		\bigg\}.
		\end{align*}
		\noindent This together with \eqref{eq:dominates} yields
		\begin{align}
		\begin{split}
		R^{\dag}(\hat\sigma_1^2,\hat{v}^2
		=&~ O_p(n^{-1}).
		\end{split}
		\label{eq:uniform R dag}
		\end{align}

		We are now ready to compare the asymptotic behaviors between the LS
		predictors and the empirical BLUPs. Note that for $i=1,\ldots,m$, we have
		\begin{align}
		\tilde{b}_{i,1}
		=&~	(\bm{z}_{i,1}'\bm{z}_{i,1})^{-1}\bm{z}_{i,1}'\bm{y}_i,
		\notag\\
		\hat{b}_{i,1}(\sigma_1^2,v^2)
		=&~	\sigma_1^2
		\bm{z}_{i,1}'(\sigma_1^2\bm{z}_{i,1}\bm{z}_{i,1}'+v^2\bm{I}_{n})^{-1}\bm{y}_i.
		\notag
		\end{align}
		\noindent Hence
		\begin{align}
		\bm{z}_{i,1}\big(\tilde{b}_{i,1}-b_{i,1}\big)
		=&~ \frac{\bm{z}_{i,1}\bm{z}_{i,1}'\bm\epsilon_i}{\bm{z}_{i,1}'\bm{z}_{i,1}},
		\label{proof:exm:eq01}
		\end{align}
		\noindent and 
		\begin{align}
		\begin{split}
		\hat{b}_{i,1}(\hat\sigma_1^2,\hat{v}^2) - b_{i,1}
		=&~	\hat{\sigma}_1^2
		\bm{z}_{i,1}'(\hat{\sigma}_1^2\bm{z}_{i,1}\bm{z}_{i,1}'+\hat{v}^2\bm{I}_{n})^{-1}(\bm{z}_{i,1}b_{i,1}+\bm\epsilon_i)-b_{i,1}\\
		=&~	\big\{\hat{\sigma}_1^2
		\bm{z}_{i,1}'(\hat{\sigma}_1^2\bm{z}_{i,1}\bm{z}_{i,1}'+\hat{v}^2\bm{I}_{n})^{-1}\bm{z}_{i,1}-1\big\}b_{i,1}+\hat{\sigma}_1^2
		\bm{z}_{i,1}'(\hat{\sigma}_1^2\bm{z}_{i,1}\bm{z}_{i,1}'+\hat{v}^2\bm{I}_{n})^{-1}\bm\epsilon_i\\
		=&~	\bigg\{(\hat{\sigma}_1^2/\hat{v}^2)
		\bm{z}_{i,1}'\bigg(
		\bm{I}_{n} 
		-\frac{(\hat{\sigma}_1^2/\hat{v}^2)\bm{z}_{i,1}\bm{z}_{i,1}'}
		{1+(\hat\sigma_1^2/\hat{v}^2)\bm{z}_{i,1}'\bm{z}_{i,1}}
		\bigg)\bm{z}_{i,1}-1\bigg\}b_i\\
		&~+(\hat{\sigma}_1^2/\hat{v}^2)
		\bm{z}_{i,1}'\bigg(
		\bm{I}_{n} 
		-\frac{(\hat{\sigma}_1^2/\hat{v}^2)\bm{z}_{i,1}\bm{z}_{i,1}'}
		{1+(\hat\sigma_1^2/\hat{v}^2)\bm{z}_{i,1}'\bm{z}_{i,1}}
		\bigg)\bm\epsilon_i\\
		=&~	\frac{(\hat{\sigma}_1^2/\hat{v}^2)\bm{z}_{i,1}'\bm\epsilon_i-b_i}{1+(\hat{\sigma}_1^2/\hat{v}^2)\bm{z}_{i,1}'\bm{z}_{i,1}},
		\end{split}
		\notag
		\end{align}
		\noindent which implies that
		\begin{align}
		\bm{z}_{i,1}\big(
		\hat{b}_{i,1}(\hat\sigma_1^2,\hat{v}^2) - b_{i,1}\big)
		=&~	\frac{\bm{z}_{i,1}\{(\hat{\sigma}_1^2/\hat{v}^2)\bm{z}_{i,1}'\bm\epsilon_i-b_i\}}{1+(\hat{\sigma}_1^2/\hat{v}^2)\bm{z}_{i,1}'\bm{z}_{i,1}}.
		\label{proof:exm:eq02}
		\end{align}
		\noindent Note that by \eqref{proof:exm:eq01},
	\begin{align*}
	\sum_{i=1}^m\big\|\bm{z}_{i,1}
	\big(\tilde{b}_{i,1}-b_{i,1}\big)\big\|^2
	=&~	\sum_{i=1}^m (\tilde{b}_{i,1}-b_{i,1})^2\bm{z}_{i,1}'\bm{z}_{i,1}
	=	\sum_{i=1}^m\frac{(\bm{z}_{i,1}'\bm\epsilon_i)^2}{\bm{z}_{i,1}'\bm{z}_{i,1}},
	\end{align*}
	\noindent and by \eqref{proof:exm:eq02},
	\begin{align*}
	\sum_{i=1}^m\big\|\bm{z}_{i,1}\big(
	\hat{b}_{i,1}(\hat\sigma_1^2,\hat{v}^2) - b_{i,1}\big)\big\|^2
	=&~ \sum_{i=1}^m\frac{\{(\hat{\sigma}_1^2/\hat{v}^2)\bm{z}_{i,1}'\bm\epsilon_i-b_{i,1}\}^2\bm{z}_{i,1}'\bm{z}_{i,1}}{\{1+(\hat{\sigma}_1^2/\hat{v}^2)\bm{z}_{i,1}'\bm{z}_{i,1}\}^2},
	\end{align*}
	\noindent which implies that 
	\begin{align}
	\begin{split}
	D(\hat\sigma^2,\hat{v}^2)
	=&~	\sum_{i=1}^m
	\bigg(
	\frac{(\bm{z}_{i,1}'\bm\epsilon_i)^2}{\bm{z}_{i,1}'\bm{z}_{i,1}}
	-\frac{\{(\hat{\sigma}_1^2/\hat{v}^2)\bm{z}_{i,1}'\bm\epsilon_i-b_{i,1}\}^2\bm{z}_{i,1}'\bm{z}_{i,1}}{\{1+(\hat{\sigma}_1^2/\hat{v}^2)\bm{z}_{i,1}'\bm{z}_{i,1}\}^2}
	\bigg)\\
	=&~	
	\sum_{i=1}^m
	\frac{(\bm{z}_{i,1}'\bm\epsilon_i)^2\{1+(\hat{\sigma}_1^2/\hat{v}^2)\bm{z}_{i,1}'\bm{z}_{i,1}\}^2
	-\{(\hat{\sigma}_1^2/\hat{v}^2)\bm{z}_{i,1}'\bm\epsilon_i-b_{i,1}\}^2(\bm{z}_{i,1}'\bm{z}_{i,1})^2}{\bm{z}_{i,1}'\bm{z}_{i,1}\{1+(\hat{\sigma}_1^2/\hat{v}^2)\bm{z}_{i,1}'\bm{z}_{i,1}\}^2}\\
	=&~\sum_{i=1}^m \frac{(\bm{z}_{i,1}'\bm\epsilon_i)^2+2(\hat{\sigma}_1^2/\hat{v}^2)(\bm{z}_{i,1}'\bm\epsilon_i)^2\bm{z}_{i,1}'\bm{z}_{i,1}+2(\hat{\sigma}_1^2/\hat{v}^2)b_{i,1}\bm{z}_{i,1}'\bm\epsilon_i(\bm{z}_{i,1}'\bm{z}_{i,1})^2-b_{i,1}^2(\bm{z}_{i,1}'\bm{z}_{i,1})^2}{\bm{z}_{i,1}'\bm{z}_{i,1}\{1+(\hat{\sigma}_1^2/\hat{v}^2)\bm{z}_{i,1}'\bm{z}_{i,1}\}^2}.
	\end{split}
	\label{proof:exm:D function}
	\end{align}
	\noindent Note that by \eqref{proof:exm:D function} and
	\begin{align*}
	\frac{2(\hat{\sigma}_1^2/\hat{v}^2)(\bm{z}_{i,1}'\bm\epsilon_i)^2}{\{1+(\hat{\sigma}_1^2/\hat{v}^2)\bm{z}_{i,1}'\bm{z}_{i,1}\}^2}
	=&~	\frac{2(\bm{z}_{i,1}'\bm\epsilon_i)^2}{(\bm{z}_{i,1}'\bm{z}_{i,1})^2(\hat\sigma_1^2/\hat{v}^2)}
	-\frac{2(\bm{z}_{i,1}'\bm\epsilon_i)^2+4(\hat\sigma_1^2/\hat{v}^2)(\bm{z}_{i,1}'\bm{z}_{i,1})(\bm{z}_{i,1}'\bm\epsilon_i)^2}{\{1+(\hat\sigma_1^2/\hat{v}^2)\bm{z}_{i,1}'\bm{z}_{i,1}\}^2
	(\hat\sigma_1^2/\hat{v}^2)(\bm{z}_{i,1}'\bm{z}_{i,1})^2},\\
	\frac{2(\hat{\sigma}_1^2/\hat{v}^2)b_{i,1}\bm{z}_{i,1}'\bm\epsilon_i(\bm{z}_{i,1}'\bm{z}_{i,1})}{\{1+(\hat{\sigma}_1^2/\hat{v}^2)\bm{z}_{i,1}'\bm{z}_{i,1}\}^2}
	=&~	\frac{2b_{i,1}\bm{z}_{i,1}'\bm\epsilon_i}{(\hat\sigma_1^2/\hat{v}^2)\bm{z}_{i,1}'\bm{z}_{i,1}}
	-\frac{2b_{i,1}\bm{z}_{i,1}'\bm\epsilon_i+4b_{i,1}\bm{z}_{i,1}'\bm\epsilon_i(\hat\sigma_1^2/\hat{v}^2)(\bm{z}_{i,1}'\bm{z}_{i,1})}{\{1+(\hat{\sigma}_1^2/\hat{v}^2)\bm{z}_{i,1}'\bm{z}_{i,1}\}^2(\hat\sigma_1^2/\hat{v}^2)(\bm{z}_{i,1}'\bm{z}_{i,1})},\\
	\frac{b_{i,1}^2(\bm{z}_{i,1}'\bm{z}_{i,1})}{\{1+(\hat{\sigma}_1^2/\hat{v}^2)\bm{z}_{i,1}'\bm{z}_{i,1}\}^2}
	=&~ \frac{b_{i,1}^2}{(\hat\sigma_1^2/\hat{v}^2)^2(\bm{z}_{i,1}'\bm{z}_{i,1})}
	-\frac{b_{i,1}^2+2(\hat\sigma_1^2/\hat{v}^2)\bm{z}_{i,1}'\bm{z}_{i,1}}{\{1+(\hat{\sigma}_1^2/\hat{v}^2)\bm{z}_{i,1}'\bm{z}_{i,1}\}^2(\hat\sigma_1^2/\hat{v}^2)^2(\bm{z}_{i,1}'\bm{z}_{i,1})}, 
	\end{align*}	
	\noindent we have 
	\begin{align}
	D(\hat\sigma^2,\hat{v}^2)
	=&~	
	\sum_{i=1}^m
	\bigg(
	\frac{2(\bm{z}_{i,1}'\bm\epsilon_i)^2}{(\bm{z}_{i,1}'\bm{z}_{i,1})^2(\hat\sigma_1^2/\hat{v}^2)}
	+\frac{2b_{i,1}\bm{z}_{i,1}'\bm\epsilon_i}{(\hat\sigma_1^2/\hat{v}^2)\bm{z}_{i,1}'\bm{z}_{i,1}}
	-\frac{b_{i,1}^2}{(\hat\sigma_1^2/\hat{v}^2)^2(\bm{z}_{i,1}'\bm{z}_{i,1})}
	\bigg)+R^{\ddag}(\hat\sigma_1^2,\hat{v}^2)
	\label{proof:exm:D function 2}
	\end{align}
	\noindent with
	\begin{align*}
	R^{\ddag}(\hat\sigma_1^2,\hat{v}^2)
	=&~	\sum_{i=1}^m
	\bigg(
	\frac{(\bm{z}_{i,1}'\bm\epsilon_i)^2}{\bm{z}_{i,1}'\bm{z}_{i,1}\{1+(\hat{\sigma}_1^2/\hat{v}^2)\bm{z}_{i,1}'\bm{z}_{i,1}\}^2}
	-\frac{2(\bm{z}_{i,1}'\bm\epsilon_i)^2+4(\hat\sigma_1^2/\hat{v}^2)(\bm{z}_{i,1}'\bm{z}_{i,1})(\bm{z}_{i,1}'\bm\epsilon_i)^2}{\{1+(\hat\sigma_1^2/\hat{v}^2)\bm{z}_{i,1}'\bm{z}_{i,1}\}^2
		(\hat\sigma_1^2/\hat{v}^2)(\bm{z}_{i,1}'\bm{z}_{i,1})^2}\\
	&~	-\frac{2b_{i,1}\bm{z}_{i,1}'\bm\epsilon_i+4b_{i,1}\bm{z}_{i,1}'\bm\epsilon_i(\hat\sigma_1^2/\hat{v}^2)(\bm{z}_{i,1}'\bm{z}_{i,1})}{\{1+(\hat{\sigma}_1^2/\hat{v}^2)\bm{z}_{i,1}'\bm{z}_{i,1}\}^2(\hat\sigma_1^2/\hat{v}^2)(\bm{z}_{i,1}'\bm{z}_{i,1})}
	+\frac{b_{i,1}^2+2(\hat\sigma_1^2/\hat{v}^2)\bm{z}_{i,1}'\bm{z}_{i,1}}{\{1+(\hat{\sigma}_1^2/\hat{v}^2)\bm{z}_{i,1}'\bm{z}_{i,1}\}^2(\hat\sigma_1^2/\hat{v}^2)^2(\bm{z}_{i,1}'\bm{z}_{i,1})}
	\bigg).
	\end{align*}
	\noindent Note that by \eqref{eq:dominates},
	\begin{align}
	\begin{split}
	R^{\ddag}(\hat\sigma_1^2,\hat{v}^2)
	=&~O_p(n^{-3/2}).
	\end{split}
	\label{eq:uniform bound R dag}
	\end{align}
	\noindent Further, by \eqref{proof:exm:mle1} and \eqref{proof:exm:mle2}, we
	have
	\begin{align}
	\begin{split}
		\frac{2(\bm{z}_{i,1}'\bm\epsilon_i)^2}{(\bm{z}_{i,1}'\bm{z}_{i,1})^2(\hat\sigma_1^2/\hat{v}^2)}
	=&~
	\frac{2(\bm{z}_{i,1}'\bm\epsilon_i)^2(\sum_{k=1}^m\bm\epsilon_k'\bm\epsilon_k/N)}{(\bm{z}_{i,1}'\bm{z}_{i,1})^2(\sum_{k=1}^mb_{k,1}^2/m)}+ \frac{2(\bm{z}_{i,1}'\bm\epsilon_i)^2
	}{(\bm{z}_{i,1}'\bm{z}_{i,1})^2\hat\sigma_1^2(\sum_{k=1}^mb_{k,1}^2/m)}\\
&~	\times \bigg\{R^{\dag}(\hat\sigma_1^2,\hat{v}^2)\frac{\sum_{k=1}^mb_{k,1}^2}{m}-\bigg(\frac{\sum_{k=1}^m\bm\epsilon_k'\bm\epsilon_k}{N}\bigg)\bigg(\sum_{i=1}^m\frac{2b_{i,1}\bm{z}_{i,1}'\bm\epsilon_i}{m\bm{z}_{i,1}'\bm{z}_{i,1}}\bigg)-R^{*}(\hat\sigma_1^2,\hat{v}^2)\frac{\sum_{k=1}^m\bm\epsilon_k'\bm\epsilon_k}{N}\bigg\}\\
	\equiv &~
	\frac{2(\bm{z}_{i,1}'\bm\epsilon_i)^2\sum_{k=1}^m\bm\epsilon_k'\bm\epsilon_k}{n(\bm{z}_{i,1}'\bm{z}_{i,1})^2\sum_{k=1}^mb_{k,1}^2}+ R_{i,1}(\hat\sigma_1^2,\hat{v}^2),
	\end{split}
	\label{proof:exm:D function 21}
	\end{align}
	\noindent with 
	\begin{align*}
	R_{i,1}(\hat\sigma_1^2,\hat{v}^2)
	=&~	\frac{2(\bm{z}_{i,1}'\bm\epsilon_i)^2
	}{(\bm{z}_{i,1}'\bm{z}_{i,1})^2\hat\sigma_1^2(\sum_{k=1}^mb_{k,1}^2/m)}
\bigg\{R^{\dag}(\hat\sigma_1^2,\hat{v}^2)\frac{\sum_{k=1}^mb_{k,1}^2}{m}\\
&~-\bigg(\frac{\sum_{k=1}^m\bm\epsilon_k'\bm\epsilon_k}{N}\bigg)\bigg(\sum_{i=1}^m\frac{2b_{i,1}\bm{z}_{i,1}'\bm\epsilon_i}{m\bm{z}_{i,1}'\bm{z}_{i,1}}\bigg)-R^{*}(\hat\sigma_1^2,\hat{v}^2)\frac{\sum_{k=1}^m\bm\epsilon_k'\bm\epsilon_k}{N}\bigg\}.
	\end{align*}
	\noindent 
	Similarly,
	\begin{align}
	\frac{2b_{i,1}\bm{z}_{i,1}'\bm\epsilon_i}{(\hat\sigma_1^2/\hat{v}^2)\bm{z}_{i,1}'\bm{z}_{i,1}}
	=&~
	\frac{2b_{i,1}\bm{z}_{i,1}'\bm\epsilon_i\sum_{k=1}^m\bm\epsilon_k'\bm\epsilon_k}
	{n\{\sum_{k=1}^mb_{k,1}^2+2\sum_{k=1}^mb_{k,1}\bm{z}_{k,1}'\bm\epsilon_k/(\bm{z}_{k,1}'\bm{z}_{k,1})\}(\bm{z}_{i,1}'\bm{z}_{i,1})}+R_{i,2}(\hat\sigma_1^2,\hat{v}^2),
	\label{proof:exm:D function 22}\\
	\frac{b_{i,1}^2}{(\hat\sigma_1^2/\hat{v}^2)^2(\bm{z}_{i,1}'\bm{z}_{i,1})}
	=&~	\frac{b_{i,1}^2(\sum_{k=1}^m\bm\epsilon_k'\bm\epsilon_k)^2}{n^2(\sum_{k=1}^mb_{k,1}^2)^2\bm{z}_{i,1}'\bm{z}_{i,1}}+R_{i,3}(\hat\sigma_1^2,\hat{v}^2)
	\label{proof:exm:D function 23}
	\end{align} 
	\noindent with
	\begin{align*}
	R_{i,2}(\hat\sigma_1^2,\hat{v}^2)
	=&~\frac{2b_{i,1}\bm{z}_{i,1}'\bm\epsilon_i}
	{\hat\sigma_1^2\bm{z}_{i,1}'\bm{z}_{i,1}
	\{\sum_{k=1}^mb_{k,1}^2/m+2\sum_{k=1}^mb_{k,1}\bm{z}_{k,1}'\bm\epsilon_k/(m\bm{z}_{k,1}'\bm{z}_{k,1})\}}\\
&~\times \bigg\{
R^{\dag}(\hat\sigma_1^2,\hat{v}^2)
\bigg(\frac{\sum_{k=1}^mb_{k,1}^2}{m}
+\sum_{k=1}^m\frac{2b_{k,1}\bm{z}_{k,1}'\bm\epsilon_k}{m\bm{z}_{k,1}'\bm{z}_{k,1}}\bigg)
-R^{*}(\hat\sigma_1^2,\hat{v}^2)\sum_{k=1}^m\frac{\bm\epsilon_k'\bm\epsilon_k}{N}
\bigg\},\\
	R_{i,3}(\hat\sigma_1^2,\hat{v}^2)
	=&~\frac{b_{i,1}^2}{\hat\sigma_1^4(\bm{z}_{i,1}'\bm{z}_{i,1})\{\sum_{k=1}^mb_{k,1}^2/m\}^2}
	\bigg(\hat{v}^2\sum_{k=1}^m\frac{b_{k,1}^2}{m}
	+\hat\sigma_1^2\sum_{k=1}^m\frac{\bm\epsilon_k'\bm\epsilon_k}{N}\bigg)\\
	&~	\times 
	\bigg\{R^{\dag}(\hat\sigma_1^2,\hat{v}^2)\frac{\sum_{k=1}^mb_{k,1}^2}{m}-\bigg(\frac{\sum_{k=1}^m\bm\epsilon_k'\bm\epsilon_k}{N}\bigg)\bigg(\sum_{i=1}^m\frac{2b_{i,1}\bm{z}_{i,1}'\bm\epsilon_i}{m\bm{z}_{i,1}'\bm{z}_{i,1}}\bigg)-R^{*}(\hat\sigma_1^2,\hat{v}^2)\frac{\sum_{k=1}^m\bm\epsilon_k'\bm\epsilon_k}{N}\bigg\}.
	\end{align*}
	\noindent Hence by \eqref{eq:dominates}, \eqref{eq:uniform bound R*}, and
	\eqref{eq:uniform R dag}, we have 
	\begin{align}
	\begin{split}
	R_{i,1}(\hat\sigma_1^2,\hat{v}^2)
	=&~O_p(n^{-3/2}),\quad i=1,2,3.
	\end{split}
	\label{eq:uniform bound R123}
	\end{align}
	\noindent Furthermore, we have
	\begin{align}
	\begin{split}
		\frac{1}{n}&\frac{2b_{i,1}\bm{z}_{i,1}'\bm\epsilon_i\sum_{k=1}^m\bm\epsilon_k'\bm\epsilon_k}
	{\{\sum_{k=1}^mb_{k,1}^2+2\sum_{k=1}^mb_{k,1}\bm{z}_{k,1}'\bm\epsilon_k/(m\bm{z}_{k,1}'\bm{z}_{k,1})\}(\bm{z}_{i,1}'\bm{z}_{i,1})}\\
	=&~ \frac{1}{n(\bm{z}_{i,1}'\bm{z}_{i,1})}\bigg\{
	\frac{2b_{i,1}\bm{z}_{i,1}'\bm\epsilon_i\sum_{k=1}^m\bm\epsilon_k'\bm\epsilon_k}
	{\sum_{k=1}^mb_{k,1}^2}-\frac{4b_{i,1}\bm{z}_{i,1}'\bm\epsilon_i\{\sum_{k=1}^m\bm\epsilon_k'\bm\epsilon_k\}\{\sum_{k=1}^mb_{k,1}\bm{z}_{k,1}'\bm\epsilon_k/(m\bm{z}_{k,1}'\bm{z}_{k,1})\}}{\{\sum_{k=1}^mb_{k,1}^2\}^2}\\
	&~ +\frac{8b_{i,1}\bm{z}_{i,1}'\bm\epsilon_i\sum_{k=1}^m\bm\epsilon_k'\bm\epsilon_k\{\sum_{k=1}^mb_{k,1}\bm{z}_{k,1}'\bm\epsilon_k/(m\bm{z}_{k,1}'\bm{z}_{k,1})\}^2}{\{\sum_{k=1}^mb_{k,1}^2+2\sum_{k=1}^mb_{k,1}\bm{z}_{k,1}'\bm\epsilon_k/(m\bm{z}_{k,1}'\bm{z}_{k,1})\}\{\sum_{k=1}^mb_{k,1}^2\}^2}\bigg\}\\
	\equiv &~ \bigg\{
	\frac{2b_{i,1}\bm{z}_{i,1}'\bm\epsilon_i\sum_{k=1}^m\bm\epsilon_k'\bm\epsilon_k}
	{n(\bm{z}_{i,1}'\bm{z}_{i,1})\sum_{k=1}^mb_{k,1}^2}-\frac{4b_{i,1}\bm{z}_{i,1}'\bm\epsilon_i\{\sum_{k=1}^m\bm\epsilon_k'\bm\epsilon_k\}\{\sum_{k=1}^mb_{k,1}\bm{z}_{k,1}'\bm\epsilon_k/(m\bm{z}_{k,1}'\bm{z}_{k,1})\}}{n(\bm{z}_{i,1}'\bm{z}_{i,1})\{\sum_{k=1}^mb_{k,1}^2\}^2}\bigg\}\\
	&~+R_{i,4},
	\end{split}
	\label{proof:exm:D function 24}
	\end{align}
	\noindent with 
	\begin{align}
	R_{i,4}
	=&~	\frac{8b_{i,1}\bm{z}_{i,1}'\bm\epsilon_i\sum_{k=1}^m\bm\epsilon_k'\bm\epsilon_k\{\sum_{k=1}^mb_{k,1}\bm{z}_{k,1}'\bm\epsilon_k/(m\bm{z}_{k,1}'\bm{z}_{k,1})\}^2}{n(\bm{z}_{i,1}'\bm{z}_{i,1})\{\sum_{k=1}^mb_{k,1}^2+2\sum_{k=1}^mb_{k,1}\bm{z}_{k,1}'\bm\epsilon_k/(m\bm{z}_{k,1}'\bm{z}_{k,1})\}\{\sum_{k=1}^mb_{k,1}^2\}^2}.
	\notag
	\end{align}
	\noindent Note that
	\begin{align}
	R_{i,4} = O_p(n^{-3/2}).
	\label{eq:uniform bound R4}
	\end{align}
	\noindent By \eqref{proof:exm:D function 2}, \eqref{proof:exm:D function
	21}, \eqref{proof:exm:D function 22}, \eqref{proof:exm:D function 23}, and
	\eqref{proof:exm:D function 24}, we have
	\begin{align}
	\begin{split}
	nD(\hat\sigma_1^2,\hat{v}^2)
	=&~	A_{n,m}+nR^{\ddag}(\hat\sigma_1^2,\hat{v}^2)+n\sum_{i=1}^m
	\bigg\{	R_{i,1}(\hat\sigma_1^2,\hat{v}^2)
	+	R_{i,2}(\hat\sigma_1^2,\hat{v}^2)	-R_{i,3}(\hat\sigma_1^2,\hat{v}^2)
	+	R_{i,4}
	\bigg\}\\
	\equiv&~A_{n,m}+O_p(n^{-1/2})
	\end{split}
	\notag
	\end{align}
	\noindent with 
	\begin{align*}
	A_{n,m}=&~\sum_{i=1}^m
	\bigg\{
	\frac{2(\bm{z}_{i,1}'\bm\epsilon_i)^2\sum_{k=1}^m\bm\epsilon_k'\bm\epsilon_k}{(\bm{z}_{i,1}'\bm{z}_{i,1})^2\sum_{k=1}^mb_{k,1}^2} - \frac{b_{i,1}^2(\sum_{k=1}^m\bm\epsilon_k'\bm\epsilon_k)^2}{n(\sum_{k=1}^mb_{k,1}^2)^2\bm{z}_{i,1}'\bm{z}_{i,1}}
	\frac{2b_{i,1}\bm{z}_{i,1}'\bm\epsilon_i\sum_{k=1}^m\bm\epsilon_k'\bm\epsilon_k}
	{(\bm{z}_{i,1}'\bm{z}_{i,1})\sum_{k=1}^mb_{k,1}^2}\\
	&~-\frac{4b_{i,1}\bm{z}_{i,1}'\bm\epsilon_i\{\sum_{k=1}^m\bm\epsilon_k'\bm\epsilon_k\}\{\sum_{k=1}^mb_{k,1}\bm{z}_{k,1}'\bm\epsilon_k/(\bm{z}_{k,1}'\bm{z}_{k,1})\}}{(\bm{z}_{i,1}'\bm{z}_{i,1})\{\sum_{k=1}^mb_{k,1}^2\}^2}
	\bigg\},
	\end{align*}
	\noindent where the last equality follows from 
	 \eqref{eq:uniform bound R dag}, \eqref{eq:uniform bound R123}, and
	 \eqref{eq:uniform bound R4}. Note that
	 $\big(\sum_{k=1}^mb_{i,k}^2/\sigma_{1,0}^2\big)^{-1}$ follows the
	 inverse-chi-squared distribution with $m$ degrees of freedom. We have   
	\begin{align}
	\begin{split}
	\mathrm{E}\bigg(\frac{1}{\sum_{i=1}^m b_{i,1}^2}\bigg)
	=&~	\frac{1}{(m-2)\sigma_{1,0}^2}, \quad \mbox{provided }m>2,\\
	\mathrm{E}
	\bigg(\frac{b_{i,1}^2}{\{\sum_{k=1}^mb_{k,1}^2\}^2}\bigg)
	=&~	\frac{1}{m(m-2)\sigma_{1,0}^2}\quad \mbox{provided }m>4.
	\end{split}
	\label{proof:exm:inverse chisquare}
	\end{align}
	\noindent By \eqref{proof:exm:inverse chisquare} and 
	\begin{align}
	\begin{split}
	\mathrm{E}
	\bigg(\bigg\{\sum_{i=1}^m \bm\epsilon_i'\bm\epsilon_i\bigg\}^2\bigg)
	=&~	(2mn+m^2n^2)v_0^4,\\
	\mathrm{E}\big(\bm\epsilon_i'\bm\epsilon_i(\bm{z}_{i,1}\bm\epsilon_i)\big)
	=&~ 0,\\
	\mathrm{E}\big(\bm\epsilon_i'\bm\epsilon_i(\bm{z}_{i,1}\bm\epsilon_i)^2\big)
	=&~ n^2v_0^4 +o(n^2),
	\end{split}
	\notag
	\end{align}
	\noindent we have, for $m>4$, 
	\begin{align}
	\begin{split}
	\mathrm{E}(A_{n,m})
	=&~\mathrm{E}\sum_{i=1}^m
	\bigg\{
	\frac{2(\bm{z}_{i,1}'\bm\epsilon_i)^2\sum_{k=1}^m\bm\epsilon_k'\bm\epsilon_k}{(\bm{z}_{i,1}'\bm{z}_{i,1})^2\sum_{k=1}^mb_{k,1}^2} - \frac{b_{i,1}^2(\sum_{k=1}^m\bm\epsilon_k'\bm\epsilon_k)^2}{n(\sum_{k=1}^mb_{k,1}^2)^2\bm{z}_{i,1}'\bm{z}_{i,1}}+
	\frac{2b_{i,1}\bm{z}_{i,1}'\bm\epsilon_i\sum_{k=1}^m\bm\epsilon_k'\bm\epsilon_k}
	{(\bm{z}_{i,1}'\bm{z}_{i,1})\sum_{k=1}^mb_{k,1}^2}\\
	&~-\frac{4b_{i,1}\bm{z}_{i,1}'\bm\epsilon_i\{\sum_{k=1}^m\bm\epsilon_k'\bm\epsilon_k\}\{\sum_{k=1}^mb_{k,1}\bm{z}_{k,1}'\bm\epsilon_k/(\bm{z}_{k,1}'\bm{z}_{k,1})\}}{(\bm{z}_{i,1}'\bm{z}_{i,1})\{\sum_{k=1}^mb_{k,1}^2\}^2}
	\bigg\}\\
	=&~ \frac{2m^2v_0^4}{(m-2)\sigma_{1,0}^2}
	-\frac{m^2v_0^4}{(m-2)\sigma_{1,0}^2}+o(1)\\
	&~-\mathrm{E}
	\bigg(\mathrm{E}\bigg(\sum_{i=1}^m\frac{4b_{i,1}\bm{z}_{i,1}'\bm\epsilon_i\{\sum_{k=1}^m\bm\epsilon_k'\bm\epsilon_k\}\{\sum_{k=1}^mb_{k,1}\bm{z}_{k,1}'\bm\epsilon_k/(\bm{z}_{k,1}'\bm{z}_{k,1})\}}{(\bm{z}_{i,1}'\bm{z}_{i,1})\{\sum_{k=1}^mb_{k,1}^2\}^2}
	\bigg|b_{1,1},\ldots,b_{m,1}\bigg)\bigg)\\
	=&~\frac{2m^2v_0^4}{(m-2)\sigma_{1,0}^2}
	-\frac{m^2v_0^4}{(m-2)\sigma_{1,0}^2}
	-\mathrm{E}
	\bigg(\sum_{i=1}^m 
	\frac{4mv_0^4b_{i,1}^2}{\{\sum_{k=1}^mb_{k,1}^2\}^2}
	\bigg)+o(1)\\
	=&~	\frac{2m^2v_0^4}{(m-2)\sigma_{1,0}^2}
	-\frac{m^2v_0^4}{(m-2)\sigma_{1,0}^2}
	-\frac{4mv_0^4}{(m-2)\sigma_{1,0}^2}+o(1)\\
	=&~ 
	\frac{m(m-4)v_0^4}{(m-2)\sigma_{1,0}^2}+o(1).
	\end{split}
		\notag
	\end{align}
	\noindent This completes the proofs.

	\subsection{Proof of Theorem \ref{appendix:theorem:MLE 4}}
	
	In this section, we first prove Theorem \ref{appendix:theorem:MLE 4} to
	simplify the proofs of Theorems \ref{appendix:theorem:MLE 2} and
	\ref{appendix:theorem:MLE 3}. 
	As with the proof of Theorem \ref{appendix:theorem:MLE},
	we shall focus on the asymptotic properties of $\hat{v}^2(\alpha,\gamma)$
	and $\{\hat{\theta}_k(\alpha,\gamma):k\in\gamma\}$,
	and derive them by solving the likelihood equations directly.
	
	We first prove \eqref{appendix:thm:mle4:correct:eq1} using \eqref{partial:v}.
	For $(\alpha,\gamma)\in(\mathcal{A}\setminus\mathcal{A}_0)\times\mathcal{G}$, we have
	\begin{align}
	(\bm{I}_{N}-\bm{M}(\alpha,\gamma;\bm\theta))\bm\mu_{0}
	=&~ (\bm{I}_{N}-\bm{M}(\alpha,\gamma;\bm\theta))
	\bm{X}(\alpha_0\setminus\alpha)\bm{\beta}_0(\alpha_0\setminus\alpha),
	\label{proof:thm:mle part3:eq0}
	\end{align}
	
	\noindent where $\bm{\beta}_0(\alpha_0\setminus\alpha)$ denotes the
	sub-vector of $\bm{\beta}_0$ corresponding to $\alpha_0\setminus\alpha$.
	Note that by the Cauchy--Schwarz inequality, we have
	\begin{align}
	\bigg(\sum_{i=1}^m n_i^{(\xi+\ell)/2}\bigg)^2=O\bigg(\sum_{i=1}^m n_i^{\xi}\sum_{i^*=1}^m n_i^{\ell}\bigg).
	\label{cauchy2}
	\end{align}
	\noindent Hence
	by \eqref{cauchy2} and Lemma \ref{appendix:lemma:xze}, we have
	\begin{align*}
	\begin{split}
	\big(\bm{X}&(\alpha_0\setminus\alpha)\bm{\beta}_0(\alpha_0\setminus\alpha)
	+\bm{Z}(\gamma_0)\bm{b}(\gamma_0)+\bm\epsilon\big)'
	\bm{H}^{-1}(\gamma,\bm\theta)\bm{M}(\alpha,\gamma;\bm\theta)\\
	&~  \times\big(\bm{X}(\alpha_0\setminus\alpha)\bm{\beta}_0(\alpha_0\setminus\alpha)
	+\bm{Z}(\gamma_0)\bm{b}(\gamma_0)+\bm\epsilon\big)\\
	=&~ \bm{\beta}_0(\alpha_0\setminus\alpha)'\big(\bm{X}(\alpha_0\setminus\alpha)'
	\bm{H}^{-1}(\gamma,\bm\theta)\bm{M}(\alpha,\gamma;\bm\theta)
	\bm{X}(\alpha_0\setminus\alpha)\bm{\beta}_0(\alpha_0\setminus\alpha)\\
	&~  +  \bm{b}(\gamma_0)'\bm{Z}(\gamma_0)'
	\bm{H}^{-1}(\gamma,\bm\theta)\bm{M}(\alpha,\gamma;\bm\theta)\bm{Z}(\gamma_0)\bm{b}(\gamma_0)\\
	&~  +\bm{\epsilon}'\bm{H}^{-1}(\gamma,\bm\theta)\bm{M}(\alpha,\gamma;\bm\theta)\bm{\epsilon}\\
	&~  +2\bm{b}(\gamma_0)'\bm{Z}(\gamma_0)'
	\bm{H}^{-1}(\gamma,\bm\theta)\bm{M}(\alpha,\gamma;\bm\theta)
	\bm{X}(\alpha_0\setminus\alpha)\bm{\beta}_0(\alpha_0\setminus\alpha)\\
	&~  +2\bm\epsilon'\bm{H}^{-1}(\gamma,\bm\theta)\bm{M}(\alpha,\gamma;\bm\theta)
	\bm{X}(\alpha_0\setminus\alpha)\bm{\beta}_0(\alpha_0\setminus\alpha)\\
	&~  +2\bm{b}(\gamma_0)'\bm{Z}(\gamma_0)'
	\bm{H}^{-1}(\gamma,\bm\theta)\bm{M}(\alpha,\gamma;\bm\theta)\bm{\epsilon}\\
	=&~ \bigg(\sum_{i=1}^m\sum_{k\in\gamma_0}b_{i,k}\bm{h}_{i,k}\bigg)'
	\bm{H}^{-1}(\gamma,\bm\theta)\bm{M}(\alpha,\gamma;\bm\theta)\bigg(\sum_{i=1}^m\sum_{k\in\gamma_0}b_{i,k}\bm{h}_{i,k}\bigg)\\
	&~  +2\bigg(\sum_{i=1}^m\sum_{k\in\gamma_0}b_{i,k}\bm{h}_{i,k}\bigg)'
	\bm{H}^{-1}(\gamma,\bm\theta)\bm{M}(\alpha,\gamma;\bm\theta)
	\bm{X}(\alpha_0\setminus\alpha)\bm{\beta}_0(\alpha_0\setminus\alpha)\\
	&~  +2\bigg(\sum_{i=1}^m\sum_{k\in\gamma_0}b_{i,k}\bm{h}_{i,k}\bigg)'
	\bm{H}^{-1}(\gamma,\bm\theta)\bm{M}(\alpha,\gamma;\bm\theta)\bm{\epsilon}
	+o_p\bigg(\sum_{i=1}^m n_i^{\xi-\tau}\bigg)\\
	&~  +O_p(p)\\
	\end{split}
	\end{align*}
	\begin{align*}
	\begin{split}
	=&~ o_p\bigg(\sum_{i=1}^mn_i^{\ell-\tau}\bigg)
	+o_p\bigg(\sum_{i=1}^mn_i^{(\xi+\ell)/2-\tau}\bigg)
	+o_p\bigg(\sum_{i=1}^m n_i^{\xi-\tau}\bigg)+O_p(p)\\
	=&~ o_p\bigg(\sum_{i=1}^mn_i^\xi\bigg)+
	o_p\bigg(\sum_{i=1}^mn_i^\ell\bigg)
	+O_p(p)
	\end{split}
	%\label{proof:thm:mle part2:eq16}
	\end{align*}
	
	\noindent uniformly over $\bm\theta\in\Theta_\gamma$.
	This and (\ref{proof:thm:mle part3:eq0}) imply
	\begin{align*}
	\begin{split}
	\bm{y}'\bm{H}^{-1}&(\gamma,\bm\theta)(\bm{I}_{N}-\bm{M}(\alpha,\gamma;\bm\theta))\bm{y}\\
	=&~ \big(\bm{X}(\alpha_0\setminus\alpha)\bm{\beta}_0(\alpha_0\setminus\alpha)
	+\bm{Z}(\gamma_0)\bm{b}(\gamma_0)+\bm\epsilon\big)'\\
	&~    \times\bm{H}^{-1}(\gamma,\bm\theta)
	(\bm{I}_{N}-\bm{M}(\alpha,\gamma;\bm\theta))\\
	&~  \times\big(\bm{X}(\alpha_0\setminus\alpha)\bm{\beta}_0(\alpha_0\setminus\alpha)
	+\bm{Z}(\gamma_0)\bm{b}(\gamma_0)+\bm\epsilon\big)\\
	=&~  \big(\bm{X}(\alpha_0\setminus\alpha)\bm{\beta}_0(\alpha_0\setminus\alpha)
	+\bm{Z}(\gamma_0)\bm{b}(\gamma_0)+\bm\epsilon\big)'
	\bm{H}^{-1}(\gamma,\bm\theta)\\
	&~  \times \big(\bm{X}(\alpha_0\setminus\alpha)\bm{\beta}_0(\alpha_0\setminus\alpha)
	+\bm{Z}(\gamma_0)\bm{b}(\gamma_0)+\bm\epsilon\big)\\
	&~  +o_p\bigg(\sum_{i=1}^mn_i^\xi\bigg)
	+ o_p\bigg(\sum_{i=1}^mn_i^\ell\bigg)
	+O_p(p)\\
	=&~ \sum_{i=1}^m \bm{\beta}_0(\alpha_0\setminus\alpha)'\bm{X}_i(\alpha_0\setminus\alpha)'
	\bm{H}_i^{-1}(\gamma,\bm\theta)\bm{X}_i(\alpha_0\setminus\alpha)\bm{\beta}_0(\alpha_0\setminus\alpha)\\
	&~  +2\sum_{i=1}^m\bm{\beta}_0(\alpha_0\setminus\alpha)'\bm{X}_i(\alpha_0\setminus\alpha)'
	\bm{H}_i^{-1}(\gamma,\bm\theta)(\bm{Z}_i(\gamma_0)\bm{b}_i(\gamma_0) +\bm{\epsilon}_i)\\
	&~  +\sum_{i=1}^m\bigg(\sum_{k\in\gamma_0}\bm{z}_{i,k}b_{i,k}+\bm{\epsilon}_i\bigg)'\bm{H}_i^{-1}(\gamma,\bm\theta)
	\bigg(\sum_{k\in\gamma_0}\bm{z}_{i,k}b_{i,k}+\bm{\epsilon}_i\bigg)\\
	&~  +o_p\bigg(\sum_{i=1}^mn_i^\xi\bigg)
	+o_p\bigg(\sum_{i=1}^mn_i^\ell\bigg)
	+O_p(p)
	\end{split}
	\end{align*}
	\begin{align*}
	\begin{split}
	=&~\sum_{i=1}^m\bm\epsilon_i'\bm\epsilon_i
	+\sum_{i=1}^m\sum_{j\in\alpha_0\setminus\alpha}\beta_{j,0}^2d_{i,j}n_i^{\xi}
	+\sum_{i=1}^m\sum_{k\in\gamma_0\setminus\gamma}b_{i,k}^2c_{i,k}n_i^\ell\\
	&~
	+o_p\bigg(\sum_{k,k^*\in\gamma\cap\gamma_0}\frac{m}{\theta_k\theta_{k^*}}\bigg)
	+O_p\bigg(\sum_{k\in\gamma\cap\gamma_0}\frac{m}{\theta_k}\bigg)\\
	&~    +o_p\bigg(\sum_{i=1}^mn_i^\xi\bigg)
	+o_p\bigg(\sum_{i=1}^mn_i^\ell\bigg)
	+O_p(p+mq)
	\end{split}
	\end{align*}
	
	\noindent uniformly over $\bm\theta\in\Theta_\gamma$,
	where the last equality follows from \eqref{cauchy} and Lemmas
	\ref{appendix:lemma:z x}--\ref{appendix:lemma:epsilon}.
	Hence by (\ref{partial:v}), we have, for $v^2\in(0,\infty)$,
	\begin{align*}
	\begin{split}
	v^4&\bigg\{\frac{\partial}{\partial v^2}\{-2\log L(\bm\theta,v^2;\alpha,\gamma)\}\bigg\}\\
	=&~ N\bigg(v^2-\frac{\bm\epsilon'\bm\epsilon}{N}
	+\frac{1}{N}\sum_{i=1}^m\sum_{j\in\alpha_0\setminus\alpha}\beta_{j,0}^2d_{i,j}n_i^{\xi}
	+\frac{1}{N}\sum_{i=1}^m\sum_{k\in\gamma_0\setminus\gamma}b_{i,k}^2c_{i,k}n_i^\ell\bigg)\\
	&~  +o_p\bigg(\sum_{i=1}^mn_i^\xi\bigg)+ o_p\bigg(\sum_{i=1}^mn_i^\ell\bigg)
	+O_p\bigg(\sum_{k,k^*\in\gamma\cap\gamma_0}\frac{m}{\theta_k\theta_{k^*}}\bigg)\\
	&~  +O_p\bigg(\sum_{k\in\gamma\cap\gamma_0}\frac{m}{\theta_k}\bigg)
	+ O_p(p+mq)
	\end{split}
	\end{align*}
	
	\noindent uniformly over $\bm\theta\in\Theta_\gamma$. This and Lemma
	\ref{appendix:prop:compact space} imply that
	\begin{align}
	\begin{split}
	\hat{v}^2(\alpha,\gamma)
	=&~ \frac{\bm\epsilon'\bm\epsilon}{N}
	+\frac{1}{N}\sum_{i=1}^m\sum_{j\in\alpha_0\setminus\alpha}\beta_{j,0}^2d_{i,j}n_i^{\xi}
	+\frac{1}{N}\sum_{i=1}^m\sum_{k\in\gamma_0\setminus\gamma}b_{i,k}^2c_{i,k}n_i^{\ell}\\
	&~  +o_p\bigg(\frac{1}{N}\sum_{i=1}^mn_i^{\xi}\bigg)
	+o_p\bigg(\frac{1}{N}\sum_{i=1}^mn_i^{\ell}\bigg)
	+O_p\Big(\frac{p+mq}{N}\Big).
	\end{split}
	\label{eq:v incorrect mixed}
	\end{align}
	
	\noindent Thus (\ref{appendix:thm:mle4:correct:eq1}) follows by applying the
	law of large numbers to $\bm\epsilon'\bm\epsilon/N$.
	In addition, if $(\xi,\ell)\in(0,1/2)\times (0,1/2)$, the asymptotic
	normality of $\hat{v}^2(\alpha,\gamma)$ follows by
	$p+mq=o(N^{1/2})$ and an application of the central limit theorem
	to $\bm\epsilon'\bm\epsilon/N$ in (\ref{eq:v incorrect mixed}).
	
	Next, we prove (\ref{appendix:thm:mle4:correct:eq2}), for
	$k\in\gamma\cap\gamma_0$, using \eqref{partial:theta2}.
	By \eqref{cauchy2} and Lemma \ref{appendix:lemma:xze} (i)--(iv),
	we have, for $k\in\gamma\cap\gamma_0$,
	\begin{align*}
	\theta_k&\bm{h}_{i,k}'\bm{H}^{-1}(\gamma,\bm\theta)\bm{M}(\alpha,\gamma;\bm\theta)
	\big(\bm{X}(\alpha_0\setminus\alpha)\bm{\beta}_0(\alpha_0\setminus\alpha)
	+\bm{Z}(\gamma_0)\bm{b}(\gamma_0)+\bm\epsilon\big)\\
	=&~ \theta_k\bm{h}_{i,k}'\bm{H}^{-1}(\gamma,\bm\theta)\bm{M}(\alpha,\gamma;\bm\theta)\\
	&~  \times\bigg(\bm{X}(\alpha_0\setminus\alpha)\bm{\beta}_0(\alpha_0\setminus\alpha)
	+\sum_{i^*=1}^m\sum_{k^*\in\gamma_0}b_{i^*,k^*}\bm{h}_{i^*,k^*}+\bm\epsilon\bigg)\\
	=&~ o_p\bigg(\frac{n_i^{(\xi-\ell)/2}\sum_{i^*=1}^mn_{i^*}^{(\xi+\ell)/2-\tau}}{\sum_{i=1}^mn_i^{\xi}}\bigg)
	+o_p(n_i^{(\xi-\ell)/2})+o_p(n_i^{-\ell/2})\\
	=&~	o_p\bigg(n_i^{(\xi-\ell)/2}\bigg(\frac{\sum_{i=1}^mn_i^{\ell}}{\sum_{i=1}^mn_i^{\xi}}\bigg)^{1/2}\bigg)+o_p(n_i^{(\xi-\ell)/2})+o_p(1)
	\end{align*}
	
	\noindent uniformly over $\bm\theta\in\Theta_\gamma$.
	This and (\ref{proof:thm:mle part3:eq0}) imply that for
	$k\in\gamma\cap\gamma_0$,
	\begin{align*}
	\begin{split}
	\theta_k & \bm{h}_{i,k}'\bm{H}^{-1}(\gamma,\bm\theta)(\bm{I}_{N}-\bm{M}(\alpha,\gamma;\bm\theta))\bm{y}\\
	=&~ \theta_k\bm{h}_{i,k}'\bm{H}^{-1}(\gamma,\bm\theta)(\bm{I}_{N}-\bm{M}(\alpha,\gamma;\bm\theta))\\
	&~	\times    \big(\bm{X}(\alpha_0\setminus\alpha)\bm{\beta}_0(\alpha_0\setminus\alpha)
	+\bm{Z}(\gamma_0)\bm{b}(\gamma_0)+\bm\epsilon\big)\\
	=&~ \theta_k\bm{h}_{i,k}'\bm{H}^{-1}(\gamma,\bm\theta)
	\big(\bm{X}(\alpha_0\setminus\alpha)\bm{\beta}_0(\alpha_0\setminus\alpha)
	+\bm{Z}(\gamma_0)\bm{b}(\gamma_0)+\bm\epsilon\big)\\
	&~  + o_p\bigg(n_i^{(\xi-\ell)/2}\bigg(\frac{\sum_{i=1}^mn_i^{\ell}}{\sum_{i=1}^mn_i^{\xi}}\bigg)^{1/2}\bigg)+o_p(n_i^{(\xi-\ell)/2})+o_p(1)\\
	=&~ \theta_k\bm{z}_{i,k}'\bm{H}_i^{-1}(\gamma,\bm\theta)
	\bigg(\bm{X}_i(\alpha_0\setminus\alpha)\bm{\beta}_0(\alpha_0\setminus\alpha)
	+\sum_{k^*\in\gamma_0}\bm{z}_{i,k^*}b_{i,k^*}+\bm\epsilon_i\bigg)\\
	&~  +o_p\bigg(n_i^{(\xi-\ell)/2}\bigg(\frac{\sum_{i=1}^mn_i^{\ell}}{\sum_{i=1}^mn_i^{\xi}}\bigg)^{1/2}\bigg)+o_p(n_i^{(\xi-\ell)/2})+o_p(1)\\
	=&~ b_{i,k}
	+o_p\bigg(n_i^{(\xi-\ell)/2}\bigg(\frac{\sum_{i=1}^mn_i^{\ell}}{\sum_{i=1}^mn_i^{\xi}}\bigg)^{1/2}\bigg)+o_p(n_i^{(\xi-\ell)/2})+o_p(1)
	\end{split}
	\end{align*}
	
	\noindent uniformly over $\bm\theta\in\Theta_\gamma$,
	where the last equality follows from Lemma \ref{appendix:lemma:z x} (iii),
	Lemma \ref{appendix:lemma:z} (ii)--(iv), and Lemma
	\ref{appendix:lemma:epsilon} (i).
	It follows that for $k\in\gamma\cap\gamma_0$,
	\begin{align*}
	\theta_k^2&\{\bm{h}_{i,k}'\bm{H}^{-1}(\gamma,\bm\theta)(\bm{I}_{N}-\bm{M}(\alpha,\gamma;\bm\theta))\bm{y}\}^2\\
	=&~ b_{i,k}^2
	+ o_p\bigg(n_i^{\xi-\ell}\bigg(\frac{\sum_{i=1}^mn_i^{\ell}}{\sum_{i=1}^mn_i^{\xi}}\bigg)\bigg)+o_p(n_i^{\xi-\ell})+o_p(1)
	\end{align*}
	\noindent uniformly over $\bm\theta\in\Theta_\gamma$.
	Hence by Lemma \ref{appendix:lemma:z} (ii) and (\ref{partial:theta2}),
	we have, for $k\in\gamma\cap\gamma_0$,
	\begin{align*}
	\begin{split}
	\theta_k^2&\bigg\{\frac{\partial}{\partial \theta_k}\{-2\log L(\bm\theta,v^2;\alpha,\gamma)\}\bigg\}\\
	=&~ m\bigg(\theta_k -\frac{1}{m}\sum_{i=1}^m\frac{b_{i,k}^2}{v^2}\bigg)
	+o_p\bigg(\sum_{i=1}^m n_i^{\xi-\ell}\bigg(1+\frac{\sum_{i=1}^mn_i^{\ell}}{\sum_{i=1}^mn_i^{\xi}}\bigg)\bigg)+o_p(m)
	\end{split}
	\end{align*}
	
	\noindent uniformly over $\bm\theta\in\Theta_\gamma$. This
	implies that for $k\in\gamma\cap\gamma_0$,
	\begin{align*}
	\hat\theta_{k}(\alpha,\gamma)
	=&~ \frac{1}{m}\sum_{i=1}^m\frac{b_{i,k}^2}{\hat{v}^2(\alpha,\gamma)}
	+o_p\bigg(\frac{1}{m}\sum_{i=1}^m n_i^{\xi-\ell}\bigg(1+\frac{\sum_{i=1}^mn_i^{\ell}}{\sum_{i=1}^mn_i^{\xi}}\bigg)\bigg)+o_p(1).
	\end{align*}
	
	\noindent This proves \eqref{appendix:thm:mle4:correct:eq2}, for
	$k\in\gamma\cap\gamma_0$.
	
	It remains to prove (\ref{appendix:thm:mle4:correct:eq2}), for $k\in
	\gamma\setminus\gamma_0$.
	Let $\bm\theta^\dag$ be $\bm\theta$ except that
	$\{\theta_k:k\in\gamma\cap\gamma_0\}$ are replaced by
	$\{\hat{\theta}_k(\alpha,\gamma):k\in\gamma\cap\gamma_0\}$.
	By \eqref{cauchy2} and Lemma \ref{appendix:lemma:xze} (i)--(iv),
	we have, for $k\in\gamma\setminus\gamma_0$,
	\begin{align*}
	\theta_{k}&
	\bm{h}_{i,k}'\bm{H}^{-1}(\gamma,\bm{\theta}^\dag )\bm{M}(\alpha,\gamma;\bm{\theta}^\dag )
	\big(\bm{X}(\alpha_0\setminus\alpha)\bm{\beta}_0(\alpha_0\setminus\alpha)
	+\bm{Z}(\gamma_0)\bm{b}(\gamma_0)+\bm\epsilon\big)\\
	=&~ \theta_k\bm{h}_{i,k}'\bm{H}^{-1}(\gamma,\bm{\theta}^\dag )\bm{M}(\alpha,\gamma;\bm{\theta}^\dag)\\
	&~  \times\bigg(\bm{X}(\alpha_0\setminus\alpha)\bm{\beta}_0(\alpha_0\setminus\alpha)
	+\sum_{i^*=1}^m\sum_{k^*\in\gamma_0}b_{i^*,k^*}\bm{h}_{i^*,k^*}+\bm\epsilon\bigg)\\
	=&~  o_p\bigg(n_i^{(\xi-\ell)/2}\bigg(\frac{\sum_{i=1}^mn_i^{\ell}}{\sum_{i=1}^mn_i^{\xi}}\bigg)^{1/2}\bigg)
	+o_p(n_i^{(\xi-\ell)/2-\tau})+o_p(n_i^{-\ell/2})\\
	=&~	o_p\bigg(n_i^{(\xi-\ell)/2}\bigg(\frac{\sum_{i=1}^mn_i^{\ell}}{\sum_{i=1}^mn_i^{\xi}}\bigg)^{1/2}\bigg)
	+o_p(n_i^{(\xi-\ell)/2})+o_p(1)
	\end{align*}
	\noindent uniformly over $\bm{\theta}(\gamma\setminus\gamma_0) \in[0,\infty)^{q(\gamma\setminus\gamma_0)}$.
	This and (\ref{proof:thm:mle part3:eq0}) imply that for $k\in\gamma\setminus\gamma_0$,
	\begin{align*}
	\begin{split}
	\theta_{k}&\bm{h}_{i,k}'\bm{H}^{-1}(\gamma,\bm{\theta}^\dag )(\bm{I}_{N}-\bm{M}(\alpha,\gamma;\bm{\theta}^\dag ))\bm{y}\\
	=&~ \theta_{k}
	\bm{h}_{i,k}'\bm{H}^{-1}(\gamma,\bm{\theta}^\dag )(\bm{I}_{N}-\bm{M}(\alpha,\gamma;\bm{\theta}^\dag ))
	\big(\bm{X}(\alpha_0\setminus\alpha)\bm{\beta}_0(\alpha_0\setminus\alpha)\\
	~& +\bm{Z}(\gamma_0)\bm{b}(\gamma_0)+\bm\epsilon\big)\\
	=&~ \theta_{k}\bm{h}_{i,k}'\bm{H}^{-1}(\gamma,\bm{\theta}^\dag )
	\big(\bm{X}(\alpha_0\setminus\alpha)\bm{\beta}_0(\alpha_0\setminus\alpha)
	+\bm{Z}(\gamma_0)\bm{b}(\gamma_0)+\bm\epsilon\big)\\
	&~  
	+ o_p(n_i^{(\xi-\ell)/2}n_{\max}^{(\ell-\xi)/2})
	+o_p(n_i^{(\xi-\ell)/2})+o_p(1)\\
	=&~ \theta_{k}\bm{z}_{i,k}'\bm{H}_i^{-1}(\gamma,\bm{\theta}^\dag )
	\bigg(\bm{X}_i(\alpha_0\setminus\alpha)\bm{\beta}_0(\alpha_0\setminus\alpha)
	+\sum_{i^*=1}^m\sum_{k^*\in\gamma_0}b_{i^*,k^*}\bm{h}_{i^*,k^*}
	+\bm\epsilon_i\bigg)\\
	&~ +o_p\bigg(n_i^{(\xi-\ell)/2}\bigg(\frac{\sum_{i=1}^mn_i^{\ell}}{\sum_{i=1}^mn_i^{\xi}}\bigg)^{1/2}\bigg)
	+o_p(n_i^{(\xi-\ell)/2})+o_p(1)\\
	=&~  
	o_p\bigg(n_i^{(\xi-\ell)/2}\bigg(\frac{\sum_{i=1}^mn_i^{\ell}}{\sum_{i=1}^mn_i^{\xi}}\bigg)^{1/2}\bigg)
	+o_p(n_i^{(\xi-\ell)/2})+o_p(1)
	\end{split}
	\end{align*}
	
	\noindent uniformly over $\bm{\theta}(\gamma\setminus\gamma_0)
	\in[0,\infty)^{q(\gamma\setminus\gamma_0)}$,
	where the last equality follows from Lemma \ref{appendix:lemma:z x} (iii),
	Lemma \ref{appendix:lemma:z} (iii)--(iv), and Lemma
	\ref{appendix:lemma:epsilon} (i).
	Therefore,
	\begin{align*}
	\theta_{k}^2&\{\bm{h}_{i,k}'\bm{H}^{-1}(\gamma,\bm{\theta}^\dag )(\bm{I}_{N}-\bm{M}(\alpha,\gamma;\bm{\theta}^\dag ))\bm{y}\}^2\\
	=&~
	o_p\bigg(n_i^{\xi-\ell}\bigg(\frac{\sum_{i=1}^mn_i^{\ell}}{\sum_{i=1}^mn_i^{\xi}}\bigg)\bigg)
	+o_p(n_i^{\xi-\ell})+o_p(1)
	\end{align*}
	uniformly over $\bm{\theta}(\gamma\setminus\gamma_0)
	\in[0,\infty)^{q(\gamma\setminus\gamma_0)}$.
	Hence by Lemma \ref{appendix:lemma:z} (ii) and (\ref{partial:theta2}), we
	have
	for $k\in\gamma\setminus\gamma_0$,
	\begin{align*}
	\begin{split}
	\theta_k^2&\bigg\{\frac{\partial}{\partial \theta_k}\{-2\log L(\bm\theta^\dag,v^2;\alpha,\gamma)\}\bigg\}\\
	=&~ m\theta_k+
	o_p\bigg(\sum_{i=1}^mn_i^{\xi-\ell}\bigg(1+\frac{\sum_{i=1}^mn_i^{\ell}}{\sum_{i=1}^mn_i^{\xi}}\bigg)\bigg)+o_p(m)
	\end{split}
	\end{align*}
	
	\noindent uniformly over $\bm\theta(\gamma\setminus\gamma_0)\in[0,\infty)^{q(\gamma\setminus\gamma_0)}$.
	This implies that for $k\in\gamma\setminus\gamma_0$,
	\[
	\hat\theta_{k}(\alpha,\gamma)
	=o_p\bigg(\frac{1}{m}\sum_{i=1}^mn_i^{\xi-\ell}\bigg(1+\frac{\sum_{i=1}^mn_i^{\ell}}{\sum_{i=1}^mn_i^{\xi}}\bigg)\bigg)+o_p(1).
	\]
	This completes the proof of \eqref{appendix:thm:mle4:correct:eq2}.
	Thus the proof of Theorem~\ref{appendix:theorem:MLE 4} is complete.

	\subsection{Proof of Theorem \ref{appendix:theorem:MLE 2}}
	
	As with the proof of Theorem \ref{appendix:theorem:MLE},
	we shall focus on the asymptotic properties of $\hat{v}^2(\alpha,\gamma)$
	and $\{\hat{\theta}_k(\alpha,\gamma):k\in\gamma\}$,
	and derive them by solving the likelihood equations directly.
	
	We first prove (\ref{appendix:thm:mle2:correct:eq1}) using
	(\ref{partial:v}). Hence by \eqref{cauchy2},
	Lemma~\ref{appendix:lemma:xze}~(i)--(iii), Lemma~\ref{appendix:lemma:xze}~(v)--(vi), and
	Lemma~\ref{appendix:lemma:xze}~(viii), we have
	\begin{align*}
	\begin{split}
	(\bm{Z}(\gamma_0)&\bm{b}(\gamma_0)+\bm\epsilon)'
	\bm{H}^{-1}(\gamma,\bm\theta)\bm{M}(\alpha,\gamma;\bm\theta)
	(\bm{Z}(\gamma_0)\bm{b}(\gamma_0)+\bm\epsilon)\\
	=&~ \bigg(\sum_{i=1}^m\sum_{k\in\gamma_0}b_{i,k}\bm{h}_{i,k}+\bm\epsilon\bigg)'
	\bm{H}^{-1}(\gamma,\bm\theta)\bm{M}(\alpha,\gamma;\bm\theta)
	\bigg(\sum_{i=1}^m\sum_{k\in\gamma_0}b_{i,k}\bm{h}_{i,k}+\bm\epsilon\bigg)\\
	=&~ o_p\bigg(\sum_{i=1}^m n_i^{\ell-\tau}\bigg)
	+o_p\bigg(\sum_{i=1}^mn_i^{\ell/2}\bigg)
	+O_p(p)\\
	=&~ o_p\bigg(\sum_{i=1}^mn_i^\ell\bigg)+ O_p(p)
	\end{split}
	\end{align*}
	
	\noindent uniformly over $\bm\theta\in\Theta_\gamma$.
	This and (\ref{proof:prop:compact space:eq2}) imply
	\begin{align*}
	\begin{split}
	\bm{y}'&\bm{H}^{-1}(\gamma,\bm\theta)(\bm{I}_{N}-\bm{M}(\alpha,\gamma;\bm\theta))\bm{y}\\
	=&~ (\bm{Z}(\gamma_0)\bm{b}(\gamma_0)+\bm\epsilon)'\bm{H}^{-1}(\gamma,\bm\theta)
	(\bm{I}_{N}-\bm{M}(\alpha,\gamma;\bm\theta))(\bm{Z}(\gamma_0)\bm{b}(\gamma_0)+\bm\epsilon)\\
	=&~ (\bm{Z}(\gamma_0)\bm{b}(\gamma_0)+\bm\epsilon)'\bm{H}^{-1}(\gamma,\bm\theta)
	(\bm{Z}(\gamma_0)\bm{b}(\gamma_0)+\bm\epsilon)
	+o_p\bigg(\sum_{i=1}^mn_i^\ell\bigg)
	+O_p(p)\\
	=&~ \sum_{i=1}^m (\bm{Z}_i(\gamma_0)\bm{b}_i(\gamma_0)+\bm\epsilon_i)'
	\bm{H}_i^{-1}(\gamma,\bm\theta)(\bm{Z}_i(\gamma_0)\bm{b}_i(\gamma_0)+\bm\epsilon_i)
	+o_p\bigg(\sum_{i=1}^mn_i^\ell\bigg)\\
	&~  +O_p(p)\\
	=&~ \sum_{i=1}^m\bm\epsilon_i'\bm\epsilon_i
	+ \sum_{i=1}^m\sum_{k\in\gamma_0\setminus\gamma}b_{i,k}^2c_{i,k}n_i^\ell
	+o_p\bigg(\sum_{i=1}^mn_i^\ell\bigg)
	+O_p\bigg(\sum_{k\in\gamma\cap\gamma_0}\frac{m}{\theta_k}\bigg)\\
	&~  +o_p\bigg(\sum_{k,k^*\in\gamma\cap\gamma_0}\frac{m}{\theta_k\theta_{k^*}}\bigg)
	+O_p(p+mq)
	\end{split}
	\end{align*}
	
	\noindent uniformly over $\bm\theta\in\Theta_\gamma$, where the last
	equality follows from
	Lemma \ref{appendix:lemma:z},
	Lemma \ref{appendix:lemma:epsilon} (i)--(ii), and 
	Lemma~\ref{appendix:lemma:epsilon}~(iv).
	Hence by (\ref{partial:v}), we have, for $v^2\in(0,\infty)$,
	\begin{align*}
	\begin{split}
	v^4&~\bigg\{\frac{\partial}{\partial v^2}\{-2\log L(\bm\theta,v^2;\alpha,\gamma)\}\bigg\}\\
	=&~ N\bigg(v^2-\frac{\bm\epsilon'\bm\epsilon}{N}
	+\frac{1}{N}\sum_{k\in\gamma_0\setminus\gamma}b_{i,k}^2c_{i,k}n_i^\ell\bigg)
	+ o_p\bigg(\sum_{i=1}^mn_i^\ell\bigg)\\
	&~
	+O_p\bigg(\sum_{k\in\gamma\cap\gamma_0}\frac{m}{\theta_k}\bigg)	+o_p\bigg(\sum_{k,k^*\in\gamma\cap\gamma_0}
	\frac{m}{\theta_k\theta_{k^*}}\bigg)
	+ O_p(p+mq)
	\end{split}
	\end{align*}
	
	\noindent uniformly over $\bm\theta\in\Theta_\gamma$. This and Lemma
	\ref{appendix:prop:compact space} imply that for
	$(\xi,\ell)\in(0,1]\times(0,1]$,
	\begin{align}
	\begin{split}
	\hat{v}^2(\alpha,\gamma)
	=&~ \frac{\bm\epsilon'\bm\epsilon}{N}
	+\frac{1}{N}\sum_{i=1}^m\sum_{k\in\gamma_0\setminus\gamma}b_{i,k}^2c_{i,k}n_i^{\ell}\\
	&~  +o_p\bigg(\frac{1}{N}\sum_{i=1}^mn_i^{\ell}\bigg)
	+O_p\bigg(\frac{p+mq}{N}\bigg).
	\end{split}
	\label{eq:v incorrect random}
	\end{align}
	
	\noindent Thus (\ref{appendix:thm:mle2:correct:eq1}) follows by applying the
	law of large numbers to $\bm\epsilon'\bm\epsilon/N$.
	In addition, if $\ell\in(0,1/2)$, the asymptotic normality of
	$\hat{v}^2(\alpha,\gamma)$ follows by
	$p+mq=o(N^{1/2})$ and an application of the central limit theorem
	to $\bm\epsilon'\bm\epsilon/N$ in (\ref{eq:v incorrect random}).
	
	Next, we prove (\ref{appendix:thm:mle2:correct:eq2}), for
	$k\in\gamma\cap\gamma_0$, using (\ref{partial:theta2}).
	By \eqref{cauchy2} and Lemma \ref{appendix:lemma:xze} (i)--(iii), we have, for
	$k\in\gamma\cap\gamma_0$,
	\begin{align*}
	\theta_k
	\bm{h}_{i,k}'&\bm{H}^{-1}(\gamma,\bm\theta)\bm{M}(\alpha,\gamma;\bm\theta)(\bm{Z}(\gamma_0)\bm{b}(\gamma_0)+\bm\epsilon)\\
	=&~ \theta_k\bm{h}_{i,k}'\bm{H}^{-1}(\gamma,\bm\theta)\bm{M}(\alpha,\gamma;\bm\theta)
	\bigg(\sum_{i^*=1}^m\sum_{k^*\in\gamma_0}b_{i^*,k^*}\bm{h}_{i^*,k^*}+\bm\epsilon\bigg)\\
	=&~ o_p\bigg(n_i^{(\xi-\ell)/2}\bigg(\frac{\sum_{i=1}^mn_i^{\ell}}{\sum_{i=1}^mn_i^{\xi}}\bigg)^{1/2}\bigg)
	+ o_p(n_i^{-\ell/2})\\
	=&~	o_p\bigg(n_i^{(\xi-\ell)/2}\bigg(\frac{\sum_{i=1}^mn_i^{\ell}}{\sum_{i=1}^mn_i^{\xi}}\bigg)^{1/2}\bigg) +o_p(1)
	\end{align*}
	
	\noindent uniformly over $\bm\theta\in\Theta_\gamma$. This and
	(\ref{proof:prop:compact space:eq2}) imply that for
	$k\in\gamma\cap\gamma_0$,
	\begin{align*}
	\begin{split}
	\theta_k&
	\bm{h}_{i,k}'\bm{H}^{-1}(\gamma,\bm\theta)(\bm{I}_{N}-\bm{M}(\alpha,\gamma;\bm\theta))\bm{y}\\
	=&~ \theta_k
	\bm{h}_{i,k}'\bm{H}^{-1}(\gamma,\bm\theta)(\bm{I}_{N}-\bm{M}(\alpha,\gamma;\bm\theta))
	(\bm{Z}(\gamma_0)\bm{b}(\gamma_0)+\bm\epsilon)\\
	=&~ \theta_k\bm{h}_{i,k}'\bm{H}^{-1}(\gamma,\bm\theta)
	(\bm{Z}(\gamma_0)\bm{b}(\gamma_0)+\bm\epsilon)
	+ o_p\bigg(n_i^{(\xi-\ell)/2}\bigg(\frac{\sum_{i=1}^mn_i^{\ell}}{\sum_{i=1}^mn_i^{\xi}}\bigg)^{1/2}\bigg)+o_p(1)\\
	=&~ \theta_k\bm{z}_{i,k}'\bm{H}_i^{-1}(\gamma,\bm\theta)
	\bigg(\sum_{k^*\in\gamma_0}\bm{z}_{i,k^*}b_{i,k^*}+\bm\epsilon_i\bigg)\\
	&~	+ o_p\bigg(n_i^{(\xi-\ell)/2}\bigg(\frac{\sum_{i=1}^mn_i^{\ell}}{\sum_{i=1}^mn_i^{\xi}}\bigg)^{1/2}\bigg)+o_p(1)\\
	=&~ b_{i,k}+o_p\bigg(n_i^{(\xi-\ell)/2}\bigg(\frac{\sum_{i=1}^mn_i^{\ell}}{\sum_{i=1}^mn_i^{\xi}}\bigg)^{1/2}\bigg)+o_p(1)
	\end{split}
	\end{align*}
	
	\noindent uniformly over $\bm\theta\in\Theta_\gamma$, where the last
	equality follows from Lemma \ref{appendix:lemma:z} (ii)--(iv) and
	Lemma \ref{appendix:lemma:epsilon} (i). Hence, for $k\in\gamma\cap\gamma_0$,
	\begin{align*}
	\theta_k^2&\{\bm{h}_{i,k}'\bm{H}^{-1}(\gamma,\bm\theta)(\bm{I}_{N}-\bm{M}(\alpha,\gamma;\bm\theta))\bm{y}\}^2\\
	=&~ b_{i,k}^2
	+o_p\bigg(n_i^{\xi-\ell}\bigg(\frac{\sum_{i=1}^mn_i^{\ell}}{\sum_{i=1}^mn_i^{\xi}}\bigg)\bigg)+o_p(1)
	\end{align*}
	\noindent uniformly over $\bm\theta\in\Theta_\gamma$. Hence by Lemma
	\ref{appendix:lemma:z} (ii) and (\ref{partial:theta2}),
	we have, for $k\in\gamma\cap\gamma_0$,
	\begin{align*}
	\begin{split}
	\theta_k^2&\bigg\{\frac{\partial}{\partial \theta_k}\{-2\log L(\bm\theta,v^2;\alpha,\gamma)\}\bigg\}\\
	=&~ m\bigg(\theta_k -\frac{1}{m}\sum_{i=1}^m\frac{b_{i,k}^2}{v^2}\bigg)
	+o_p\bigg(\sum_{i=1}^mn_i^{\xi-\ell}\bigg(\frac{\sum_{i=1}^mn_i^{\ell}}{\sum_{i=1}^mn_i^{\xi}}\bigg)\bigg)+o_p(m)
	\end{split}
	\end{align*}
	\noindent uniformly over $\bm\theta\in\Theta_\gamma$. Hence we have,
	for $k\in\gamma\cap\gamma_0$,
	\begin{align*}
	\hat\theta_{k}(\alpha,\gamma)
	=&~ \frac{1}{m}\sum_{i=1}^m\frac{b_{i,k}^2}{\hat{v}^2(\alpha,\gamma)}
	+o_p\bigg(\frac{1}{m}\sum_{i=1}^mn_i^{\xi-\ell}\bigg(\frac{\sum_{i=1}^mn_i^{\ell}}{\sum_{i=1}^mn_i^{\xi}}\bigg)\bigg)+o_p(1).
	\end{align*}
	\noindent This completes the proof of
	(\ref{appendix:thm:mle2:correct:eq2}), for $k\in\gamma\cap\gamma_0$.
	
	It remains to prove (\ref{appendix:thm:mle2:correct:eq2}), for $k\in\gamma\setminus\gamma_0$.
	Let $\bm\theta^\dag$ be $\bm\theta$ except that $\{\theta_k:k\in\gamma\cap\gamma_0\}$
	are replaced by $\{\hat{\theta}_k(\alpha,\gamma):k\in\gamma\cap\gamma_0\}$.
	By \eqref{cauchy2} and Lemma \ref{appendix:lemma:xze} (i)--(iii), we have, for
	$k\in\gamma\setminus\gamma_0$,
	\begin{align*}
	\begin{split}
	\theta_k&\bm{h}_{i,k}'\bm{H}^{-1}(\gamma,\bm{\theta}^\dag )\bm{M}(\alpha,\gamma;\bm{\theta}^\dag )
	(\bm{Z}(\gamma_0)\bm{b}(\gamma_0)+\bm\epsilon)\\
	=&~ \theta_k\bm{h}_{i,k}'\bm{H}^{-1}(\gamma,\bm{\theta}^\dag )\bm{M}(\alpha,\gamma;\bm{\theta}^\dag )
	\bigg(\sum_{i^*=1}^m\sum_{k^*\in\gamma_0}b_{i^*,k^*}\bm{h}_{i^*,k^*}+\bm\epsilon\bigg)\\
	=&~ o_p\bigg(n_i^{(\xi-\ell)/2}\bigg(\frac{\sum_{i=1}^mn_i^{\ell}}{\sum_{i=1}^mn_i^{\xi}}\bigg)^{1/2}\bigg) +o_p(1)
	\end{split}
	\end{align*}
	\noindent uniformly over $\bm{\theta}(\gamma\setminus\gamma_0) \in[0,\infty)^{q(\gamma\setminus\gamma_0)}$.
	This and (\ref{proof:prop:compact space:eq2}) imply that for
	$k\in\gamma\setminus\gamma_0$,
	\begin{align*}
	\begin{split}
	\theta_{k}&\bm{h}_{i,k}'\bm{H}^{-1}(\gamma,\bm{\theta}^\dag )
	(\bm{I}_{N}-\bm{M}(\alpha,\gamma;\bm{\theta}^\dag ))\bm{y}\\
	=&~ \theta_{k}\bm{h}_{i,k}'\bm{H}^{-1}(\gamma,\bm{\theta}^\dag )(\bm{I}_{N}-\bm{M}(\alpha,\gamma;\bm{\theta}^\dag ))
	(\bm{Z}(\gamma_0)\bm{b}(\gamma_0)+\bm\epsilon)\\
	=&~ \theta_{k}\bm{h}_{i,k}'\bm{H}^{-1}(\gamma,\bm{\theta}^\dag )
	(\bm{Z}(\gamma_0)\bm{b}(\gamma_0)+\bm\epsilon)
	+o_p\bigg(n_i^{(\xi-\ell)/2}\bigg(\frac{\sum_{i=1}^mn_i^{\ell}}{\sum_{i=1}^mn_i^{\xi}}\bigg)^{1/2}\bigg) +o_p(1)\\
	=&~ \theta_k\bm{z}_{i,k}'\bm{H}_i^{-1}(\gamma,\bm{\theta}^\dag )
	\bigg(\sum_{k^*\in\gamma_0}\bm{z}_{i,k^*}b_{i,k^*}+\bm\epsilon_i\bigg)\\
	&~+o_p\bigg(n_i^{(\xi-\ell)/2}\bigg(\frac{\sum_{i=1}^mn_i^{\ell}}{\sum_{i=1}^mn_i^{\xi}}\bigg)^{1/2}\bigg) +o_p(1)\\
	=&~
	o_p\bigg(n_i^{(\xi-\ell)/2}\bigg(\frac{\sum_{i=1}^mn_i^{\ell}}{\sum_{i=1}^mn_i^{\xi}}\bigg)^{1/2}\bigg) +o_p(1)
	\end{split}
	\end{align*}
	
	\noindent uniformly over $\bm{\theta}(\gamma\setminus\gamma_0) \in[0,\infty)^{q(\gamma\setminus\gamma_0)}$,
	where the last equality follows from Lemma \ref{appendix:lemma:z} (iii)--(iv)
	and Lemma \ref{appendix:lemma:epsilon} (i).
	Therefore,
	\begin{align*}
	\theta_{k}^2\{\bm{h}_{i,k}'\bm{H}^{-1}(\gamma,\bm{\theta}^\dag )
	(\bm{I}_{N}-\bm{M}(\alpha,\gamma;\bm{\theta}^\dag ))\bm{y}\}^2
	=&~ o_p\bigg(n_i^{\xi-\ell}\bigg(\frac{\sum_{i=1}^mn_i^{\ell}}{\sum_{i=1}^mn_i^{\xi}}\bigg)\bigg)+o_p(1)
	\end{align*}
	\noindent uniformly over $\bm{\theta}(\gamma\setminus\gamma_0) \in[0,\infty)^{q(\gamma\setminus\gamma_0)}$.
	Hence by Lemma \ref{appendix:lemma:z} (ii) and (\ref{partial:theta2}), we
	have,
	for $k\in\gamma\setminus\gamma_0$,
	\begin{align*}
	\begin{split}
	\theta_k^2&\bigg\{\frac{\partial}{\partial \theta_k}\{-2\log L(\bm\theta^\dag,v^2;\alpha,\gamma)\}\bigg\}
	= m\theta_k+o_p\bigg(\sum_{i=1}^mn_i^{\xi-\ell}\bigg(\frac{\sum_{i=1}^mn_i^{\ell}}{\sum_{i=1}^mn_i^{\xi}}\bigg)\bigg)+o_p(m)
	\end{split}
	\end{align*}
	\noindent uniformly over $\bm\theta(\gamma\setminus\gamma_0)\in[0,\infty)^{q(\gamma\setminus\gamma_0)}$.
	This implies that, for $k\in\gamma\setminus\gamma_0$,
	\begin{align*}
	\hat\theta_k(\alpha,\gamma) = o_p\bigg(\frac{1}{m}\sum_{i=1}^mn_i^{\xi-\ell}\bigg(\frac{\sum_{i=1}^mn_i^{\ell}}{\sum_{i=1}^mn_i^{\xi}}\bigg)\bigg)+o_p(1).
	\end{align*}
	\noindent This completes the proof of (\ref{appendix:thm:mle2:correct:eq2}).
	Hence the proof of Theorem \ref{appendix:theorem:MLE 2}
	is complete.

	\subsection{Proof of Theorem \ref{appendix:theorem:MLE 3}}
	
	As with the proof of Theorem \ref{appendix:theorem:MLE},
	we shall focus on the asymptotic properties of $\hat{v}^2(\alpha,\gamma)$
	and $\{\hat{\theta}_k(\alpha,\gamma):k\in\gamma\}$,
	and derive them by solving the likelihood equations directly.
	
	We first prove (\ref{appendix:thm:mle3:correct:eq1}) using (\ref{partial:v}).
	By Lemma \ref{appendix:lemma:xze} (i), Lemma \ref{appendix:lemma:xze} (iii)--(v),
	Lemma \ref{appendix:lemma:xze} (vii), and Lemma~\ref{appendix:lemma:xze}~(x), we have
	\begin{align*}
	\big(\bm{X}&(\alpha_0\setminus\alpha)\bm{\beta}_0(\alpha_0\setminus\alpha)
	+\bm{Z}(\gamma_0)\bm{b}(\gamma_0)+\bm\epsilon\big)'
	\bm{H}^{-1}(\gamma,\bm\theta)\bm{M}(\alpha,\gamma;\bm\theta)\\
	&~  \times\big(\bm{X}(\alpha_0\setminus\alpha)\bm{\beta}_0(\alpha_0\setminus\alpha)
	+\bm{Z}(\gamma_0)\bm{b}(\gamma_0)+\bm\epsilon\big)\\
	=&~\bigg(\bm{X}(\alpha_0\setminus\alpha)\bm{\beta}_0(\alpha_0\setminus\alpha)
	+\sum_{i=1}^m\sum_{k\in\gamma_0}b_{i,k}\bm{h}_{i,k}+\bm\epsilon\bigg)'
	\bm{H}^{-1}(\gamma,\bm\theta)\bm{M}(\alpha,\gamma;\bm\theta)\\
	&~  \times\bigg(\bm{X}(\alpha_0\setminus\alpha)\bm{\beta}_0(\alpha_0\setminus\alpha)
	+\sum_{i=1}^m\sum_{k\in\gamma_0}b_{i,k}\bm{h}_{i,k}+\bm\epsilon\bigg)\\
	=&~ o\bigg(\sum_{i=1}^n n_i^{\xi}\bigg)
	+o_p\bigg(\sum_{k,k^*\in\gamma_0}\frac{m}{\theta_k\theta_{k^*}}\bigg)
	+o_p\bigg(\sum_{k\in\gamma_0}\frac{m}{\theta_k}\bigg)
	+O_p(p)
	\end{align*}
	
	\noindent uniformly over $\bm\theta\in\Theta_\gamma$. This and
	(\ref{proof:thm:mle part3:eq0}) imply
	\begin{align*}
	\bm{y}'&\bm{H}^{-1}(\gamma,\bm\theta)(\bm{I}_{N}-\bm{M}(\alpha,\gamma;\bm\theta))\bm{y}\\
	=&~ \big(\bm{X}(\alpha_0\setminus\alpha)\bm{\beta}_0(\alpha_0\setminus\alpha)
	+\bm{Z}(\gamma_0)\bm{b}(\gamma_0)+\bm\epsilon\big)'
	\bm{H}^{-1}(\gamma,\bm\theta)\\
	&~  \times(\bm{I}_{N}-\bm{M}(\alpha,\gamma;\bm\theta))
	\big(\bm{X}(\alpha_0\setminus\alpha)\bm{\beta}_0(\alpha_0\setminus\alpha)
	+\bm{Z}(\gamma_0)\bm{b}(\gamma_0)+\bm\epsilon\big)
	\end{align*}
	\begin{align*}
	=&~ \big(\bm{X}(\alpha_0\setminus\alpha)\bm{\beta}_0(\alpha_0\setminus\alpha)
	+\bm{Z}(\gamma_0)\bm{b}(\gamma_0)+\bm\epsilon\big)'
	\bm{H}^{-1}(\gamma,\bm\theta)\\
	&~  \times\big(\bm{X}(\alpha_0\setminus\alpha)\bm{\beta}_0(\alpha_0\setminus\alpha)
	+\bm{Z}(\gamma_0)\bm{b}(\gamma_0)+\bm\epsilon\big)\\
	&~  +o\bigg(\sum_{i=1}^n n_i^{\xi}\bigg)
	+o_p\bigg(\sum_{k,k^*\in\gamma_0}\frac{m}{\theta_k\theta_{k^*}}\bigg)
	+o_p\bigg(\sum_{k\in\gamma_0}\frac{m}{\theta_k}\bigg)
	+O_p(p)\\
	=&~ \sum_{i=1}^m \bm{\beta}_0(\alpha_0\setminus\alpha)'\bm{X}_i(\alpha_0\setminus\alpha)'
	\bm{H}_i^{-1}(\gamma,\bm\theta)\bm{X}_i(\alpha_0\setminus\alpha)\bm{\beta}_0(\alpha_0\setminus\alpha)\\
	&~  +2\sum_{i=1}^m\bm{\beta}_0(\alpha_0\setminus\alpha)'\bm{X}_i(\alpha_0\setminus\alpha)'
	\bm{H}_i^{-1}(\gamma,\bm\theta)(\bm{Z}_i(\gamma_0)\bm{b}_i(\gamma_0) +\bm{\epsilon}_i)\\
	&~  +\sum_{i=1}^m(\bm{Z}_i(\gamma_0)\bm{b}_i(\gamma_0) +\bm{\epsilon}_i)'\bm{H}_i^{-1}(\gamma,\bm\theta)
	(\bm{Z}_i(\gamma_0)\bm{b}_i(\gamma_0) +\bm{\epsilon}_i)\\
	&~  +o\bigg(\sum_{i=1}^n n_i^{\xi}\bigg)
	+o_p\bigg(\sum_{k,k^*\in\gamma_0}\frac{m}{\theta_k\theta_{k^*}}\bigg)
	+o_p\bigg(\sum_{k\in\gamma_0}\frac{m}{\theta_k}\bigg)
	+O_p(p)\\
	=&~ \sum_{i=1}^m\bm\epsilon_i'\bm\epsilon_i +\sum_{i=1}^m\sum_{j\in\alpha_0\setminus\alpha}
	\beta_{j,0}^2d_{i,j}n_i^{\xi}+o_p\bigg(\sum_{i=1}^mn_i^{\xi}\bigg)
	+o_p\bigg(\sum_{k,k^*\in\gamma_0}\frac{m}{\theta_k\theta_{k^*}}\bigg)\\
	&~  +O_p\bigg(\sum_{k\in\gamma_0}\frac{m}{\theta_k}\bigg)
	+O_p(p+mq)
	\end{align*}
	
	\noindent uniformly over $\bm\theta\in\Theta_\gamma$, where the last
	equality follows from
	Lemma \ref{appendix:lemma:z} (ii)--(iv) and Lemma \ref{appendix:lemma:epsilon}.
	Hence by (\ref{partial:v}), we have, for $v^2\in(0,\infty)$,
	\begin{align*}
	\begin{split}
	v^4&\bigg\{\frac{\partial}{\partial v^2}\{-2\log L(\bm\theta,v^2;\alpha,\gamma)\}\bigg\}\\
	=&~ N\bigg(v^2-\frac{\bm\epsilon'\bm\epsilon}{N} +\frac{1}{N}\sum_{i=1}^m\sum_{j\in\alpha_0\setminus\alpha}
	\beta_{j,0}^2d_{i,j}n_i^{\xi}\bigg)
	+o_p\bigg(\sum_{i=1}^mn_i^{\xi}\bigg)\\
	&~
	+o_p\bigg(\sum_{k,k^*\in\gamma_0}\frac{m}{\theta_k\theta_{k^*}}\bigg)
	+O_p\bigg(\sum_{k\in\gamma_0}\frac{m}{\theta_k}\bigg)
	+O_p(p+mq)
	\end{split}
	\end{align*}
	
	\noindent uniformly over $\bm\theta\in\Theta_\gamma$. This and Lemma
	\ref{appendix:prop:compact space} imply
	that for $(\xi,\ell)\in(0,1]\times(0,1]$,
	\begin{align}
	\begin{split}
	\hat{v}^2(\alpha,\gamma)
	=&~ \frac{\bm\epsilon'\bm\epsilon}{N}
	+\frac{1}{N}\sum_{i=1}^m\sum_{j\in\alpha_0\setminus\alpha}\beta_{j,0}^2d_{i,j}n_i^{\xi} \\
	&~  +o_p\bigg(\frac{1}{N}\sum_{i=1}^mn_i^{\xi}\bigg)
	+O_p\bigg(\frac{p+mq}{N}\bigg).
	\end{split}
	\label{eq:v incorrect fixed}
	\end{align}
	
	\noindent Thus (\ref{appendix:thm:mle3:correct:eq1}) follows by applying the
	law of large numbers to $\bm\epsilon'\bm\epsilon/N$.
	In addition, if $\xi\in(0,1/2)$, the asymptotic normality of
	$\hat{v}^2(\alpha,\gamma)$ follows by
	$p+mq=o(N^{1/2})$ and an application of the central limit theorem
	to $\bm\epsilon'\bm\epsilon/N$ in (\ref{eq:v incorrect fixed}).
	
	Next, we prove (\ref{appendix:thm:mle3:correct:eq2}), for
	$k\in\gamma\cap\gamma_0$, using (\ref{partial:theta2}).
	By Lemma \ref{appendix:lemma:xze} (i) and Lemma \ref{appendix:lemma:xze}
	(iii)--(iv), we have, for $k\in\gamma\cap\gamma_0$,
	\begin{align*}
	\theta_k&\bm{h}_{i,k}'\bm{H}^{-1}(\gamma,\bm\theta)\bm{M}(\alpha,\gamma;\bm\theta)\big(
	\bm{X}(\alpha_0\setminus\alpha)\bm{\beta}_0(\alpha_0\setminus\alpha)
	+\bm{Z}(\gamma_0)\bm{b}(\gamma_0)+\bm\epsilon\big)\\
	=&~ \theta_k\bm{h}_{i,k}'\bm{H}^{-1}(\gamma,\bm\theta)\bm{M}(\alpha,\gamma;\bm\theta)\\
	&~  \times\bigg(\bm{X}(\alpha_0\setminus\alpha)\bm{\beta}_0(\alpha_0\setminus\alpha)
	+\sum_{i^*=1}^m\sum_{k^*\in\gamma_0}b_{i^*,k^*}\bm{h}_{i^*,k^*}
	+\bm\epsilon\bigg)\\
	=&~ o_p\bigg(\sum_{k^*\in\gamma_0}\frac{n_i^{(\xi-\ell)/2}}{\theta_{k^*}}\bigg)
	+o_p(n_i^{(\xi-\ell)/2-\tau})+o_p(n_i^{-\ell/2})\\
	=&~o_p\bigg(\sum_{k^*\in\gamma_0}\frac{n_i^{(\xi-\ell)/2}}{\theta_{k^*}}\bigg)+o_p(1)
	\end{align*}
	
	\noindent uniformly over $\bm\theta\in\Theta_\gamma$. This and
	(\ref{proof:thm:mle part3:eq0}) imply
	that for $k\in\gamma\cap\gamma_0$,
	\begin{align*}
	\begin{split}
	\theta_k\bm{h}_{i,k}'&\bm{H}^{-1}(\gamma,\bm\theta)(\bm{I}_{N}-\bm{M}(\alpha,\gamma;\bm\theta))\bm{y}\\
	=&~ \theta_k\bm{h}_{i,k}'\bm{H}^{-1}(\gamma,\bm\theta)
	(\bm{I}_{N}-\bm{M}(\alpha,\gamma;\bm\theta))\\
	&~\times
	\big(\bm{X}(\alpha_0\setminus\alpha)\bm{\beta}_0(\alpha_0\setminus\alpha)
	+\bm{Z}(\gamma_0)\bm{b}(\gamma_0)+\bm\epsilon\big)\\
	=&~ \theta_k\bm{h}_{i,k}'\bm{H}^{-1}(\gamma,\bm\theta)
	\big(\bm{X}(\alpha_0\setminus\alpha)\bm{\beta}_0(\alpha_0\setminus\alpha)
	+\bm{Z}(\gamma_0)\bm{b}(\gamma_0)+\bm\epsilon\big)\\
	&~  +o_p\bigg(\sum_{k^*\in\gamma_0}\frac{n_i^{(\xi-\ell)/2}}{\theta_{k^*}}\bigg)+o_p(1)\\
	=&~ \theta_k\bm{z}_{i,k}'\bm{H}_i^{-1}(\gamma,\bm\theta)
	\bigg(\bm{X}_i(\alpha_0\setminus\alpha)\bm{\beta}_0(\alpha_0\setminus\alpha)
	+\sum_{k^*\in\gamma_0}\bm{z}_{i,k^*}b_{i,k^*}+\bm\epsilon_i\bigg)\\
	&~  +o_p\bigg(\sum_{k^*\in\gamma_0}\frac{n_i^{(\xi-\ell)/2}}{\theta_{k^*}}\bigg)+o_p(1)
	\\
	=&~ b_{i,k}
	+o_p\bigg(\sum_{k^*\in\gamma_0}\frac{n_i^{(\xi-\ell)/2}}{\theta_{k^*}}\bigg)+o_p(1)
	\end{split}
	\end{align*}
	
	\noindent uniformly over $\bm\theta\in\Theta_\gamma$, where the last
	equality follows from Lemma \ref{appendix:lemma:z x} (iii), Lemma
	\ref{appendix:lemma:z} (ii)--(iii),
	and Lemma \ref{appendix:lemma:epsilon} (i).
	Hence, for $k\in\gamma\cap\gamma_0$,
	\begin{align*}
	\theta_k^2\{\bm{h}_{i,k}'\bm{H}^{-1}(\gamma,\bm\theta)(\bm{I}_{N}-\bm{M}(\alpha,\gamma;\bm\theta))\bm{y}\}^2
	=&~ b_{i,k}^2
	+o_p\bigg(\sum_{k,k^*\in\gamma_0}\frac{n_i^{\xi-\ell}}{\theta_k\theta_{k^*}}\bigg)+o_p(1)
	\end{align*}
	
	\noindent uniformly over $\bm\theta\in\Theta_\gamma$.
	Hence by Lemma \ref{appendix:lemma:z} (ii) and (\ref{partial:theta2}),
	we have, for $k\in\gamma\cap\gamma_0$,
	\begin{align*}
	\begin{split}
	\theta_k^2&\bigg\{\frac{\partial}{\partial \theta_k}\{-2\log L(\bm\theta,v^2;\alpha,\gamma)\}\bigg\}\\
	=&~ m\bigg(\theta_k -\frac{1}{m}\sum_{i=1}^m\frac{b_{i,k}^2}{v^2}\bigg)
	+o_p\bigg(\sum_{i=1}^m\sum_{k,k^*\in\gamma_0}\frac{n_i^{\xi-\ell}}{\theta_k\theta_{k^*}}\bigg)+o_p(1)
	\end{split}
	\end{align*}
	\noindent uniformly over $\bm\theta\in\Theta_\gamma$. This and Lemma
	\ref{appendix:prop:compact space} imply that for $k\in\gamma\cap\gamma_0$,
	\begin{align*}
	\hat\theta_{k}(\alpha,\gamma)
	=&~ \frac{1}{m}\sum_{i=1}^m\frac{b_{i,k}^2}{\hat{v}^2(\alpha,\gamma)}
	+o_p\bigg(\frac{1}{m}\sum_{i=1}^mn_i^{\xi-\ell}\bigg)+o_p(1).
	\end{align*}
	\noindent This completes the proof of (\ref{appendix:thm:mle3:correct:eq2})
	when $k\in\gamma\cap\gamma_0$.
	
	It remains to prove (\ref{appendix:thm:mle3:correct:eq2}), for
	$k\in\gamma\setminus\gamma_0$.
	Let $\bm\theta^\dag$ be $\bm\theta$ except that
	$\{\theta_k:k\in\gamma\cap\gamma_0\}$ are
	replaced by $\{\hat{\theta}_k(\alpha,\gamma):k\in\gamma\cap\gamma_0\}$.
	By Lemma \ref{appendix:lemma:xze} (i) and Lemma \ref{appendix:lemma:xze}
	(iii)--(iv),
	we have, for $k\in\gamma\setminus\gamma_0$,
	\begin{align*}
	\theta_{k}&\bm{h}_{i,k}'\bm{H}^{-1}(\gamma,\bm{\theta}^\dag )\bm{M}(\alpha,\gamma;\bm{\theta}^\dag )
	\big(\bm{X}(\alpha_0\setminus\alpha)\bm{\beta}_0(\alpha_0\setminus\alpha)
	+\bm{Z}(\gamma_0)\bm{b}(\gamma_0)+\bm\epsilon\big)\\
	=&~ \theta_{k}\bm{h}_{i,k}'\bm{H}^{-1}(\gamma,\bm{\theta}^\dag )\bm{M}(\alpha,\gamma;\bm{\theta}^\dag)\\
	&~  \times\bigg(\bm{X}(\alpha_0\setminus\alpha)\bm{\beta}_0(\alpha_0\setminus\alpha)
	+\sum_{i^*=1}^m\sum_{k^*\in\gamma_0}b_{i^*,k^*}\bm{h}_{i^*,k^*}+\bm\epsilon\bigg)\\
	=&~
	o_p(n_i^{(\xi-\ell)/2-\tau})+o_p(n_i^{-\ell/2})\\
	=&~o_p(n_i^{(\xi-\ell)/2})+o_p(1)
	\end{align*}
	
	\noindent uniformly over $\bm{\theta}(\gamma\setminus\gamma_0) \in[0,\infty)^{q(\gamma\setminus\gamma_0)}$.
	This and (\ref{proof:thm:mle part3:eq0}) imply that for
	$k\in\gamma\setminus\gamma_0$,
	\begin{align*}
	\begin{split}
	\theta_{k}&\bm{h}_{i,k}'\bm{H}^{-1}(\gamma,\bm{\theta}^\dag )(\bm{I}_{N}-\bm{M}(\alpha,\gamma;\bm{\theta}^\dag ))\bm{y}\\
	=&~ \theta_{k}
	\bm{h}_{i,k}'\bm{H}^{-1}(\gamma,\bm{\theta}^\dag )(\bm{I}_{N}-\bm{M}(\alpha,\gamma;\bm{\theta}^\dag ))\\
	&~\times\big(\bm{X}(\alpha_0\setminus\alpha)
	\bm{\beta}_0(\alpha_0\setminus\alpha)
	+\bm{Z}(\gamma_0)\bm{b}(\gamma_0)+\bm\epsilon\big)\\
	=&~ \theta_{k}\bm{h}_{i,k}'\bm{H}^{-1}(\gamma,\bm{\theta}^\dag )
	\big(\bm{X}(\alpha_0\setminus\alpha)\bm{\beta}_0(\alpha_0\setminus\alpha)
	+\bm{Z}(\gamma_0)\bm{b}(\gamma_0)+\bm\epsilon\big)\\
	&~    +o_p(n_i^{(\xi-\ell)/2})
	+o_p(1)\\
	=&~  \theta_{k}\bm{z}_{i,k}'\bm{H}_i^{-1}(\gamma,\bm{\theta}^\dag)
	\bigg(\bm{X}_i(\alpha_0\setminus\alpha)\bm{\beta}_0(\alpha_0\setminus\alpha)
	+\sum_{k^*\in\gamma_0}\bm{z}_{i,k^*}b_{i,k^*}+\bm\epsilon_i\bigg)\\
	&~  +o_p(n_i^{(\xi-\ell)/2})
	+o_p(1)\\
	=&~ o_p(n_i^{(\xi-\ell)/2})+o_p(1)
	\end{split}
	\end{align*}
	
	\noindent uniformly over $\bm{\theta}(\gamma\setminus\gamma_0) \in[0,\infty)^{q(\gamma\setminus\gamma_0)}$,
	where the last equality follows from Lemma \ref{appendix:lemma:z x} (iii),
	Lemma \ref{appendix:lemma:z} (iii), and Lemma \ref{appendix:lemma:epsilon} (i).
	Therefore,
	\begin{align*}
	\theta_{k}^2\{\bm{h}_{i,k}'\bm{H}^{-1}(\gamma,\bm{\theta}^\dag )
	(\bm{I}_{N}-\bm{M}(\alpha,\gamma;\bm{\theta}^\dag ))\bm{y}\}^2
	=&~o_p(n_i^{\xi-\ell})+o_p(1)
	\end{align*}
	\noindent uniformly over $\bm{\theta}(\gamma\setminus\gamma_0) \in[0,\infty)^{q(\gamma\setminus\gamma_0)}$.
	Hence by Lemma \ref{appendix:lemma:z} (ii) and (\ref{partial:theta2}), we
	have,
	for $k\in\gamma\setminus\gamma_0$,
	\begin{align*}
	\begin{split}
	\theta_k^2\bigg\{\frac{\partial}{\partial \theta_k}\{-2\log L(\bm\theta^\dag,v^2;\alpha,\gamma)\}\bigg\}
	=&~ m\theta_k + o_p\bigg(\sum_{i=1}^mn_i^{\xi-\ell}\bigg)
	+o_p(m)
	\end{split}
	\end{align*}
	\noindent uniformly over $\bm\theta(\gamma\setminus\gamma_0)\in[0,\infty)^{q(\gamma\setminus\gamma_0)}$.
	This and Lemma \ref{appendix:prop:compact space} imply that for
	$k\in\gamma\setminus\gamma_0$,
	\begin{align*}
	\hat\theta_{k}(\alpha,\gamma) = o_p\bigg(\frac{1}{m}\sum_{i=1}^mn_i^{\xi-\ell}\bigg)
	+o_p(1).
	\end{align*}
	\noindent This completes the proof of (\ref{appendix:thm:mle3:correct:eq2}),
	for $k\in\gamma\setminus\gamma_0$.
	Hence the proof of Theorem \ref{appendix:theorem:MLE 3} is complete.
	
	\setcounter{equation}{0}
	\section{Proofs of Auxiliary Lemmas}
	\label{appendix:lemma proofs}
	\subsection{Proof of Lemma \ref{appendix:lemma:z x}}
	Let $\bm{z}_{i,(s)}$; $s=1,\dots,q(\gamma)$ be the $s$-th
	column of $\bm{Z}_i(\gamma)$ and $\bm{H}_{i,t}(\gamma,\bm\theta)$ defined in
	(\ref{proof:lem:z:eq0}). For Lemma~\ref{appendix:lemma:z x}~(i)--(ii) to
	hold, it suffices to prove
	that for $k\notin\gamma$ and $j,j^*=1,\dots,p$,
	\begin{align}
	\bm{x}_{i,j}'\bm{H}_{i,t}^{-1}(\gamma,\bm\theta)\bm{x}_{i,j}
	=&~ d_{i,j}n_i^\xi + o(n_i^\xi) + o(tn_i^{\xi-2\tau}),
	\label{proof:lem:zx:eq1}\\
	\bm{x}_{i,j}'\bm{H}_{i,t}^{-1}(\gamma,\bm\theta)\bm{x}_{i,j^*}
	=&~ o(n_i^{\xi-\tau}) + o(tn_i^{\xi-2\tau}),
	\label{proof:lem:zx:eq2}\\
	\bm{x}_{i,j}'\bm{H}_{i,t}^{-1}(\gamma,\bm\theta)\bm{z}_{i,k}
	=&~ o(n_i^{(\xi+\ell)/2-\tau}) + o(tn_i^{(\xi+\ell)/2-2\tau})
	\label{proof:lem:zx:eq3}
	\end{align}
	
	\noindent uniformly over $\bm\theta\in[0,\infty)^{q(\gamma)}$. We
	prove (\ref{proof:lem:zx:eq1})--(\ref{proof:lem:zx:eq3}) by
	induction. For $j=1,\dots,p$ and $t=1$, by (\ref{eq:inv decompose}) and
	(A1)--(A3), we have
	\begin{align*}
	\bm{x}_{i,j}'\bm{H}_{i,1}^{-1}(\gamma,\bm\theta)\bm{x}_{i,j}
	=&~ \bm{x}_{i,j}'\bm{x}_{i,j} - \frac{\theta_{(1)}
		\bm{x}_{i,j}'\bm{z}_{i,(1)}\bm{z}_{i,(1)}'\bm{x}_{i,j}}{1+\theta_{(1)}\bm{z}_{i,(1)}'\bm{z}_{i,(1)}}\\
	=&~ d_{i,j}n_i^\xi + o(n_i^\xi) + o(n_i^{\xi-2\tau})
	\end{align*}
	
	\noindent uniformly over $\bm\theta\in[0,\infty)^{q(\gamma)}$. For $j,j^*=1,\dots,p$,
	$j\neq j^*$ and $t=1$, by (\ref{eq:inv decompose}) and (A1)--(A3), we have
	\begin{align*}
	\bm{x}_{i,j}'\bm{H}_{i,1}^{-1}(\gamma,\bm\theta)\bm{x}_{i,j^*}
	=&~ \bm{x}_{i,j}'\bm{x}_{i,j^*} - \frac{\theta_{(1)}\bm{x}_{i,j}'\bm{z}_{i,(1)}
		\bm{z}_{i,(1)}'\bm{x}_{i,j^*}}{1+\theta_{(1)}\bm{z}_{i,(1)}'\bm{z}_{i,(1)}}\\
	=&~ o(n_i^{\xi-\tau})+o(n_i^{\xi-2\tau})
	\end{align*}
	
	\noindent uniformly over $\bm\theta\in[0,\infty)^{q(\gamma)}$.
	For $j=1,\dots,p$, $k\notin\gamma$ and $t=1$, by (\ref{eq:inv decompose})
	and (A1)--(A3), we have
	\begin{align*}
	\bm{x}_{i,j}'\bm{H}_{i,1}^{-1}(\gamma,\bm\theta)\bm{z}_{i,k}
	=&~ \bm{x}_{i,j}'\bm{z}_{i,k} - \frac{\theta_{(1)}
		\bm{x}_{i,j}'\bm{z}_{i,(1)}\bm{z}_{i,(1)}'\bm{z}_{i,k}}{1+\theta_{(1)}\bm{z}_{i,(1)}'\bm{z}_{i,(1)}}\\
	=&~ o(n_i^{(\xi+\ell)/2-\tau})+o(n_i^{(\xi+\ell)/2-2\tau})
	\end{align*}
	
	\noindent uniformly over $\bm\theta\in[0,\infty)^{q(\gamma)}$. Now
	suppose that (\ref{proof:lem:zx:eq1})--(\ref{proof:lem:zx:eq3}) hold
	for $t=r$. Then for $j=1,\dots,p$ and $t=r+1$,
	by (\ref{eq:inv decompose}) and
	(\ref{proof:lem:zx:eq1})--(\ref{proof:lem:zx:eq3}) with $t=r$, and Lemma
	\ref{appendix:lemma:z} (i),
	we have
	\begin{align*}
	\bm{x}_{i,j}'&\bm{H}_{i,r+1}^{-1}(\gamma,\bm\theta)\bm{x}_{i,j}\\
	=&~ \bm{x}_{i,j}'\bm{H}_{i,r}^{-1}(\gamma,\bm\theta)\bm{x}_{i,j}
	- \frac{\theta_{(r+1)}\bm{x}_{i,j}'\bm{H}_{i,r}^{-1}(\gamma,\bm\theta)
		\bm{z}_{i,(r+1)}\bm{z}_{i,(r+1)}'\bm{H}_{i,r}^{-1}(\gamma,\bm\theta)\bm{x}_{i,j}}
	{1+\theta_{(r+1)}\bm{z}_{i,(r+1)}'\bm{H}_{i,r}^{-1}(\gamma,\bm\theta)\bm{z}_{i,(r+1)}}\\
	=&~ d_{i,j}n_i^\xi + o(n_i^\xi) +o(\{r+1\}n_i^{\xi-2\tau})
	\end{align*}
	
	\noindent uniformly over $\bm\theta\in[0,\infty)^{q(\gamma)}$. For $j,j^*=1,\dots,p$, $j\neq j^*$,
	and $t=r+1$, by (\ref{eq:inv decompose}) and
	(\ref{proof:lem:zx:eq1})--(\ref{proof:lem:zx:eq3}) with $t=r$,
	and Lemma \ref{appendix:lemma:z} (i), we have
	\begin{align*}
	\bm{x}_{i,j}'&\bm{H}_{i,r+1}^{-1}(\gamma,\bm\theta)\bm{x}_{i,j^*}\\
	=&~ \bm{x}_{i,j}'\bm{H}_{i,r}^{-1}(\gamma,\bm\theta)\bm{x}_{i,j^*}
	- \frac{\theta_{(r+1)}\bm{x}_{i,j}'\bm{H}_{i,r}^{-1}(\gamma,\bm\theta)
		\bm{z}_{i,(r+1)}\bm{z}_{i,(r+1)}'\bm{H}_{i,r}^{-1}(\gamma,\bm\theta)\bm{x}_{i,j^*}}
	{1+\theta_{(r+1)}\bm{z}_{i,(r+1)}'\bm{H}_{i,r}^{-1}(\gamma,\bm\theta)\bm{z}_{i,(r+1)}}\\
	=&~ o(n_i^{\xi-\tau}) + o(\{r+1\}n_i^{\xi-2\tau})
	\end{align*}
	
	\noindent uniformly over $\bm\theta\in[0,\infty)^{q(\gamma)}$.
	For $j,j^*=1,\dots,p$, $k\notin\gamma$,
	and $t=r+1$, by (\ref{eq:inv decompose}) and
	(\ref{proof:lem:zx:eq1})--(\ref{proof:lem:zx:eq3}) with $t=r$, and Lemma
	\ref{appendix:lemma:z} (i),
	we have
	\begin{align*}
	\bm{x}_{i,j}'&\bm{H}_{i,r+1}^{-1}(\gamma,\bm\theta)\bm{z}_{i,k}\\
	=&~ \bm{x}_{i,j}'\bm{H}_{i,r}^{-1}(\gamma,\bm\theta)\bm{z}_{i,k}
	- \frac{\theta_{(r+1)}\bm{x}_{i,j}'\bm{H}_{i,r}^{-1}(\gamma,\bm\theta)
		\bm{z}_{i,(r+1)}\bm{z}_{i,(r+1)}'\bm{H}_{i,r}^{-1}(\gamma,\bm\theta)
		\bm{z}_{i,k}}{1+\theta_{(r+1)}\bm{z}_{i,(r+1)}'\bm{H}_{i,r}^{-1}(\gamma,\bm\theta)\bm{z}_{i,(r+1)}}\\
	=&~ o(n_i^{(\xi+\ell)/2-\tau}) + o(\{r+1\}n_i^{(\xi+\ell)/2-2\tau})
	\end{align*}
	
	\noindent uniformly over $\bm\theta\in[0,\infty)^{q(\gamma)}$. This
	completes the proofs of (\ref{proof:lem:zx:eq1})--(\ref{proof:lem:zx:eq3}).
	Hence the proofs of Lemma~\ref{appendix:lemma:z x}~(i)--(ii) are complete.
	
	We finally prove Lemma \ref{appendix:lemma:z x} (iii). Without loss of
	generality, we assume that $q(\gamma)=q$, $t=q$,
	and $k=(q)$. Then by (\ref{eq:inv decompose}),
	\begin{align*}
	\theta_{(q)}\bm{x}_{i,j}'\bm{H}_{i,q}^{-1}(\gamma,\bm\theta)\bm{z}_{i,(q)}
	=&~ \theta_{(q)}\bigg\{\bm{x}_{i,j}'\bm{H}_{i,q-1}^{-1}(\gamma,\bm\theta)\bm{z}_{i,(q)}\\
	&~  -\frac{\theta_{(q)}\bm{x}_{i,j}'\bm{H}_{i,q-1}^{-1}(\gamma,\bm\theta)\bm{z}_{i,(q)}\bm{z}_{i,(q)}'
		\bm{H}_{i,q-1}^{-1}(\gamma,\bm\theta)\bm{z}_{i,q}}
	{1+\theta_{(q)} \bm{z}_{i,(q)}'\bm{H}_{i,q-1}^{-1}(\gamma,\bm\theta)\bm{z}_{i,(q)}}\bigg\}\\
	=&~ \frac{\theta_{(q)}\bm{x}_{i,j}'\bm{H}_{i,q-1}^{-1}(\gamma,\bm\theta)\bm{z}_{i,(q)}}
	{1+\theta_{(q)} \bm{z}_{i,(q)}'\bm{H}_{i,q-1}^{-1}(\gamma,\bm\theta)\bm{z}_{i,(q)}},
	\end{align*}
	
	\noindent where we note that $\theta_{(q)}$ can be arbitrarily small and the
	dominant term of the denominator of the last equation can be equal to (i)
	$\theta_{(q)}
	\bm{z}_{i,(q)}'\bm{H}_{i,q-1}^{-1}(\gamma,\bm\theta)\bm{z}_{i,(q)}$ or (ii)
	$1$. For the case of (i), $\theta_{(q)}n_i^{\ell}\rightarrow\infty$ by Lemma
	\ref{appendix:lemma:z} (i); hence, using Lemma \ref{appendix:lemma:z x}
	(ii) and Lemma \ref{appendix:lemma:z} (i), we have
	\begin{align*}
	\theta_{(q)}\bm{x}_{i,j}'\bm{H}_{i,q}^{-1}(\gamma,\bm\theta)\bm{z}_{i,(q)}
	=&~	o(n_i^{(\xi-\ell)/2-\tau}),
	\end{align*}
%	\noindent which also gives that
	\noindent and thus               % AMH; check
	\begin{align*}
	\bm{x}_{i,j}'\bm{H}_{i,q}^{-1}(\gamma,\bm\theta)\bm{z}_{i,(q)}
	=&~	o(n_i^{(\xi+\ell)/2-\tau}).
	\end{align*}
	\noindent For the case of (ii), $\theta_{(q)}=O(n_i^{-\ell})$ by Lemma
	\ref{appendix:lemma:z} (i); hence, using Lemma \ref{appendix:lemma:z} (i),
	we have
	\begin{align*}
	\theta_{(q)}\bm{x}_{i,j}'\bm{H}_{i,q}^{-1}(\gamma,\bm\theta)\bm{z}_{i,(q)}
	=&~	o(\theta_{(q)}n_i^{(\xi+\ell)/2-\tau}),
	\end{align*}
	\noindent which also gives the following two results:
	\begin{align*}
	\bm{x}_{i,j}'\bm{H}_{i,q}^{-1}(\gamma,\bm\theta)\bm{z}_{i,(q)}
	=&~	o(n_i^{(\xi+\ell)/2-\tau}),\\
	\theta_{(q)}\bm{x}_{i,j}'\bm{H}_{i,q}^{-1}(\gamma,\bm\theta)\bm{z}_{i,(q)}
	=&~	o(n_i^{(\xi-\ell)/2-\tau}).
	\end{align*}
	\noindent In conclusion, we have
	\begin{align}
	\begin{split}
	\theta_{(q)}\bm{x}_{i,j}'\bm{H}_{i,q}^{-1}(\gamma,\bm\theta)\bm{z}_{i,(q)}
	=&~	o(n_i^{(\xi-\ell)/2-\tau}),\\
	\bm{x}_{i,j}'\bm{H}_{i,q}^{-1}(\gamma,\bm\theta)\bm{z}_{i,(q)}
	=&~	o(n_i^{(\xi+\ell)/2-\tau})
	\end{split}
	\label{proof:revised one}
	\end{align}
	\noindent uniformly over $\bm\theta\in[0,\infty)^{q}$. This completes the proof.
	\subsection{Proof of Lemma \ref{appendix:lemma:z}}
	Let $\bm{z}_{i,(s)}$; $s=1,\dots,q(\gamma)$ be the $s$-th
	column of $\bm{Z}_i(\gamma)$ and $\bm{H}_{i,t}(\gamma,\bm\theta)$ defined in
	(\ref{proof:lem:z:eq0}). We first prove Lemma \ref{appendix:lemma:z} (i).
	By (\ref{proof:lem:z:eq0}), it suffices to prove that for $k\notin\gamma$,
	\begin{align}
	\bm{z}_{i,k}'\bm{H}_{i,t}^{-1}(\gamma,\bm\theta)\bm{z}_{i,k}
	=&~  c_{i,k}n_i^\ell +o(n_i^\ell) + o(tn_i^{\ell-2\tau}),
	\label{proof:lem:z:eq1}
	\end{align}
	\noindent and for $k,k^*\notin\gamma$ and $k\neq k^*$,
	\begin{align}
	\bm{z}_{i,k}'\bm{H}_{i,t}^{-1}(\gamma,\bm\theta)\bm{z}_{i,k^*}
	=&~  o(n_i^{\ell-\tau})+ o(tn_i^{\ell-2\tau})
	\label{proof:lem:z:eq2}
	\end{align}
	
	\noindent uniformly over $\bm\theta\in[0,\infty)^{q(\gamma)}$ by induction.
	For $t=1$ and $k\notin\gamma$, by (\ref{eq:inv decompose}) and (A2), we have
	\begin{align*}
	\bm{z}_{i,k}'\bm{H}_{i,1}^{-1}(\gamma,\bm\theta)\bm{z}_{i,k}
	=&~ \bm{z}_{i,k}'\bigg(\bm{I}_{n_i} - \frac{\theta_{(1)}\bm{z}_{i,(1)}\bm{z}_{i,(1)}'}
	{1+\theta_{(1)}\bm{z}_{i,(1)}'\bm{z}_{i,(1)}}\bigg)\bm{z}_{i,k}\\
	=&~ \bm{z}_{i,k}'\bm{z}_{i,k}- \frac{\theta_{(1)}\bm{z}_{i,k}'\bm{z}_{i,(1)}\bm{z}_{i,(1)}'
		\bm{z}_{i,k}}{1+\theta_{(1)}\bm{z}_{i,(1)}'\bm{z}_{i,(1)}}\\
	=&~ c_{i,k} n_i^\ell +o(n_i^\ell)+ o(n_i^{\ell-2\tau})
	\end{align*}
	
	\noindent uniformly over $\bm\theta\in[0,\infty)^{q(\gamma)}$. For
	$k,k^*\notin\gamma$ and $k\neq k^*$, by (\ref{eq:inv decompose}) and (A2),
	we have
	\begin{align*}
	\bm{z}_{i,k}'\bm{H}_{i,1}^{-1}(\gamma,\bm\theta)\bm{z}_{i,k^*}
	=&~ \bm{z}_{i,k}'\bm{z}_{i,k^*}
	- \frac{\theta_{(1)}\bm{z}_{i,k}'\bm{z}_{i,(1)}\bm{z}_{i,(1)}'
		\bm{z}_{i,k^*}}{1+\theta_{(1)}\bm{z}_{i,(1)}'\bm{z}_{i,(1)}}\\
	=&~ o(n_i^{\ell-\tau}) + o(n_i^{\ell-2\tau})
	\end{align*}
	
	\noindent uniformly over $\bm\theta\in[0,\infty)^{q(\gamma)}$. Now
	suppose that (\ref{proof:lem:z:eq1}) and (\ref{proof:lem:z:eq2})
	hold for $t=r$. Then for $k\notin\gamma$ and $t=r+1$, by (\ref{eq:inv decompose}),
	and (\ref{proof:lem:z:eq1}) and (\ref{proof:lem:z:eq2}) with $t=r$, we have
	\begin{align*}
	\bm{z}_{i,k}'\bm{H}_{i,r+1}^{-1}(\gamma,\bm\theta)\bm{z}_{i,k}
	=&~ \bm{z}_{i,k}'\bm{H}_{i,r}^{-1}(\gamma,\bm\theta)\bm{z}_{i,k}\\
	&~  - \frac{\theta_{(r+1)}\bm{z}_{i,k}'\bm{H}_{i,r}^{-1}(\gamma,\bm\theta)
		\bm{z}_{i,(r+1)}\bm{z}_{i,(r+1)}'\bm{H}_{i,r}^{-1}(\gamma,\bm\theta)\bm{z}_{i,k}}
	{1+\theta_{(r+1)}\bm{z}_{i,(r+1)}'\bm{H}_{i,r}^{-1}(\gamma,\bm\theta)\bm{z}_{i,(r+1)}}\\
	=&~ c_{i,k} n_i^\ell + o(n_i^\ell) + o(\{r+1\}n_i^{\ell-2\tau})
	\end{align*}
	
	\noindent uniformly over $\bm\theta\in[0,\infty)^{q(\gamma)}$.
	For $k,k^*\notin\gamma$ and $t=r+1$, by (\ref{eq:inv decompose}),
	and (\ref{proof:lem:z:eq1}) and (\ref{proof:lem:z:eq2}) with $t=r$, we have
	\begin{align*}
	\bm{z}_{i,k}'\bm{H}_{i,r+1}^{-1}(\gamma,\bm\theta)\bm{z}_{i,k^*}
	=&~ \bm{z}_{i,k}'\bm{H}_{i,r}^{-1}(\gamma,\bm\theta)\bm{z}_{i,k^*}\\
	&~  - \frac{\theta_{(r+1)}\bm{z}_{i,k}'\bm{H}_{i,r}^{-1}(\gamma,\bm\theta)
		\bm{z}_{i,(r+1)}\bm{z}_{i,(r+1)}'\bm{H}_{i,r}^{-1}(\gamma,\bm\theta)\bm{z}_{i,k^*}}
	{1+\theta_{(r+1)}\bm{z}_{i,(r+1)}'\bm{H}_{i,r}^{-1}(\gamma,\bm\theta)\bm{z}_{i,(r+1)}}\\
	=&~ o(n_i^{\ell-\tau}) + o(\{r+1\}n_i^{\ell-2\tau})
	\end{align*}
	
	\noindent uniformly over $\bm\theta\in[0,\infty)^{q(\gamma)}$. This
	completes the proof of (\ref{proof:lem:z:eq1}) and
	(\ref{proof:lem:z:eq2}). Hence Lemma \ref{appendix:lemma:z} (i) follows from (\ref{proof:lem:z:eq1}),
	(\ref{proof:lem:z:eq2}) with $t=q(\gamma)$ and $q=o(n_{\min}^\tau)$.
	This completes the proof of Lemma \ref{appendix:lemma:z} (i).
	
	We now prove Lemma \ref{appendix:lemma:z} (ii). Without loss of
	generality, we assume that $q(\gamma)=q$ and
	$k=(q)$. Then by Lemma \ref{appendix:lemma:z} (i) and (\ref{eq:inv decompose}),
	\begin{align*}
	\theta_{(q)}^2\bm{z}_{i,(q)}'\bm{H}_{i,q}^{-1}(\gamma,\bm\theta)\bm{z}_{i,(q)}
	=&~ \theta_{(q)}^2\bigg\{\bm{z}_{i,(q)}'\bm{H}_{i,q-1}^{-1}(\gamma,\bm\theta)\bm{z}_{i,(q)}\\
	&~  - \frac{\theta_{(q)}(\bm{z}_{i,(q)}'\bm{H}_{i,q-1}^{-1}
		(\gamma,\bm\theta)\bm{z}_{i,(q)})^2}
	{1+\theta_{(q)} \bm{z}_{i,(q)}'\bm{H}_{i,q-1}^{-1}(\gamma,\bm\theta)\bm{z}_{i,(q)}}\bigg\}\\
	=&~ \frac{\theta_{(q)}^2\bm{z}_{i,(q)}'\bm{H}_{i,q-1}^{-1}(\gamma,\bm\theta)\bm{z}_{i,(q)}}
	{1+\theta_{(q)} \bm{z}_{i,(q)}'\bm{H}_{i,q-1}^{-1}(\gamma,\bm\theta)\bm{z}_{i,(q)}}
	=   O(\theta_{(q)}^2n_i^{\ell})
	\end{align*}
	
	\noindent uniformly over $\bm\theta\in[0,\infty)^q$.
	Again, by Lemma \ref{appendix:lemma:z} (i), we have
	\begin{align*}
	\theta_{(q)}^2\bm{z}_{i,(q)}'\bm{H}_{i,q}^{-1}(\gamma,\bm\theta)\bm{z}_{i,(q)}
	=&~ \frac{\theta_{(q)}^2\bm{z}_{i,(q)}'\bm{H}_{i,q-1}^{-1}(\gamma,\bm\theta)\bm{z}_{i,(q)}}
	{1+\theta_{(q)} \bm{z}_{i,(q)}'\bm{H}_{i,q-1}^{-1}(\gamma,\bm\theta)\bm{z}_{i,(q)}}\\
	=&~ \theta_{(q)} - \frac{\theta_{(q)}}{1+\theta_{(q)} \bm{z}_{i,(q)}'\bm{H}_{i,q-1}^{-1}(\gamma,\bm\theta)\bm{z}_{i,(q)}}\\
	=&~ \theta_{(q)} + O(n_i^{-\ell})
	\end{align*}
	
	\noindent uniformly over $\bm\theta\in[0,\infty)^q$. This completes the
	proof of Lemma \ref{appendix:lemma:z} (ii).
	
	We now prove Lemma \ref{appendix:lemma:z} (iii). Without loss of
	generality, we assume that $q(\gamma)=q$, $k=(q)$, and $k^*=(q-1)$. Then by
	(\ref{eq:inv decompose}),
	\begin{align*}
	\theta_{(q)}&\theta_{(q-1)}\bm{z}_{i,(q)}'\bm{H}_{i,q}^{-1}(\gamma,\bm\theta)\bm{z}_{i,(q-1)}\\
	=&~ \theta_{(q)}\theta_{(q-1)}\bigg\{\bm{z}_{i,(q)}'\bm{H}_{i,q-1}^{-1}(\gamma,\bm\theta)\bm{z}_{i,(q-1)}\\
	&~  -\frac{\theta_{(q)}\bm{z}_{i,(q)}'\bm{H}_{i,q-1}^{-1}(\gamma,\bm\theta)
		\bm{z}_{i,(q)}\bm{z}_{i,(q)}'\bm{H}_{i,q-1}^{-1}(\gamma,\bm\theta)\bm{z}_{i,(q-1)}}
	{1+\theta_{(q)} \bm{z}_{i,(q)}'\bm{H}_{i,q-1}^{-1}(\gamma,\bm\theta)\bm{z}_{i,(q)}}\bigg\}\\
	=&~ \frac{\theta_{(q)}\theta_{(q-1)}\bm{z}_{i,(q)}'\bm{H}_{i,q-1}^{-1}(\gamma,\bm\theta)\bm{z}_{i,(q-1)}}
	{1+\theta_{(q)} \bm{z}_{i,(q)}'\bm{H}_{i,q-1}^{-1}(\gamma,\bm\theta)\bm{z}_{i,(q)}}\\
	=&~ \frac{\theta_{(q)}\theta_{(q-1)}\bm{z}_{i,(q)}'\bm{H}_{i,q-2}^{-1}(\gamma,\bm\theta)\bm{z}_{i,(q-1)}}
	{(1+\theta_{(q)} \bm{z}_{i,(q)}'\bm{H}_{i,q-1}^{-1}(\gamma,\bm\theta)\bm{z}_{i,(q)})
		(1+\theta_{(q-1)} \bm{z}_{i,(q-1)}'\bm{H}_{i,q-2}^{-1}(\gamma,\bm\theta)\bm{z}_{i,(q-1)})},
	\end{align*}
	
	\noindent where we note that $\theta_{(q)}$ and $\theta_{(q-1)}$ can be
	arbitrarily small and the dominant term of the denominator of the last equation
	can be equal to 
	\begin{enumerate}
		\item[(i)] $\theta_{(q)}\theta_{(q-1)} \bm{z}_{i,(q)}'\bm{H}_{i,q-1}^{-1}(\gamma,\bm\theta)\bm{z}_{i,(q)} \bm{z}_{i,(q-1)}'\bm{H}_{i,q-2}^{-1}(\gamma,\bm\theta)\bm{z}_{i,(q-1)}$; 
		\item[(ii)] $\theta_{(q)} \bm{z}_{i,(q)}'\bm{H}_{i,q-1}^{-1}(\gamma,\bm\theta)\bm{z}_{i,(q)}+\theta_{(q-1)} \bm{z}_{i,(q-1)}'\bm{H}_{i,q-2}^{-1}(\gamma,\bm\theta)\bm{z}_{i,(q-1)}$; 
		\item[(iii)] $1$.
	\end{enumerate} 
	\noindent For the case of (i), $\theta_{(q)}n_i^{\ell}\rightarrow\infty$ and
	$\theta_{(q-1)}n_i^{\ell}\rightarrow\infty$ by Lemma \ref{appendix:lemma:z}
	(i); hence, using Lemma \ref{appendix:lemma:z} (i), we have
	\begin{align*}
	\theta_{(q)}\theta_{(q-1)}\bm{z}_{i,(q)}'\bm{H}_{i,q}^{-1}(\gamma,\bm\theta)\bm{z}_{i,(q-1)}
	=&~	o_p(n_i^{-\ell-\tau}),
	\end{align*}	
	\noindent which also gives the following two results:
	\begin{align*}
	\theta_{(q)}\bm{z}_{i,(q)}'\bm{H}_{i,q}^{-1}(\gamma,\bm\theta)\bm{z}_{i,(q-1)}
	=&~	o_p(n_i^{-\tau}),\\
	\bm{z}_{i,(q)}'\bm{H}_{i,q}^{-1}(\gamma,\bm\theta)\bm{z}_{i,(q-1)}
	=&~	o_p(n_i^{\ell-\tau}).
	\end{align*}
	\noindent For the case of (ii), $\theta_{(q)}n_i^{\ell}\rightarrow\infty$ and $\theta_{(q)}=O(n_i^{-\ell})$
	(or vice versa) by Lemma \ref{appendix:lemma:z} (i); hence, using Lemma
	\ref{appendix:lemma:z} (i), we have
	\begin{align*}
	\theta_{(q)}\theta_{(q-1)}\bm{z}_{i,(q)}'\bm{H}_{i,q}^{-1}(\gamma,\bm\theta)\bm{z}_{i,(q-1)}
	=&~	o_p(\theta_{(q-1)}n_i^{-\tau}),
	\end{align*}
	\noindent which gives the following three results:
	\begin{align*}
	\theta_{(q)}\theta_{(q-1)}\bm{z}_{i,(q)}'\bm{H}_{i,q}^{-1}(\gamma,\bm\theta)\bm{z}_{i,(q-1)}
	=&~	o_p(n_i^{-\ell-\tau}),\\
	\theta_{(q)}\bm{z}_{i,(q)}'\bm{H}_{i,q}^{-1}(\gamma,\bm\theta)\bm{z}_{i,(q-1)}
	=&~	o_p(n_i^{-\tau}),\\
	\bm{z}_{i,(q)}'\bm{H}_{i,q}^{-1}(\gamma,\bm\theta)\bm{z}_{i,(q-1)}
	=&~	o_p(n_i^{\ell-\tau}).
	\end{align*}
	\noindent For the case of (iii), $\theta_{(q)}=O(n_i^{-\ell})$ and
	$\theta_{(q)}=O(n_i^{-\ell})$ by Lemma \ref{appendix:lemma:z} (i); hence,
	using Lemma \ref{appendix:lemma:z} (i), we have
	\begin{align*}
	\theta_{(q)}\theta_{(q-1)}\bm{z}_{i,(q)}'\bm{H}_{i,q}^{-1}(\gamma,\bm\theta)\bm{z}_{i,(q-1)}
	=&~	o_p(\theta_{(q)}\theta_{(q-1)}n_i^{\ell-\tau}),
	\end{align*}
	\noindent which also gives the following three results:
	\begin{align*}
	\theta_{(q)}\theta_{(q-1)}\bm{z}_{i,(q)}'\bm{H}_{i,q}^{-1}(\gamma,\bm\theta)\bm{z}_{i,(q-1)}
	=&~	o_p(n_i^{-\ell-\tau}),\\
	\theta_{(q)}\bm{z}_{i,(q)}'\bm{H}_{i,q}^{-1}(\gamma,\bm\theta)\bm{z}_{i,(q-1)}
	=&~	o_p(n_i^{-\tau}),\\
	\bm{z}_{i,(q)}'\bm{H}_{i,q}^{-1}(\gamma,\bm\theta)\bm{z}_{i,(q-1)}
	=&~	o_p(n_i^{\ell-\tau}).
	\end{align*}
	\noindent In conclusion, we have
	\begin{align}
	\begin{split}
	\theta_{(q)}\theta_{(q-1)}\bm{z}_{i,(q)}'\bm{H}_{i,q}^{-1}(\gamma,\bm\theta)\bm{z}_{i,(q-1)}
	=&~	o_p(n_i^{-\ell-\tau}),\\
	\theta_{(q)}\bm{z}_{i,(q)}'\bm{H}_{i,q}^{-1}(\gamma,\bm\theta)\bm{z}_{i,(q-1)}
	=&~	o_p(n_i^{-\tau}),\\
	\bm{z}_{i,(q)}'\bm{H}_{i,q}^{-1}(\gamma,\bm\theta)\bm{z}_{i,(q-1)}
	=&~	o_p(n_i^{\ell-\tau})
	\end{split}
	\label{proof:revised two}
	\end{align}
	\noindent uniformly over $\bm\theta\in[0,\infty)^q$. This completes the
	proof of Lemma \ref{appendix:lemma:z} (iii).
	
	We finally prove Lemma \ref{appendix:lemma:z} (iv). Without loss of
	generality, it suffices to prove Lemma \ref{appendix:lemma:z} (iv) by
	replacing $\bm{H}_i(\gamma,\bm\theta)$ with $\bm{H}_{i,q-1}(\gamma,\bm\theta)$
	with $q(\gamma)=q$, $k=(q-1)$, and $k^*=(q)$. Then by (\ref{eq:inv decompose}),
	\begin{align*}
	\theta_{(q-1)}&\bm{z}_{i,(q-1)}'\bm{H}_{i,q-1}^{-1}(\gamma,\bm\theta)\bm{z}_{i,(q)}\\
	=&~ \theta_{(q-1)}\bigg\{\bm{z}_{i,(q-1)}'\bm{H}_{i,q-2}^{-1}(\gamma,\bm\theta)\bm{z}_{i,(q)}\\
	&~  -\frac{\theta_{(q-1)}\bm{z}_{i,(q-1)}'\bm{H}_{i,q-2}^{-1}(\gamma,\bm\theta)\bm{z}_{i,(q-1)}
		\bm{z}_{i,(q-1)}'\bm{H}_{i,q-2}^{-1}(\gamma,\bm\theta)\bm{z}_{i,(q)}}
	{1+\theta_{(q-1)} \bm{z}_{i,(q-1)}'\bm{H}_{i,q-2}^{-1}(\gamma,\bm\theta)\bm{z}_{i,(q-1)}}\bigg\}\\
	=&~ \frac{\theta_{(q-1)}\bm{z}_{i,(q-1)}'\bm{H}_{i,q-2}^{-1}(\gamma,\bm\theta)\bm{z}_{i,(q)}}
	{1+\theta_{(q-1)} \bm{z}_{i,(q-1)}'\bm{H}_{i,q-2}^{-1}(\gamma,\bm\theta)\bm{z}_{i,(q-1)}}.
	\end{align*}
	
	\noindent Hence, Lemma \ref{appendix:lemma:z} (iv) follows from Lemma
	\ref{appendix:lemma:z} (i) and arguments similar to the proof of
	\eqref{proof:revised one}. This completes the proof.
	
	\subsection{Proof of Lemma \ref{appendix:lemma:epsilon}}
	
	Note that for $k=1,\ldots,q$ and $j=1,\ldots,p$,
	\begin{align*}
	\bm\epsilon_i'\bm{z}_{i,k}
	=&~ O_p(n_i^{\ell/2}),\\
	\bm\epsilon_i'\bm{x}_{i,j}
	=&~ O_p(n_i^{\xi/2}).
	\end{align*}
	\noindent Lemma \ref{appendix:lemma:epsilon} (ii)--(iii) then follow
	arguments similarly from the induction and the proofs of
	Lemma \ref{appendix:lemma:z x} (iii) are hence omitted.
	
	We next prove Lemma \ref{appendix:lemma:epsilon} (iv). Let $\bm{z}_{i,(s)}$
	be the $s$-th column of $\bm{Z}_i(\gamma)$ and
	$\bm{H}_{i,t}(\gamma,\bm\theta)$ be defined in (\ref{proof:lem:z:eq0}).
	Without loss of
	generality, we assume $q(\gamma)=q$. Hence by (\ref{fn:expand inv H}), Lemma
	\ref{appendix:lemma:z} (i), and Lemma~\ref{appendix:lemma:epsilon}~(ii),
	we have \begin{align*}
	\bm\epsilon_i'\bm{H}_{i,q}^{-1}(\gamma,\bm\theta)\bm\epsilon_i
	=&~ \bm\epsilon_i'\bm\epsilon_i
	-\sum_{k=1}^{q}\frac{\theta_{(k)} \bm\epsilon_i'\bm{H}_{i,k-1}^{-1}(\gamma,\bm\theta)\bm{z}_{i,(k)}\bm{z}_{i,(k)}'\bm{H}_{i,k-1}^{-1}(\gamma,\bm\theta)\bm\epsilon_i}
	{1+\theta_{(k)}\bm{z}_{i,(k)}'\bm{H}_{i,k-1}^{-1}(\gamma,\bm\theta)\bm{z}_{i,(k)}}\\
	=&~ \bm\epsilon_i'\bm\epsilon_i+ O_p(q)
	\end{align*}
	
	\noindent uniformly over $\bm\theta\in[0,\infty)^{q}$. This completes the
	proof of Lemma \ref{appendix:lemma:epsilon} (iv).
	
	It remains to prove Lemma \ref{appendix:lemma:epsilon} (i). Again, without
	loss of generality, it suffices to prove
	Lemma \ref{appendix:lemma:epsilon} (i) for $q(\gamma)=q$ and $k=(q)$.
	Then by (\ref{eq:inv decompose}),
	\begin{align*}
	\theta_{(q)}\bm\epsilon_i'\bm{H}_{i,q}^{-1}(\gamma,\bm\theta)\bm{z}_{i,(q)}
	=&~ \frac{\theta_{(q)}\bm\epsilon_i'\bm{H}_{i,q-1}^{-1}(\gamma,\bm\theta)\bm{z}_{i,(q)}}
	{1+\theta_{(q)}\bm{z}_{i,(q)}'\bm{H}_{i,q-1}^{-1}(\gamma,\bm\theta)\bm{z}_{i,(q)}}.
	\end{align*}
	\noindent Hence, Lemma \ref{appendix:lemma:epsilon} (i) follows from  Lemma
	\ref{appendix:lemma:z} (i), Lemma \ref{appendix:lemma:epsilon} (ii), and
	arguments similar to the proof of \eqref{proof:revised one}. This completes
	the proof.
	
	\subsection{Proof of Lemma \ref{appendix:prop:compact space}}
	We show the lemma for
	$(\alpha,\gamma)\in\mathcal{A}_0\times\mathcal{G}_0$, where the
	proofs with respect to the remaining models are similar and are hence
	omitted.
	
	Let $\bm{z}_{i,(s)}$ be the $s$-th column of
	$\bm{Z}_i(\gamma)$ and $\bm{H}_{i,t}(\gamma,\bm\theta)$ be defined in
	(\ref{proof:lem:z:eq0}). Without loss of generality, we assume that
	$q(\gamma)=q$ and $\bm{Z}_i(\gamma_0)\bm{b}_i(\gamma_0) =
	\sum_{s=q-q_0+1}^{q}\bm{z}_{i,(s)}b_{i,(s)}$. It then suffices to prove that for
	$(\alpha,\gamma)\in\mathcal{A}_0\times\mathcal{G}_0$ and $v^2>0$
	\begin{align}
	-2\log L(\bm\theta,v^2;\alpha,\gamma)
	- \{-2\log L(\bm\theta_0^\dag,v^2;\alpha,\gamma)\}\xrightarrow{p}\infty,
	\label{proof:prop:compact space:eq15}
	\end{align}
	\noindent as both $N\rightarrow\infty$ and $\theta_{(k)}\rightarrow0$ for some $k\in\{q-q_0+1,\dots,q\}$,
	where $\bm\theta_0^\dag\equiv(0,\dots,0,\theta_{(q-q_0+1),0},\dots,\theta_{(q),0})'$,
	and $\theta_{(s),0}$ being the true
	value of $\theta_{(s)}$; $s=q-q_0+1,\dots,q$.
	Note that by
	(\ref{matrix:H:lmm}) and (\ref{eq:det decompose}), we have
	\begin{align*}
	\det(\bm{H}_i(\gamma,\bm\theta))
	=&~ \det\bigg(\bm{I}_{n_i}+\sum_{s=1}^{q}\theta_{(s)}\bm{z}_{i,(s)}\bm{z}_{i,(s)}'\bigg)\\
	=&~ \det(\bm{H}_{i,q-1}(\gamma,\bm\theta)+\theta_{(q)}\bm{z}_{i,(q)}\bm{z}_{i,(q)}')\\
	=&~ \det(\bm{H}_{i,q-1}(\gamma,\bm\theta))(1+\theta_{(q)}\bm{z}_{i,(q)}'\bm{H}_{i,q-1}^{-1}(\gamma,\bm\theta)\bm{z}_{i,(q)}).
	\end{align*}
	\noindent Continuously expanding the above equation by (\ref{eq:det
		decompose}) yields
	\begin{align*}
	\begin{split}
	\log\det(\bm{H}_i(\gamma,\bm\theta))
	=&~ \log\bigg\{\prod_{s=1}^q
	(1+\theta_{(s)}\bm{z}_{i,(s)}'\bm{H}_{i,s-1}^{-1}(\gamma,\bm\theta)\bm{z}_{i,(s)})\bigg\}\\
	=&~ \sum_{s=1}^q\log(1+\theta_{(s)}\bm{z}_{i,(s)}'\bm{H}_{i,s-1}^{-1}(\gamma,\bm\theta)\bm{z}_{i,(s)}),
	\end{split}
	\end{align*}
	
	\noindent where $\bm{H}_{i,0}(\gamma,\bm\theta)=\bm{I}_{n_i}$.
	This together with (\ref{fn:likelihood}) and (\ref{proof:prop:compact
	space:eq2}) yields
	for $(\alpha,\gamma)\in\mathcal{A}_0\times\mathcal{G}_0$ and fixed
	$v^2>0$,
	\begin{align*}
	\begin{split}
	-2&\log L(\bm\theta,v^2;\alpha,\gamma)\\
	=&~ N\log(2\pi) + N\log(v^2) + \log\det(\bm{H}(\gamma,\bm\theta))
	+ \frac{\bm{y}'\bm{H}^{-1}(\gamma,\bm\theta)\bm{A}(\alpha,\gamma;\bm\theta)\bm{y}}{v^2}\\
	=&~ N\log(2\pi) + N\log(v^2) + \sum_{i=1}^m\sum_{s=1}^q
	\log(1+\theta_{(s)}\bm{z}_{i,(s)}'\bm{H}_{i,s-1}^{-1}(\gamma,\bm\theta)\bm{z}_{i,(s)})\\
	&~ + \frac{(\bm{Z}(\gamma_0)\bm{b}(\gamma_0)+\bm\epsilon)'\bm{H}^{-1}(\gamma,\bm\theta)
		(\bm{I}_{N}-\bm{M}(\alpha,\gamma;\bm\theta))(\bm{Z}(\gamma_0)\bm{b}(\gamma_0)+\bm\epsilon)}{v^2}.
	\end{split}
	\end{align*}
	
	\noindent Hence, we have, for
	$(\alpha,\gamma)\in\mathcal{A}_0\times\mathcal{G}_0$,
	\begin{align*}
	-2\log& L(\bm\theta,v^2;\alpha,\gamma)
	- \{-2\log L(\bm\theta_0^\dag,v^2;\alpha,\gamma)\}\\
	=&~ \sum_{i=1}^m\bigg\{ \sum_{s=q-q_0+1}^{q}\log\bigg(\frac{1+\theta_{(s)}\bm{z}_{i,(s)}'\bm{H}_{i,s-1}^{-1}(\gamma,\bm\theta)\bm{z}_{i,(s)}}
	{1+\theta_{(s),0}\bm{z}_{i,(s)}'\bm{H}_{i,s-1}^{-1}(\gamma,\bm\theta_0^\dag)\bm{z}_{i,(s)}}\bigg)\bigg\}\\
	&~  +\frac{1}{v^2}(\bm{Z}(\gamma_0)\bm{b}(\gamma_0)+\bm\epsilon)'
	\big\{\bm{H}^{-1}(\gamma,\bm\theta)(\bm{I}_{N}-\bm{M}(\alpha,\gamma;\bm\theta))\\
	&~    -\bm{H}^{-1}(\gamma,\bm\theta_0^\dag)(\bm{I}_{N}-\bm{M}(\alpha,\gamma;\bm\theta_0^\dag))
	\big\}(\bm{Z}(\gamma_0)\bm{b}(\gamma_0)+\bm\epsilon),
	\end{align*}
	\noindent where
	\begin{align*}
	(\bm{Z}(\gamma_0)\bm{b}&(\gamma_0)+\bm\epsilon)'
	\big\{\bm{H}^{-1}(\gamma,\bm\theta)(\bm{I}_{N}-\bm{M}(\alpha,\gamma;\bm\theta))\\
	&~    -\bm{H}^{-1}(\gamma,\bm\theta_0^\dag)(\bm{I}_{N}-\bm{M}(\alpha,\gamma;\bm\theta_0^\dag))
	\big\}(\bm{Z}(\gamma_0)\bm{b}(\gamma_0)+\bm\epsilon)\\
	=&~ \bm{b}(\gamma_0)'\bm{Z}(\gamma_0)'\{\bm{H}^{-1}(\gamma,\bm\theta)
	(\bm{I}_{N}-\bm{M}(\alpha,\gamma;\bm\theta))\\
	&~  - \bm{H}^{-1}(\gamma,\bm\theta_0^\dag)
	(\bm{I}_{N}-\bm{M}(\alpha,\gamma;\bm\theta_0^\dag))\}\bm{Z}(\gamma_0)\bm{b}(\gamma_0)\\
	&~    +2\bm{b}(\gamma_0)'\bm{Z}(\gamma_0)'\{\bm{H}^{-1}(\gamma,\bm\theta)
	(\bm{I}_{N}-\bm{M}(\alpha,\gamma;\bm\theta))\\
	&~  - \bm{H}^{-1}(\gamma,\bm\theta_0^\dag)
	(\bm{I}_{N}-\bm{M}(\alpha,\gamma;\bm\theta_0^\dag))\}\bm\epsilon\\
	&~  +\bm\epsilon'\{\bm{H}^{-1}(\gamma,\bm\theta)
	(\bm{I}_{N}-\bm{M}(\alpha,\gamma;\bm\theta))\\
	&~  -\bm{H}^{-1}(\gamma,\bm\theta_0^\dag)
	(\bm{I}_{N}-\bm{M}(\alpha,\gamma;\bm\theta_0^\dag))\}\bm\epsilon.
	\end{align*}
	\noindent Hence, for (\ref{proof:prop:compact space:eq15}) to hold,
	it suffices to prove
	\begin{align}
	\begin{split}
	\bm\epsilon'&\{\bm{H}^{-1}(\gamma,\bm\theta)
	(\bm{I}_{N}-\bm{M}(\alpha,\gamma;\bm\theta))\\
	&~- \bm{H}^{-1}(\gamma,\bm\theta_0^\dag)
	(\bm{I}_{N}-\bm{M}(\alpha,\gamma;\bm\theta_0^\dag))\}\bm\epsilon
	=O_p(m)
	\end{split}
	\label{proof:prop:compact space:eq7}
	\end{align}
	
	\noindent uniformly over $\bm{\theta}\in[0,\infty)^{q}$ and
	\begin{align}
	\begin{split}
	\sum_{i=1}^m&~\bigg\{ \sum_{s=q-q_0+1}^{q}\log\bigg(\frac{1+\theta_{(s)}\bm{z}_{i,(s)}'\bm{H}_{i,s-1}^{-1}(\gamma,\bm\theta)\bm{z}_{i,(s)}}
	{1+\theta_{(s),0}\bm{z}_{i,(s)}'\bm{H}_{i,s-1}^{-1}(\gamma,\bm\theta_0^\dag)\bm{z}_{i,(s)}}\bigg)\bigg\}\\
	&~  +\frac{1}{v^2}\bigg(\bm{b}(\gamma_0)'\bm{Z}(\gamma_0)'\{\bm{H}^{-1}(\gamma,\bm\theta)
	(\bm{I}_{N}-\bm{M}(\alpha,\gamma;\bm\theta))\\
	&~  - \bm{H}^{-1}(\gamma,\bm\theta_0^\dag)(\bm{I}_{N}-\bm{M}(\alpha,\gamma;\bm\theta_0^\dag))\}\bm{Z}(\gamma_0)\bm{b}(\gamma_0)
	\bigg) + O_p(m)\xrightarrow{p}\infty,
	\end{split}
	\label{proof:prop:compact space:eq6}
	\end{align}
	
	\noindent as both $N\rightarrow\infty$ and $\theta_{(k)}\rightarrow0$ for
	some $k\in\{q-q_0+1,\dots,q\}$. Before proving
	(\ref{proof:prop:compact
		space:eq7}) and (\ref{proof:prop:compact space:eq6}), we prove the
		following equations, for $\bm{h}_{i,k}$ being defined in
		(\ref{matrix:H}) and	$k=q-q_0+1,\dots,q$:
	\begin{align}
	\bm\epsilon'\bm{H}^{-1}(\gamma,\bm\theta_0^\dag)\bm{h}_{i,(k)}\bm{h}_{i,(k)}'\bm{H}^{-1}(\gamma,\bm\theta)\bm\epsilon
	=&~ O_p(1),
	\label{proof:prop:compact space:eq8}\\
	\bm\epsilon'\bm{H}^{-1}(\gamma,\bm\theta)\bm{h}_{i,(k)}\bm{h}_{i,(k)}'\bm{H}^{-1}(\gamma,\bm\theta_0^\dag)\bm{M}(\alpha,\gamma;\bm\theta)\bm\epsilon
	=&~ o_p(1),
	\label{proof:prop:compact space:eq9}
	\end{align}
	\noindent and
	\begin{align}
	\begin{split}
	\bm\epsilon'&\bm{H}^{-1}(\gamma,\bm\theta_0^\dag)\bm{X}(\alpha)
	(\bm{X}(\alpha)'\bm{H}^{-1}(\gamma,\bm\theta)\bm{X}(\alpha))^{-1}\\
	&~\times\bm{X}(\alpha)'\bm{H}^{-1}(\gamma,\bm\theta_0^\dag)\bm{h}_{i,(k)}
	\bm{h}_{i,(k)}'\bm{H}^{-1}(\gamma,\bm\theta)\bm\epsilon
	=o_p(1),
	\end{split}
	\label{proof:prop:compact space:eq10}\\
	\begin{split}
	\bm\epsilon'&\bm{H}^{-1}(\gamma,\bm\theta_0^\dag)\bm{X}(\alpha)(\bm{X}(\alpha)'\bm{H}^{-1}(\gamma,\bm\theta_0^\dag)\bm{X}(\alpha))^{-1}\\
	&~  \times\bm{X}(\alpha)'\bm{H}^{-1}(\gamma,\bm\theta)\bm{h}_{i,(k)}
	\bm{h}_{i,(k)}'\bm{H}^{-1}(\gamma,\bm\theta_0^\dag)\bm{X}(\alpha)\\
	&~  \times(\bm{X}(\alpha)'\bm{H}^{-1}(\gamma,\bm\theta)\bm{X}(\alpha))^{-1}\bm{X}(\alpha)'\bm{H}^{-1}(\gamma,\bm\theta_0^\dag)\bm\epsilon=o_p(1)
	\end{split}
	\label{proof:prop:compact space:eq11}
	\end{align}
	
	\noindent uniformly over $\bm\theta\in [0,\infty)^q$. It
	suffices to prove (\ref{proof:prop:compact
		space:eq8})--(\ref{proof:prop:compact space:eq11}) for $k=q$. For
	(\ref{proof:prop:compact space:eq8}) with $k=q$, we have
	\begin{align*}
	\bm\epsilon'&\bm{H}^{-1}(\gamma,\bm\theta_0^\dag)\bm{h}_{i,(q)}\bm{h}_{i,(q)}'\bm{H}^{-1}(\gamma,\bm\theta)\bm\epsilon\\
	=&~ \{\bm\epsilon_i'\bm{H}_{i,q}^{-1}(\gamma,\bm\theta_0^\dag)\bm{z}_{i,(q)}\}\{\bm{z}_{i,(q)}'\bm{H}_{i,q}^{-1}(\gamma,\bm\theta)\bm\epsilon_i\}\\
	=&~ \{\bm\epsilon_i'\bm{H}_{i,q}^{-1}(\gamma,\bm\theta_0^\dag)\bm{z}_{i,(q)}\}\bigg(
	\frac{\bm{z}_{i,(q)}'\bm{H}_{i,q-1}^{-1}(\gamma,\bm\theta)\bm\epsilon_i}
	{1+\theta_{(q)}\bm{z}_{i,(q)}'\bm{H}_{i,q-1}^{-1}(\gamma,\bm\theta)\bm{z}_{i,(q)}}\bigg)\\
	=&~ \{O_p(n_i^{-\ell/2})\}_{1\times 1}\{O_p(n_i^{\ell/2})\}_{1\times 1}\\
	=&~ O_p(1)
	\end{align*}
	
	\noindent uniformly over $\bm\theta\in[0,\infty)^q$, where
	the second last equality follows from Lemma \ref{appendix:lemma:z} (i) and
	Lemma~\ref{appendix:lemma:epsilon}~(i)--(ii).
	For (\ref{proof:prop:compact space:eq9}) with $k=q$, we have
	\begin{align*}
	\bm\epsilon'&\bm{H}^{-1}(\gamma,\bm\theta)\bm{h}_{i,(q)}\bm{h}_{i,(q)}'\bm{H}^{-1}(\gamma,\bm\theta_0^\dag)\bm{M}(\alpha,\gamma;\bm\theta)\bm\epsilon\\
	=&~ \{\bm\epsilon_i'\bm{H}_{i,q}^{-1}(\gamma,\bm\theta)\bm{z}_{i,(q)}\}
	\bm{h}_{i,(q)}'\bm{H}^{-1}(\gamma,\bm\theta_0^\dag)\bm{M}(\alpha,\gamma;\bm\theta)\bm\epsilon\\
	=&~ \bigg(\frac{\bm\epsilon_i'\bm{H}_{i,q-1}^{-1}(\gamma,\bm\theta)\bm{z}_{i,(q)}}
	{1+\theta_{(q)}\bm{z}_{i,(q)}'\bm{H}_{i,q-1}^{-1}(\gamma,\bm\theta)\bm{z}_{i,(q)}}\bigg)
	\bigg(\frac{\bm{z}_{i,(q)}'\bm{H}_{i,q}^{-1}(\gamma,\bm\theta_0^\dag)\bm{X}_i(\alpha)}{(\sum_{i=1}^m n_i^{\xi})^{1/2}}\bigg)\\
	&~  \times\bigg(\frac{\sum_{i=1}^m\bm{X}_i(\alpha)'\bm{H}_{i,q}^{-1}(\gamma,\bm\theta)\bm{X}_i(\alpha)}
	{\sum_{i=1}^m n_i^{\xi}}\bigg)^{-1}
	\bigg(\frac{\sum_{i=1}^m\bm{X}_i(\alpha)'\bm{H}_{i,q}^{-1}(\gamma,\bm\theta)\bm\epsilon_i}
	{(\sum_{i=1}^m n_i^{\xi})^{1/2}}\bigg)\\
	=&~ \{O_p(n_i^{\ell/2})\}_{1\times1}\{o(n_i^{-\ell/2-\tau})\}_{1\times p(\alpha)}
	\{\bm{T}^{-1}(\alpha)+ \{o(n_{\min}^{-\tau})\}_{p(\alpha)\times p(\alpha)}\}\\
	&~	\times    \{O_p(1)\}_{p(\alpha)\times 1}\\
	=&~ o_p(1)
	\end{align*}
	
	\noindent uniformly over $\bm\theta\in[0,\infty)^q$, where
	the second equality follows from (\ref{fn:M}) and (\ref{matrix:H first
	step}) and
	the third equality follows from (\ref{proof:lemma:xze:eq0}), Lemma
	\ref{appendix:lemma:z x} (iii), and
	Lemma \ref{appendix:lemma:epsilon} (ii)--(iii).
	For (\ref{proof:prop:compact space:eq10}) with $k=q$, we have
	\begin{align*}
	\bm\epsilon'&\bm{H}^{-1}(\gamma,\bm\theta_0^\dag)\bm{X}(\alpha)
	(\bm{X}(\alpha)'\bm{H}^{-1}(\gamma,\bm\theta)\bm{X}(\alpha))^{-1}\\
	&~  \times\bm{X}(\alpha)'\bm{H}^{-1}(\gamma,\bm\theta_0^\dag)\bm{h}_{i,(q)}
	\bm{h}_{i,(q)}'\bm{H}^{-1}(\gamma,\bm\theta)\bm\epsilon\\
	=&~ \bigg(\frac{\sum_{i=1}^m\bm\epsilon_i'\bm{H}_{i,q}^{-1}(\gamma,\bm\theta_0^\dag)
		\bm{X}_i(\alpha)}{(\sum_{i=1}^mn_i^{\xi})^{1/2}}\bigg)\bigg(\frac{\sum_{i=1}^m\bm{X}_i(\alpha)'\bm{H}_{i,q}^{-1}(\gamma,\bm\theta)\bm{X}_i(\alpha)}{\sum_{i=1}^m n_i^{\xi}}\bigg)^{-1}\\
	&~  \times\bigg(\frac{\bm{X}_i(\alpha)'\bm{H}_{i,q}^{-1}(\gamma,\bm\theta_0^\dag)\bm{z}_{i,(q)}}{(\sum_{i=1}^mn_i^{\xi})^{1/2}}\bigg)
	\bigg(\frac{\bm{z}_{i,(q)}'\bm{H}_{i,q-1}^{-1}(\gamma,\bm\theta)\bm\epsilon_i}{1+\theta_{(q)}\bm{z}_{i,(q)}'\bm{H}_{i,q-1}^{-1}(\gamma,\bm\theta)\bm{z}_{i,(q)}}\bigg)\\
	=&~ \{O_p(1)\}_{1\times p(\alpha)}\{\bm{T}^{-1}(\alpha)+ \{o(n_{\min}^{-\tau})\}_{p(\alpha)\times p(\alpha)}\}
	\{o(n_i^{-\ell/2-\tau})\}_{p(\alpha)\times 1}\\
	&~	\times\{O_p(n_i^{\ell/2})\}_{1\times 1}\\
	=&~ o_p(1)
	\end{align*}
	
	\noindent uniformly over $\bm\theta\in[0,\infty)^q$, where
	the second equality follows from (\ref{proof:lemma:xze:eq0}), Lemma
	\ref{appendix:lemma:z x} (iii), and
	Lemma~\ref{appendix:lemma:epsilon}~(ii)--(iii).
	For (\ref{proof:prop:compact space:eq11}) with $k=q$,
	\begin{align*}
	\bm\epsilon'&\bm{H}^{-1}(\gamma,\bm\theta_0^\dag)\bm{X}(\alpha)(\bm{X}(\alpha)'\bm{H}^{-1}(\gamma,\bm\theta_0^\dag)\bm{X}(\alpha))^{-1}
	\bm{X}(\alpha)'\bm{H}^{-1}(\gamma,\bm\theta)\bm{h}_{i,(q)}\\
	&~  \times\bm{h}_{i,(q)}'\bm{H}^{-1}(\gamma,\bm\theta_0^\dag)\bm{X}(\alpha)
	(\bm{X}(\alpha)'\bm{H}^{-1}(\gamma,\bm\theta)\bm{X}(\alpha))^{-1}\bm{X}(\alpha)'\bm{H}^{-1}(\gamma,\bm\theta_0^\dag)\bm\epsilon\\
	=&~ \bigg(\frac{\sum_{i=1}^m\bm\epsilon_i'\bm{H}_{i,q}^{-1}(\gamma,\bm\theta_0^\dag)\bm{X}_i(\alpha)}{(\sum_{i=1}^mn_i^{\xi})^{1/2}}\bigg)
	\bigg(\frac{\sum_{i=1}^m\bm{X}_i(\alpha)'\bm{H}_{i,q}^{-1}(\gamma,\bm\theta_0^\dag)\bm{X}_i(\alpha)}{\sum_{i=1}^mn_i^{\xi}}\bigg)^{-1}\\
	&~  \times\bigg(\frac{\bm{X}_i(\alpha)'\bm{H}_{i,q-1}^{-1}(\gamma,\bm\theta)\bm{z}_{i,(q)}}
	{(\sum_{i=1}^mn_i^{\xi})^{1/2}(1+\theta_{(q)}\bm{z}_{i,(q)}'\bm{H}_{i,q-1}^{-1}(\gamma,\bm\theta)\bm{z}_{i,(q)})}\bigg)
	\bigg(\frac{\bm{z}_{i,(q)}'\bm{H}_{i,q}^{-1}(\gamma,\bm\theta_0^\dag)\bm{X}_i(\alpha)}{(\sum_{i=1}^mn_i^{\xi})^{1/2}}\bigg)\\
	&~  \times\bigg(\frac{\sum_{i=1}^m\bm{X}_i(\alpha)'\bm{H}_{i,q}^{-1}(\gamma,\bm\theta)\bm{X}_i(\alpha)}{\sum_{i=1}^mn_i^{\xi}}\bigg)^{-1}
	\bigg(\frac{\sum_{i=1}^m\bm{X}_i(\alpha)'\bm{H}_{i,q}^{-1}(\gamma,\bm\theta_0^\dag)\bm\epsilon_i}{(\sum_{i=1}^mn_i^{\xi})^{1/2}}\bigg)\\
	=&~ \{O_p(1)\}_{1\times p(\alpha)}\{\bm{T}^{-1}(\alpha)+ \{o(n_{\min}^{-\tau})\}_{p(\alpha)\times p(\alpha)}\}
	\{o(n_i^{\ell/2-\tau})\}_{p(\alpha)\times 1}\\
	&~  \times \{o(n_{\min}^{-\tau})\}_{1\times p(\alpha)}
	\{\bm{T}^{-1}(\alpha)+ \{o(n_{\min}^{-\tau})\}_{p(\alpha)\times p(\alpha)}\}
	\{O_p(1)\}_{p(\alpha)\times 1}\\
	=&~  o_p(1)
	\end{align*}
	
	\noindent uniformly over $\bm\theta\in[0,\infty)^q$, where
	the second equality follows from (\ref{proof:lemma:xze:eq0}), Lemma
	\ref{appendix:lemma:z x} (ii)--(iii), and
	Lemma~\ref{appendix:lemma:epsilon}~(iii).
	This completes the proofs of (\ref{proof:prop:compact
		space:eq8})--(\ref{proof:prop:compact space:eq11}). We now prove
	(\ref{proof:prop:compact space:eq7}). Note that
	\begin{align}
	\begin{split}
	\bm\epsilon'&\{\bm{H}^{-1}(\gamma,\bm\theta_0^\dag)
	\bm{M}(\alpha,\gamma;\bm\theta_0^\dag)-\bm{H}^{-1}(\gamma,\bm\theta)
	\bm{M}(\alpha,\gamma;\bm\theta)\}\bm\epsilon\\
	=&~\bm\epsilon'\{\bm{H}^{-1}(\gamma,\bm\theta_0^\dag)
	\bm{M}(\alpha,\gamma;\bm\theta_0^\dag)-\bm{H}^{-1}(\gamma,\bm\theta_0^\dag)\bm{M}(\alpha,\gamma;\bm\theta)\\
	&~  +\bm{H}^{-1}(\gamma,\bm\theta_0^\dag)\bm{M}(\alpha,\gamma;\bm\theta)
	-\bm{H}^{-1}(\gamma,\bm\theta)\bm{M}(\alpha,\gamma;\bm\theta)\}\bm\epsilon\\
	=&~ \bm\epsilon'\{\bm{H}^{-1}(\gamma,\bm\theta_0^\dag)\bm{M}(\alpha,\gamma;\bm\theta_0^\dag)
	-\bm{H}^{-1}(\gamma,\bm\theta_0^\dag)\bm{M}(\alpha,\gamma;\bm\theta)\}\bm\epsilon + o_p(m)\\
	=&~ \bm\epsilon'\bm{H}^{-1}(\gamma,\bm\theta_0^\dag)\{\bm{M}(\alpha,\gamma;\bm\theta_0^\dag)\\
	&~  -\bm{X}(\alpha)(\bm{X}(\alpha)'\bm{H}^{-1}(\gamma,\bm\theta)\bm{X}(\alpha))^{-1}\bm{X}(\alpha)'\bm{H}^{-1}(\gamma,\bm\theta_0^\dag)\\
	&~  +\bm{X}(\alpha)(\bm{X}(\alpha)'\bm{H}^{-1}(\gamma,\bm\theta)\bm{X}(\alpha))^{-1}\bm{X}(\alpha)'\bm{H}^{-1}(\gamma,\bm\theta_0^\dag)\\
	&~  -\bm{M}(\alpha,\gamma;\bm\theta)\}\bm\epsilon + o_p(m)\\
	=&~ \bm\epsilon'\bm{H}^{-1}(\gamma,\bm\theta_0^\dag)\{\bm{M}(\alpha,\gamma;\bm\theta_0^\dag)\\
	&~  -\bm{X}(\alpha)(\bm{X}(\alpha)'\bm{H}^{-1}(\gamma,\bm\theta)\bm{X}(\alpha))^{-1}
	\bm{X}(\alpha)'\bm{H}^{-1}(\gamma,\bm\theta_0^\dag)\}\bm\epsilon + o_p(m)\\
	=&~ o_p(m)
	\end{split}
	\label{proof:prop:compact space:eq12}
	\end{align}
	
	\noindent uniformly over $\bm\theta\in[0,\infty)^q$, where the second
	equality follows from
	(\ref{proof:prop:compact space:eq9}) that
	\begin{align*}
	\bm\epsilon'&\{\bm{H}^{-1}(\gamma,\bm\theta_0^\dag)-\bm{H}^{-1}(\gamma,\bm\theta)\}\bm{M}(\alpha,\gamma;\bm\theta)\bm\epsilon\\
	=&~ \bm\epsilon'\bm{H}^{-1}(\gamma,\bm\theta)\{\bm{H}(\gamma,\bm\theta)-\bm{H}(\gamma,\bm\theta_0^\dag)\}
	\bm{H}^{-1}(\gamma,\bm\theta_0^\dag)\bm{M}(\alpha,\gamma;\bm\theta)\bm\epsilon\\
	=&~ \sum_{i=1}^m\sum_{k=q-q_0+1}^{q}(\theta_{(k)}-\theta_{(k),0})\bm\epsilon'\bm{H}^{-1}(\gamma,\bm\theta)\bm{h}_{i,(k)}\bm{h}_{i,(k)}'\bm{H}^{-1}(\gamma,\bm\theta_0^\dag)
	\bm{M}(\alpha,\gamma;\bm\theta)\bm\epsilon\\
	=&~ o_p(m)
	\end{align*}
	
	\noindent uniformly over $\bm\theta\in[0,\infty)^q$, the second last
	equality follows from
	(\ref{proof:prop:compact space:eq10}) that
	\begin{align*}
	\bm\epsilon'&\bm{H}^{-1}(\gamma,\bm\theta_0^\dag)\{
	\bm{X}(\alpha)(\bm{X}(\alpha)'\bm{H}^{-1}(\gamma,\bm\theta)\bm{X}(\alpha))^{-1}\bm{X}(\alpha)'\bm{H}^{-1}(\gamma,\bm\theta_0^\dag)\\
	&~	-\bm{M}(\alpha,\gamma;\bm\theta)
	\}\bm\epsilon\\
	=&~\bm\epsilon'\bm{H}^{-1}(\gamma,\bm\theta_0^\dag)
	\bm{X}(\alpha)(\bm{X}(\alpha)'\bm{H}^{-1}(\gamma,\bm\theta)\bm{X}(\alpha))^{-1}\bm{X}(\alpha)'\\
	&~	\times\{\bm{H}^{-1}(\gamma,\bm\theta_0^\dag)-\bm{H}^{-1}(\gamma,\bm\theta)\}\bm\epsilon\\
	=&~ \sum_{i=1}^m\sum_{k=q-q_0+1}^{q}(\theta_{(k)}-\theta_{(k),0})
	\bm\epsilon'\bm{H}^{-1}(\gamma,\bm\theta_0^\dag)\bm{X}(\alpha)
	(\bm{X}(\alpha)'\bm{H}^{-1}(\gamma,\bm\theta_0^\dag)\bm{X}(\alpha))^{-1}\\
	&~  \times
	\bm{X}(\alpha)'\bm{H}^{-1}(\gamma,\bm\theta)\bm{h}_{i,(k)}
	\bm{h}_{i,(k)}'\bm{H}^{-1}(\gamma,\bm\theta_0^\dag)\bm{X}(\alpha)\\
	=&~ o_p(m)
	\end{align*}
	
	\noindent uniformly over $\bm\theta\in[0,\infty)^q$, and the last equality
	follows from
	(\ref{proof:prop:compact space:eq11}) that
	\begin{align*}
	\bm\epsilon'&~\bm{H}^{-1}(\gamma,\bm\theta_0^\dag)\{\bm{M}(\alpha,\gamma;\bm\theta_0^\dag)\\
	&~  -\bm{X}(\alpha)(\bm{X}(\alpha)'\bm{H}^{-1}(\gamma,\bm\theta)\bm{X}(\alpha))^{-1}\bm{X}(\alpha)'\bm{H}^{-1}(\gamma,\bm\theta_0^\dag)\}\bm\epsilon\\
	=&~ \bm\epsilon'\bm{H}^{-1}(\gamma,\bm\theta_0^\dag)\bm{X}(\alpha)
	\{(\bm{X}(\alpha)'\bm{H}^{-1}(\gamma,\bm\theta_0^\dag)\bm{X}(\alpha))^{-1}\\
	&~    -(\bm{X}(\alpha)'\bm{H}^{-1}(\gamma,\bm\theta)\bm{X}(\alpha))^{-1}\}\bm{X}(\alpha)'\bm{H}^{-1}(\gamma,\bm\theta_0^\dag)\bm\epsilon\\
	=&~ \bm\epsilon'\bm{H}^{-1}(\gamma,\bm\theta_0^\dag)\bm{X}(\alpha)(\bm{X}(\alpha)'\bm{H}^{-1}(\gamma,\bm\theta_0^\dag)\bm{X}(\alpha))^{-1}
	\bm{X}(\alpha)'\{\bm{H}^{-1}(\gamma,\bm\theta)\\
	&~  -\bm{H}^{-1}(\gamma,\bm\theta_0^\dag)\}
	\bm{X}(\alpha)(\bm{X}(\alpha)'\bm{H}^{-1}(\gamma,\bm\theta)\bm{X}(\alpha))^{-1}\bm{X}(\alpha)'\bm{H}^{-1}(\gamma,\bm\theta_0^\dag)\bm\epsilon\\
	=&~ \sum_{i=1}^m\sum_{k=q-q_0+1}^{q}(\theta_{(k),0}-\theta_{(k)})\bm\epsilon'\bm{H}^{-1}(\gamma,\bm\theta_0^\dag)\bm{X}(\alpha)
	(\bm{X}(\alpha)'\bm{H}^{-1}(\gamma,\bm\theta_0^\dag)\bm{X}(\alpha))^{-1}\\
	&~  \times\bm{X}(\alpha)'\bm{H}^{-1}(\gamma,\bm\theta)\bm{h}_{i,(k)}
	\bm{h}_{i,(k)}'\bm{H}^{-1}(\gamma,\bm\theta_0^\dag)\bm{X}(\alpha)
	(\bm{X}(\alpha)'\bm{H}^{-1}(\gamma,\bm\theta)\bm{X}(\alpha))^{-1}\\
	&~  \times\bm{X}(\alpha)'\bm{H}^{-1}(\gamma,\bm\theta_0^\dag)\bm\epsilon\\
	=&~ o_p(m)
	\end{align*}
	
	\noindent uniformly over $\bm\theta\in[0,\infty)^q$. Also, by
	(\ref{proof:prop:compact space:eq8}),
	\begin{align*}
	\bm\epsilon'&\{\bm{H}^{-1}(\gamma,\bm\theta)-\bm{H}^{-1}(\gamma,\bm\theta_0^\dag)\}\bm\epsilon\\
	=&~\bm\epsilon'\bm{H}^{-1}(\gamma,\bm\theta_0^\dag)\{\bm{H}(\gamma,\bm\theta_0^\dag)-\bm{H}(\gamma,\bm\theta)\}\bm{H}^{-1}(\gamma,\bm\theta)\bm\epsilon\\
	=&~ \sum_{i=1}^m\sum_{k=q-q_0+1}^{q}\{\theta_{(k),0}-\theta_{(k)}\}
	\bm\epsilon'\bm{H}^{-1}(\gamma,\bm\theta_0^\dag)\bm{h}_{i,(k)}\bm{h}_{i,(k)}'\bm{H}^{-1}(\gamma,\bm\theta)\bm\epsilon\\
	=&~  O_p(m)
	\end{align*}
	
	\noindent uniformly over $\bm\theta\in[0,\infty)^q$. This together with
	(\ref{proof:prop:compact space:eq12}) gives (\ref{proof:prop:compact
		space:eq7}). We now prove (\ref{proof:prop:compact space:eq6}).
	As with the proof of (\ref{proof:prop:compact space:eq12}), we have
	\begin{align*}
	\bm{b}(\gamma_0)'\bm{Z}(\gamma_0)'&\{\bm{H}^{-1}(\gamma,\bm\theta_0^\dag)
	\bm{M}(\alpha,\gamma;\bm\theta_0^\dag)\\
	&~  -\bm{H}^{-1}(\gamma,\bm\theta)
	\bm{M}(\alpha,\gamma;\bm\theta)\}\bm{Z}(\gamma_0)\bm{b}(\gamma_0)
	= o_p(m)
	\end{align*}
	
	\noindent uniformly over $\bm\theta\in[0,\infty)^{q}$. Hence
	\begin{align*}
	\bm{b}(\gamma_0)'&\bm{Z}(\gamma_0)'\{\bm{H}^{-1}(\gamma,\bm\theta)
	(\bm{I}_{N}-\bm{M}(\alpha,\gamma;\bm\theta))\\
	&~-\bm{H}^{-1}(\gamma,\bm\theta_0^\dag)
	(\bm{I}_{N}-\bm{M}(\alpha,\gamma;\bm\theta_0^\dag))\}\bm{Z}(\gamma_0)\bm{b}(\gamma_0)\\
	=&~ \bm{b}(\gamma_0)'\bm{Z}(\gamma_0)'\{\bm{H}^{-1}(\gamma,\bm\theta)-\bm{H}^{-1}(\gamma,\bm\theta_0^\dag)\}\bm{Z}(\gamma_0)\bm{b}(\gamma_0)+ o_p(m)\\
	=&~ \sum_{i=1}^m\sum_{s=q-q_0+1}^{q}(\theta_{(s),0}-\theta_{(s)})\bm{b}(\gamma_0)'\bm{Z}(\gamma_0)'
	\bm{H}^{-1}(\gamma,\bm\theta)\bm{h}_{i,(s)}\\
	&~  \times\bm{h}_{i,(s)}'\bm{H}^{-1}(\gamma,\bm\theta_0^\dag)\bm{Z}(\gamma_0)\bm{b}(\gamma_0) + o_p(m)
	\end{align*}
	
	\noindent uniformly over $\bm\theta\in[0,\infty)^{q}$. Hence, for
	(\ref{proof:prop:compact space:eq6}) to hold, it suffices to prove that for $k
	=q-q_0+1,\dots,q$ and $i=1,\dots,m$,
	\begin{align}
	\begin{split}
	\log\bigg(&~\frac{1+\theta_{(k)}\bm{z}_{i,(k)}'\bm{H}_{i,k-1}^{-1}(\gamma,\bm\theta)\bm{z}_{i,(k)}}
	{1+\theta_{(k),0}\bm{z}_{i,(k)}'\bm{H}_{i,k-1}^{-1}(\gamma,\bm\theta_0^\dag)\bm{z}_{i,(k)}}\bigg)\\
	=&~
	o_p\bigg(\bm{b}(\gamma_0)'\bm{Z}(\gamma_0)'
	\bm{H}^{-1}(\gamma,\bm\theta)\bm{h}_{i,(k)}
	\bm{h}_{i,(k)}'\bm{H}^{-1}(\gamma,\bm\theta_0^\dag)\bm{Z}(\gamma_0)\bm{b}(\gamma_0)\bigg),
	\end{split}
	\label{proof:prop:compact space:eq13}
	\end{align}
	\noindent as both $N\rightarrow\infty$ and $\theta_{(k)}\rightarrow0$ for
	some $k\in\{q-q_0+1,\ldots,q\}$. It suffices to prove
	(\ref{proof:prop:compact space:eq13}) for $k=q$.
	By Lemma \ref{appendix:lemma:z} (ii)--(iii), we have
	\begin{align*}
	\bm{b}(\gamma_0)'&\bm{Z}(\gamma_0)'
	\bm{H}^{-1}(\gamma,\bm\theta)\bm{h}_{i,(q)}
	\bm{h}_{i,(q)}'\bm{H}^{-1}(\gamma,\bm\theta_0^\dag)\bm{Z}(\gamma_0)\bm{b}(\gamma_0)\\
	=&~ \bigg(\frac{b_{i,(q)}\bm{z}_{i,(q)}'\bm{H}_{i,q-1}^{-1}(\gamma,\bm\theta)\bm{z}_{i,(q)}}
	{\{1+\theta_{(q)}\bm{z}_{i,(q)}'\bm{H}_{i,q-1}^{-1}(\gamma,\bm\theta)\bm{z}_{i,(q)}\}}\bigg)
	\{\bm{z}_{i,(q)}'\bm{H}_{i,q}^{-1}(\gamma,\bm\theta_0^\dag)\bm{Z}_i(\gamma_0)\bm{b}_i(\gamma_0)\}\\
	=&~ \bigg(\frac{b_{i,(q)}\bm{z}_{i,(q)}'\bm{H}_{i,q-1}^{-1}(\gamma,\bm\theta)\bm{z}_{i,(q)}}
	{1+\theta_{(q)}\bm{z}_{i,(q)}'\bm{H}_{i,q-1}^{-1}(\gamma,\bm\theta)\bm{z}_{i,(q)}}\bigg)
	\bigg(\frac{b_{i,(q)}}{\theta_{(q),0}} + o_p(n_i^{-\ell-\tau})\bigg).
	\end{align*}
	\noindent Hence, for (\ref{proof:prop:compact space:eq13}) with $k=q$ to
	hold, it suffices to prove that
	\begin{align*}
	\log\bigg(&\frac{1+\theta_{(q)}\bm{z}_{i,(q)}'\bm{H}_{i,q-1}^{-1}(\gamma,\bm\theta)\bm{z}_{i,(q)}}
	{1+\theta_{(q),0}\bm{z}_{i,(q)}'\bm{H}_{i,q-1}^{-1}(\gamma,\bm\theta_0^\dag)\bm{z}_{i,(q)}}\bigg)
	\bigg(\frac{\bm{z}_{i,(q)}'\bm{H}_{i,q-1}^{-1}(\gamma,\bm\theta)\bm{z}_{i,(q)}}
	{1+\theta_{(q)}\bm{z}_{i,(q)}'\bm{H}_{i,q-1}^{-1}(\gamma,\bm\theta)\bm{z}_{i,(q)}}\bigg)^{-1}\\
	&~	\rightarrow0,
	\end{align*}
	\noindent as both $N\rightarrow\infty$ and $\theta_{(q)}\rightarrow0$,
	which follows from Lemma \ref{appendix:lemma:z} (i) and
	L'Hospital's rule. This completes the proof of
	(\ref{proof:prop:compact space:eq13}). This completes the proof.
	
	\subsection{Proof of Lemma \ref{appendix:lemma:xze}}
	
	We first prove Lemma \ref{appendix:lemma:xze} (i). For
	$i,i^*=1,\ldots,m$, $(\alpha,\gamma)\in\mathcal{A}\times\mathcal{G}$ and $k,k^*\in\gamma$,
	we have
	\begin{align*}
	\begin{split}
	\theta_k\theta_{k^*}&
	\bm{h}_{i,k}'\bm{H}^{-1}(\gamma,\bm\theta)\bm{M}(\alpha,\gamma;\bm\theta)\bm{h}_{i^*,k^*}\\
	=&~
	\big(\theta_k\bm{z}_{i,k}'\bm{H}_i^{-1}(\gamma,\bm\theta)\bm{X}_i(\alpha)\big)
	\bigg(\frac{\sum_{i=1}^m\bm{X}_i(\alpha)'\bm{H}_i^{-1}(\gamma,\bm\theta)\bm{X}_i(\alpha)}{\sum_{i=1}^mn_i^{\xi}}\bigg)^{-1}\\
	&~  \times\bigg(\frac{
		\theta_{k^*}\bm{X}_{i^*}(\alpha)'\bm{H}_{i^*}^{-1}(\gamma,\bm\theta)\bm{z}_{i^*,k^*}}{\sum_{i=1}^mn_i^{\xi}}\bigg)\\
	=&~ \{o(n_i^{(\xi-\ell)/2-\tau})\}_{1\times p(\alpha)}\bigg\{\bm{T}^{-1}(\alpha)
	+ \{o(n_{\min}^{-\tau})\}_{p(\alpha)\times p(\alpha)}\bigg\}\\
	&~  \times\bigg\{o\bigg(\frac{n_{i^*}^{(\xi-\ell)/2-\tau}}{\sum_{i=1}^mn_i^{\xi}}\bigg)\bigg\}_{p(\alpha)\times 1}\\
	=& o\Bigg(\frac{n_i^{(\xi-\ell)/2}n_{i^*}^{(\xi-\ell)/2-\tau}}{\sum_{i=1}^mn_i^{\xi}}\Bigg)
	\end{split}
	\end{align*}
	
	\noindent uniformly over $\bm\theta\in\Theta_\gamma$, where the second
	equality follows
	from (\ref{proof:lemma:xze:eq0}) and Lemma \ref{appendix:lemma:z x} (iii).
	Similarly, by (\ref{proof:lemma:xze:eq0}) and Lemma \ref{appendix:lemma:z x} (iii), we have
	\begin{align*}
	\begin{split}
	\theta_k\bm{h}_{i,k}'&\bm{H}^{-1}(\gamma,\bm\theta)\bm{M}(\alpha,\gamma;\bm\theta)\bm{h}_{i^*,k^*}\\
	=&~ \big(\theta_k\bm{z}_{i,k}'\bm{H}_i^{-1}(\gamma,\bm\theta)\bm{X}_i(\alpha)\big)
	\bigg(\frac{\sum_{i=1}^m\bm{X}_i(\alpha)'\bm{H}_i^{-1}(\gamma,\bm\theta)\bm{X}_i(\alpha)}{\sum_{i=1}^mn_i^{\xi}}\bigg)^{-1}\\
	&~  \times\bigg(\frac{\bm{X}_{i^*}(\alpha)'\bm{H}_{i^*}^{-1}(\gamma,\bm\theta)\bm{z}_{i^*,k^*}}{\sum_{i=1}^mn_i^{\xi}}\bigg)\\
	=&~ \{o(n_i^{(\xi-\ell)/2-\tau})\}_{1\times p(\alpha)}\bigg\{\bm{T}^{-1}(\alpha)
	+\{o(n_{\min}^{-\tau})\}_{p(\alpha)\times p(\alpha)}\bigg\}\\
	&~  \times\Bigg\{o\Bigg(\frac{n_{i^*}^{(\xi+\ell)/2-\tau}}{\sum_{i=1}^mn_i^{\xi}}\Bigg)\Bigg\}_{p(\alpha)\times 1}\\
	=&~	 o\Bigg(\frac{n_i^{(\xi-\ell)/2}n_{i^*}^{(\xi+\ell)/2-\tau}}{\sum_{i=1}^mn_i^{\xi}}\Bigg)
	\end{split}
	\end{align*}
	
	\noindent uniformly over $\bm\theta\in\Theta_\gamma$. Further, by
	(\ref{proof:lemma:xze:eq0}) and Lemma \ref{appendix:lemma:z x} (iii), we have
	\begin{align*}
	\begin{split}
	\bm{h}_{i,k}'&\bm{H}^{-1}(\gamma,\bm\theta)\bm{M}(\alpha,\gamma;\bm\theta)\bm{h}_{i^*,k^*}\\
	=&~ \big(\theta_k\bm{z}_{i,k}'\bm{H}_i^{-1}(\gamma,\bm\theta)\bm{X}_i(\alpha)\big)
	\bigg(\frac{\sum_{i=1}^m\bm{X}_i(\alpha)'\bm{H}_i^{-1}(\gamma,\bm\theta)\bm{X}_i(\alpha)}{\sum_{i=1}^mn_i^{\xi}
	}\bigg)^{-1}\\
	&~  \times\bigg(\frac{\bm{X}_{i^*}(\alpha)'\bm{H}_{i^*}^{-1}(\gamma,\bm\theta)\bm{z}_{i^*,k^*}}{\sum_{i=1}^mn_i^{\xi}
	}\bigg)\\
	=&~ \{o(n_i^{(\xi+\ell)/2-\tau})\}_{1\times p(\alpha)}\bigg\{\bm{T}^{-1}(\alpha)
	+\{o(n_{\min}^{-\tau})\}_{p(\alpha)\times p(\alpha)}\bigg\}\\
	&~  \times\Bigg\{o\Bigg(\frac{n_{i^*}^{(\xi+\ell)/2-\tau}}{\sum_{i=1}^mn_i^{\xi}}\Bigg)\Bigg\}_{p(\alpha)\times 1}\\
	=&~ o\Bigg(\frac{n_i^{(\xi+\ell)/2}n_{i^*}^{(\xi+\ell)/2-\tau}}{\sum_{i=1}^mn_i^{\xi}}\Bigg)
	\end{split}
	\end{align*}
	
	\noindent uniformly over $\bm\theta\in\Theta_\gamma$. This completes the proof of
	Lemma \ref{appendix:lemma:xze} (i).
	
	We now prove Lemma \ref{appendix:lemma:xze} (ii). For
	$i,i^*=1,\ldots,m$,
	$(\alpha,\gamma)\in\mathcal{A}\times\mathcal{G}$,
	$k\in\gamma$ and $k^*\notin\gamma$,
	\begin{align*}
	\theta_k&\bm{h}_{i,k}'\bm{H}^{-1}(\gamma,\bm\theta)\bm{M}(\alpha,\gamma;\bm\theta)\bm{h}_{i^*,k^*}\\
	=&~ \big(\theta_k\bm{z}_{i,k}'\bm{H}_i^{-1}(\gamma,\bm\theta)
	\bm{X}_i(\alpha)\big)
	\bigg(\frac{\sum_{i=1}^m\bm{X}_i(\alpha)'\bm{H}_i^{-1}(\gamma,\bm\theta)\bm{X}_i(\alpha)}{\sum_{i=1}^mn_i^{\xi}}\bigg)^{-1}\\
	&~  \times\bigg(\frac{\bm{X}_{i^*}(\alpha)'\bm{H}_{i^*}^{-1}(\gamma,\bm\theta)\bm{z}_{i^*,k^*}}{\sum_{i=1}^mn_i^{\xi}}\bigg)\\
	=&~ \{o(n_i^{(\xi-\ell)/2-\tau})\}_{1\times p(\alpha)}
	\bigg\{\bm{T}^{-1}(\alpha)+\{o(n_{\min}^{-\tau})\}_{p(\alpha)\times p(\alpha)}\bigg\}\\
	&~  \times\bigg\{o\bigg(\frac{n_i^{(\xi+\ell)/2-\tau}}{\sum_{i=1}^mn_i^{\xi}}\bigg)\bigg\}_{1\times p(\alpha)}\\
	=&~o\Bigg(\frac{n_i^{(\xi-\ell)/2}n_{i^*}^{(\xi+\ell)/2-\tau}}{\sum_{i=1}^mn_i^{\xi}}\Bigg)
	\end{align*}
	
	\noindent uniformly over $\bm\theta\in\Theta_\gamma$, where the second
	equality follows from Lemma \ref{appendix:lemma:z x} (ii)--(iii) and (\ref{proof:lemma:xze:eq0}).
	Similarly, by (\ref{proof:lemma:xze:eq0}) and Lemma \ref{appendix:lemma:z x} (ii)--(iii),
	we have
	\begin{align*}
	\bm{h}_{i,k}'&\bm{H}^{-1}(\gamma,\bm\theta)\bm{M}(\alpha,\gamma;\bm\theta)\bm{h}_{i^*,k^*}\\
	=&~ \big(\theta_k\bm{z}_{i,k}'\bm{H}_i^{-1}(\gamma,\bm\theta)\bm{X}_i(\alpha)\big)
	\bigg(\frac{\sum_{i=1}^m\bm{X}_i(\alpha)'\bm{H}_i^{-1}(\gamma,\bm\theta)\bm{X}_i(\alpha)}{\sum_{i=1}^mn_i^{\xi}}\bigg)^{-1}\\
	&~  \times\bigg(\frac{\bm{X}_{i^*}(\alpha)'\bm{H}_{i^*}^{-1}(\gamma,\bm\theta)\bm{z}_{i^*,k^*}}{\sum_{i=1}^mn_i^{\xi}}\bigg)\\
	=&~ \{o(n_i^{(\xi+\ell)/2-\tau})\}_{1\times p(\alpha)}
	\bigg\{\bm{T}^{-1}(\alpha)+\{o(n_{\min}^{-\tau})\}_{p(\alpha)\times p(\alpha)}\bigg\}\\
	&~	\times     \Bigg\{o\Bigg(\frac{n_{i^*}^{(\xi+\ell)/2-\tau}}
	{\sum_{i=1}^mn_i^{\xi}}\Bigg)\Bigg\}_{p(\alpha)\times1}\\
	=&~ o\Bigg(\frac{n_i^{(\xi+\ell)/2}n_{i^*}^{(\xi+\ell)/2-\tau}}{\sum_{i=1}^mn_i^{\xi}}\Bigg)
	\end{align*}
	
	\noindent uniformly over $\bm\theta\in\Theta_\gamma$. This completes the proof of
	Lemma \ref{appendix:lemma:xze} (ii).
	
	We now prove Lemma \ref{appendix:lemma:xze} (iii). For
	$(\alpha,\gamma)\in\mathcal{A}\times\mathcal{G}$ and $k\in\gamma$,
	\begin{align*}
	\begin{split}
	\theta_k&\bm{h}_{i,k}'\bm{H}^{-1}(\gamma,\bm\theta)\bm{M}(\alpha,\gamma;\bm\theta)\bm\epsilon\\
	=&~ \bigg(\frac{\theta_k\bm{z}_{i,k}'\bm{H}_i^{-1}(\gamma,\bm\theta)\bm{X}_i(\alpha)}{(\sum_{i=1}^mn_i^{\xi})^{1/2}}\bigg)
	\bigg(\frac{\sum_{i=1}^m\bm{X}_i(\alpha)'\bm{H}_i^{-1}(\gamma,\bm\theta)\bm{X}_i(\alpha)}
	{\sum_{i=1}^mn_i^{\xi}}\bigg)^{-1}\\
	&~  \times\bigg(\frac{\sum_{i=1}^m\bm{X}_i(\alpha)'\bm{H}_i^{-1}(\gamma,\bm\theta)\bm\epsilon_i}
	{(\sum_{i=1}^mn_i^{\xi})^{1/2}}\bigg)\\
	=&~ \{o(n_i^{-\ell/2-\tau})\}_{1\times p(\alpha)}\bigg\{\bm{T}^{-1}(\alpha)
	+ \{o(n_{\min}^{-\tau})\}_{p(\alpha)\times p(\alpha)}\bigg\}
	\{O_p(1)\}_{p(\alpha)\times 1}\\
	=&~ o_p(n_i^{-\ell/2})
	\end{split}
	\end{align*}
	
	\noindent uniformly over $\bm\theta\in\Theta_\gamma$, where the second
	equality follows from (\ref{proof:lemma:xze:eq0}), Lemma \ref{appendix:lemma:z x} (iii), and
	Lemma~\ref{appendix:lemma:epsilon}~(iii).
	Similarly, by (\ref{proof:lemma:xze:eq0}), Lemma \ref{appendix:lemma:z x} (iii),
	and Lemma \ref{appendix:lemma:epsilon} (iii), we have
	\begin{align*}
	\begin{split}
	\bm{h}_{i,k}'&\bm{H}^{-1}(\gamma,\bm\theta)\bm{M}(\alpha,\gamma;\bm\theta)\bm\epsilon\\
	=&~ \bigg(\frac{\bm{z}_{i,k}'\bm{H}_i^{-1}(\gamma,\bm\theta)\bm{X}_i(\alpha)}{(\sum_{i=1}^mn_i^{\xi})^{1/2}}\bigg)
	\bigg(\frac{\sum_{i=1}^m\bm{X}_i(\alpha)'\bm{H}_i^{-1}(\gamma,\bm\theta)\bm{X}_i(\alpha)}{\sum_{i=1}^mn_i^{\xi}}\bigg)^{-1}\\
	&~  \times\bigg(\frac{\sum_{i=1}^m\bm{X}_i(\alpha)'\bm{H}_i^{-1}(\gamma,\bm\theta)\bm\epsilon_i}{(\sum_{i=1}^mn_i^{\xi})^{1/2}}\bigg)\\
	=&~ \{o(n_i^{\ell/2-\tau})\}_{1\times p(\alpha)}\bigg\{\bm{T}^{-1}(\alpha)
	+ \{o(n_{\min}^{-\tau})\}_{p(\alpha)\times p(\alpha)}\bigg\}
	\{O_p(1)\}_{p(\alpha)\times 1}\\
	=&~ o_p(n_i^{\ell/2})
	\end{split}
	\end{align*}
	
	\noindent uniformly over $\bm\theta\in\Theta_\gamma$. This
	completes the proof of Lemma \ref{appendix:lemma:xze} (iii).
	
	We now prove Lemma \ref{appendix:lemma:xze} (iv). For
	$(\alpha,\gamma)\in(\mathcal{A}\setminus\mathcal{A}_0)\times\mathcal{G}$,
	$k\in\gamma$,
	\begin{align*}
	\begin{split}
	\theta_k&\bm{h}_{i,k}'\bm{H}^{-1}(\gamma,\bm\theta)\bm{M}(\alpha,\gamma;\bm\theta)\bm{X}(\alpha_0\setminus\alpha)\bm\beta(\alpha_0\setminus\alpha)\\
	=&~ \big(\theta_k\bm{z}_{i,k}'\bm{H}_i^{-1}(\gamma,\bm\theta)\bm{X}_i(\alpha)\big)
	\bigg(\frac{\sum_{i=1}^m\bm{X}_i(\alpha)'\bm{H}_i^{-1}(\gamma,\bm\theta)\bm{X}_i(\alpha)}{\sum_{i=1}^mn_i^{\xi}}\bigg)^{-1}\\
	&~  \times\bigg(\frac{\sum_{i=1}^m\sum_{j\in\gamma_0\setminus\alpha}\bm{X}_i(\alpha)'\bm{H}_i^{-1}(\gamma,\bm\theta)\bm{x}_{i,j}\beta_{j,0}}{\sum_{i=1}^mn_i^{\xi}}\bigg)\\
	=&~ \{o(n_i^{(\xi-\ell)/2-\tau})\}_{1\times p(\alpha)}
	\bigg\{\bm{T}^{-1}(\alpha)+\{o(n_{\min}^{-\tau})\}_{p(\alpha)\times p(\alpha)}\bigg\}
	\bigg\{o\bigg(\frac{\sum_{i=1}^m n_i^{\xi-\tau}}{\sum_{i=1}^mn_i^{\xi}}\bigg)\bigg\}_{p(\alpha)\times 1}\\
	=&~ \{o(n_i^{(\xi-\ell)/2-\tau})\}_{1\times p(\alpha)}
	\bigg\{\bm{T}^{-1}(\alpha)+\{o(n_{\min}^{-\tau})\}_{p(\alpha)\times p(\alpha)}\bigg\}
	\{o(n_{\min}^{-\tau})\}_{p(\alpha)\times1}\\
	=&~ o(n_i^{(\xi-\ell)/2-\tau})
	\end{split}
	\end{align*}
	
	\noindent uniformly over $\bm\theta\in\Theta_\gamma$, where the second
	equality follows from (\ref{proof:lemma:xze:eq0}), Lemma
	\ref{appendix:lemma:z x} (i), and Lemma~\ref{appendix:lemma:z x}~(iii).
	Similarly, by (\ref{proof:lemma:xze:eq0}) and Lemma \ref{appendix:lemma:z x}
	(i) and (iii), we have
	\begin{align*}
	\begin{split}
	\bm{h}_{i,k}'&\bm{H}^{-1}(\gamma,\bm\theta)\bm{M}(\alpha,\gamma;\bm\theta)\bm{X}(\alpha_0\setminus\alpha)\bm\beta(\alpha_0\setminus\alpha)\\
	=&~ \big(\theta_k\bm{z}_{i,k}'\bm{H}_i^{-1}(\gamma,\bm\theta)\bm{X}_i(\alpha)\big)
	\bigg(\frac{\sum_{i=1}^m\bm{X}_i(\alpha)'\bm{H}_i^{-1}(\gamma,\bm\theta)\bm{X}_i(\alpha)}{\sum_{i=1}^mn_i^{\xi}}\bigg)^{-1}\\
	&~  \times\bigg(\frac{\sum_{i=1}^m\sum_{j\in\alpha_0\setminus\alpha}\bm{X}_i(\alpha)'\bm{H}_i^{-1}(\gamma,\bm\theta)\bm{x}_{i,j}\beta_{j,0}}{\sum_{i=1}^mn_i^{\xi}}\bigg)\\
	=&~ \{o(n_i^{(\xi+\ell)/2-\tau})\}_{1\times p(\alpha)}
	\bigg\{\bm{T}^{-1}(\alpha)+\{o(n_{\min}^{-\tau})\}_{p(\alpha)\times p(\alpha)}\bigg\}
	\{o(n_{\min}^{-\tau})\}_{p(\alpha)\times1}\\
	=&~ o_p(n_i^{(\xi+\ell)/2-\tau})
	\end{split}
	\end{align*}
	
	\noindent uniformly over $\bm\theta\in\Theta_\gamma$. This
	completes the proof of Lemma \ref{appendix:lemma:xze} (iv).
	
	We now prove Lemma \ref{appendix:lemma:xze} (v). For
	$(\alpha,\gamma)\in\mathcal{A}\times\mathcal{G}$, we have
	\begin{align*}
	\begin{split}
	\bm\epsilon'&\bm{H}^{-1}(\gamma,\bm\theta)\bm{M}(\alpha,\gamma;\bm\theta)\bm\epsilon\\
	=&~\bigg(\frac{\sum_{i=1}^m\bm\epsilon_i'\bm{H}_i^{-1}(\gamma,\bm\theta)\bm{X}_i(\alpha)}{(\sum_{i=1}^mn_i^{\xi})^{1/2}}\bigg)
	\bigg(\frac{\sum_{i=1}^m\bm{X}_i(\alpha)'\bm{H}_i^{-1}(\gamma,\bm\theta)\bm{X}_i(\alpha)}{\sum_{i=1}^mn_i^{\xi}}\bigg)^{-1}\\
	&~  \times\bigg(\frac{\sum_{i=1}^m\bm{X}_i(\alpha)'\bm{H}_i^{-1}(\gamma,\bm\theta)\bm\epsilon_i}{(\sum_{i=1}^mn_i^{\xi})^{1/2}}\bigg)\\
	=&~ \{O_p(1)\}_{1\times p(\alpha)}\bigg\{\bm{T}^{-1}(\alpha)
	+ \{o(n_{\min}^{-\tau})\}_{p(\alpha)\times p(\alpha)}\bigg\}\{O_p(1)\}_{p(\alpha)\times1}\\
	=&~ O_p(p(\alpha))
	\end{split}
	\end{align*}
	
	\noindent uniformly over $\bm\theta\in\Theta_\gamma$, where the second
	equality follows from (\ref{proof:lemma:xze:eq0}) and Lemma
	\ref{appendix:lemma:epsilon} (iii). This completes the proof of Lemma
	\ref{appendix:lemma:xze} (v).
	
	We now prove Lemma \ref{appendix:lemma:xze} (vi). For
	$(\alpha,\gamma)\in\mathcal{A}\times\mathcal{G}$ and
	$k\notin\gamma$, we have
	\begin{align*}
	\bm{h}_{i,k}'&\bm{H}^{-1}(\gamma,\bm\theta)\bm{M}(\alpha,\gamma;\bm\theta)\bm\epsilon\\
	=&~ \bigg(\frac{\bm{z}_{i,k}'\bm{H}_i^{-1}(\gamma,\bm\theta)\bm{X}_i(\alpha)}{(\sum_{i=1}^mn_i^{\xi})^{1/2}}\bigg)
	\bigg(\frac{\sum_{i=1}^m\bm{X}_i(\alpha)'\bm{H}_i^{-1}(\gamma,\bm\theta)\bm{X}_i(\alpha)}{\sum_{i=1}^mn_i^{\xi}}\bigg)^{-1}\\
	&~  \times\bigg(\frac{\sum_{i=1}^m\bm{X}_i(\alpha)'\bm{H}_i^{-1}(\gamma,\bm\theta)\bm\epsilon_i}{(\sum_{i=1}^mn_i^{\xi})^{1/2}}\bigg)\\
	=&~ \{o(n_i^{\ell/2-\tau})\}_{1\times p(\alpha)}
	\bigg\{\bm{T}^{-1}(\alpha)+ \{o(n_{\min}^{-\tau})\}_{p(\alpha)\times p(\alpha)}\bigg\}
	\{O_p(1)\}_{p(\alpha)\times 1}\\
	=&~ o_p(n_i^{\ell/2})
	\end{align*}
	
	\noindent uniformly over $\bm\theta\in\Theta_\gamma$, where the second
	equality follows from (\ref{proof:lemma:xze:eq0}), Lemma
	\ref{appendix:lemma:z x} (ii), and Lemma~\ref{appendix:lemma:epsilon}~(iii).
	This completes the proof of Lemma \ref{appendix:lemma:xze} (vi).
	
	We now prove Lemma \ref{appendix:lemma:xze} (vii). For
	$(\alpha,\gamma)\in(\mathcal{A}\setminus\mathcal{A}_0)\times\mathcal{G}$, we have
	\begin{align*}
	\begin{split}
	\bm\epsilon'&\bm{H}^{-1}(\gamma,\bm\theta)\bm{M}(\alpha,\gamma;\bm\theta)\bm{X}(\alpha_0\setminus\alpha)\bm\beta(\alpha_0\setminus\alpha)\\
	=&~\bigg(\frac{\sum_{i=1}^m\bm\epsilon_i'\bm{H}_i^{-1}(\gamma,\bm\theta)\bm{X}_i(\alpha)}{(\sum_{i=1}^mn_i^{\xi})^{1/2}}\bigg)
	\bigg(\frac{\sum_{i=1}^m\bm{X}_i(\alpha)'\bm{H}_i^{-1}(\gamma,\bm\theta)\bm{X}_i(\alpha)}{\sum_{i=1}^mn_i^{\xi}}\bigg)^{-1}\\
	&~  \times\bigg(\frac{\sum_{i=1}^m\sum_{j\in\alpha_0\setminus\alpha}\bm{X}_i(\alpha)'\bm{H}_i^{-1}(\gamma,\bm\theta)\bm{x}_{i,j}\beta_{j,0}}{(\sum_{i=1}^mn_i^{\xi})^{1/2}}\bigg)\\
	=&~ \{O_p(1)\}_{1\times p(\alpha)}\bigg\{\bm{T}^{-1}(\alpha)
	+ \{o(n_{\min}^{-\tau})\}_{p(\alpha)\times p(\alpha)}\bigg\}\\
	&~	\times \bigg\{o\bigg(\bigg(\sum_{i=1}^m n_i^{\xi}\bigg)^{1/2}n_{\min}^{-\tau}\bigg)\bigg\}_{p(\alpha)\times1}\\
	=&~ o_p\bigg(\bigg(\sum_{i=1}^m n_i^{\xi}\bigg)^{1/2}\bigg)
	\end{split}
	\end{align*}
	
	\noindent uniformly over $\bm\theta\in\Theta_\gamma$, where the second
	equality follows from (\ref{proof:lemma:xze:eq0}), Lemma
	\ref{appendix:lemma:z x} (i), and Lemma~\ref{appendix:lemma:epsilon}~(iii).
	This completes the proof of Lemma \ref{appendix:lemma:xze} (vii).
	
	We now prove Lemma \ref{appendix:lemma:xze} (viii). For
	$i,i^*=1,\ldots,m$, $(\alpha,\gamma)\in\mathcal{A}\times\mathcal{G}$
	and $k,k^*\notin\gamma$, we have
	\begin{align*}
	\bm{h}_{i,k}'&\bm{H}^{-1}(\gamma,\bm\theta)\bm{M}(\alpha,\gamma;\bm\theta)\bm{h}_{i^*,k^*}\\
	=&~ \big(\bm{z}_{i,k}'\bm{H}_i^{-1}(\gamma,\bm\theta)\bm{X}_i(\alpha)\big)
	\bigg(\frac{\sum_{i=1}^m\bm{X}_i(\alpha)'\bm{H}_i^{-1}(\gamma,\bm\theta)\bm{X}_i(\alpha)}{\sum_{i=1}^mn_i^{\xi}}\bigg)^{-1}\\
	&~  \times\bigg(\frac{\bm{X}_{i^*}(\alpha)'\bm{H}_{i^*}^{-1}(\gamma,\bm\theta)\bm{z}_{i^*,k^*}}
	{\sum_{i=1}^mn_i^{\xi}}\bigg)\\
	=&~ \{o(n_i^{(\xi+\ell)/2-\tau})\}_{1\times p(\alpha)}\bigg\{\bm{T}^{-1}(\alpha)
	+ \{o(n_{\min}^{-\tau})\}_{p(\alpha)\times p(\alpha)}\bigg\}\\
	&~  \times\Bigg\{o_p\Bigg(\frac{n_{i^*}^{(\xi+\ell)/2-\tau}}{\sum_{i=1}^mn_i^{\xi}}\Bigg)\Bigg\}_{p(\alpha)\times 1}\\
	=&~` o_p\Bigg(\frac{n_i^{(\xi+\ell)/2}n_{i^*}^{(\xi+\ell)/2-\tau}}{\sum_{i=1}^mn_i^{\xi}}\Bigg)
	\end{align*}
	
	\noindent uniformly over $\bm\theta\in\Theta_\gamma$, where the second
	equality follows from (\ref{proof:lemma:xze:eq0}) and Lemma
	\ref{appendix:lemma:z x} (ii). This completes the proof of
	Lemma \ref{appendix:lemma:xze} (viii).
	
	We now prove Lemma \ref{appendix:lemma:xze} (ix). For
	$(\alpha,\gamma)\in(\mathcal{A}\setminus\mathcal{A}_0)\times\mathcal{G}$,
	$k\notin\gamma$, we have
	\begin{align*}
	\begin{split}
	\bm{h}_{i,k}'\bm{H}^{-1}&(\gamma,\bm\theta)\bm{M}(\alpha,\gamma;\bm\theta)\bm{X}(\alpha_0\setminus\alpha)\bm\beta(\alpha_0\setminus\alpha)\\
	=&~\big(\bm{z}_{i,k}'\bm{H}_i^{-1}(\gamma,\bm\theta)\bm{X}_i(\alpha)\big)
	\bigg(\frac{\sum_{i=1}^m\bm{X}_i(\alpha)'\bm{H}_i^{-1}(\gamma,\bm\theta)\bm{X}_i(\alpha)}{\sum_{i=1}^mn_i^{\xi}}\bigg)^{-1}\\
	&~  \times\bigg(\frac{\sum_{i=1}^m\sum_{j\in\alpha_0\setminus\alpha}
		\bm{X}_i(\alpha)'\bm{H}_i^{-1}(\gamma,\bm\theta)\bm{x}_{i,j}\beta_{j,0}}{\sum_{i=1}^mn_i^{\xi}}\bigg)\\
	=&~ \{o(n_i^{(\xi+\ell)/2-\tau})\}_{1\times p(\alpha)}\bigg\{\bm{T}^{-1}(\alpha)
	+ \{o(n_{\min}^{-\tau})\}_{p(\alpha)\times p(\alpha)}\bigg\}\\
	&~	\times \{o(n_{\min}^{-\tau})\}_{p(\alpha)\times 1}\\
	=&~ o(n_i^{(\xi+\ell)/2-\tau})
	\end{split}
	\end{align*}
	
	\noindent uniformly over $\bm\theta\in\Theta_\gamma$, where the second
	equality follows from (\ref{proof:lemma:xze:eq0}) and Lemma \ref{appendix:lemma:z x} (i)--(ii). This
	completes the proof of Lemma \ref{appendix:lemma:xze} (ix).
	
	We finally prove Lemma \ref{appendix:lemma:xze} (x). For
	$(\alpha,\gamma)\in(\mathcal{A}\setminus\mathcal{A}_0)\times\mathcal{G}$, we have
	\begin{align*}
	\begin{split}
	\bm\beta&(\alpha_0\setminus\alpha)'\bm{X}(\alpha_0\setminus\alpha)'\bm{H}^{-1}(\gamma,\bm\theta)
	\bm{M}(\alpha,\gamma;\bm\theta)\bm{X}(\alpha_0\setminus\alpha)\bm\beta(\alpha_0\setminus\alpha)\\
	=&~ \bigg(\sum_{i=1}^m\sum_{j\in\alpha_0\setminus\alpha}\beta_{j,0}\bm{x}_{i,j}'\bm{H}_i^{-1}(\gamma,\bm\theta)\bm{X}_i(\alpha)\bigg)
	\bigg(\frac{\sum_{i=1}^m\bm{X}_i(\alpha)'\bm{H}_i^{-1}(\gamma,\bm\theta)\bm{X}_i(\alpha)}{\sum_{i=1}^mn_i^{\xi}}\bigg)^{-1}\\
	&~  \times\bigg(\frac{\sum_{i=1}^m\sum_{j\in\alpha_0\setminus\alpha}
		\bm{X}_i(\alpha)'\bm{H}_i^{-1}(\gamma,\bm\theta)\bm{X}_i(\alpha)\bm{x}_{i,j}}{\sum_{i=1}^mn_i^{\xi}}\bigg)\\
	=&~  \bigg\{o\bigg(\sum_{i=1}^mn_i^{\xi-\tau}\bigg)\bigg\}_{1\times p(\alpha)}
	\bigg\{\bm{T}^{-1}(\alpha)+\{o(n_{\min}^{-\tau})\}_{p(\alpha)\times p(\alpha)}\bigg\}
	\{o(n_{\min}^{-\tau})\}_{p(\alpha)\times1}\\
	=&~ o\bigg(\sum_{i=1}^mn_i^{\xi-\tau}\bigg)
	\end{split}
	\end{align*}
	
	\noindent uniformly over $\bm\theta\in\Theta_\gamma$, where the second
	equality follows from (\ref{proof:lemma:xze:eq0}) and Lemma \ref{appendix:lemma:z x} (i). This completes the proof.
	
\end{appendix}
\end{document}